# METHODS IN INDUSTRIAL BIOTECHNOLOGY FOR CHEMICAL ENGINEERS

# METHODS IN INDUSTRIAL BIOTECHNOLOGY FOR CHEMICAL ENGINEERS

**W. B. Vasantha Kandasamy**
e-mail: **vasanthakandasamy@gmail.com**
web: **http://mat.iitm.ac.in/~wbv**

**Florentin Smarandache**
e-mail: **smarand@unm.edu**

# CONTENTS

Chapter Three

## DETERMINATION OF TEMPERATURE SET POINTS FOR CRUDE OIL 47

Chapter Four

## STUDY OF FLOW RATES IN CHEMICAL PLANTS 79

Chapter Five

## MINMIZATION OF WASTE GAS FLOW IN CHEMICAL INDUSTRIES 89

# PREFACE

Industrial Biotechnology is an interdisciplinary topic to which tools of modern biotechnology are applied for finding proper proportion of raw mix of chemicals, determination of set points, finding the flow rates etc., This study is significant as it results in better economy, quality product and control of pollution. The authors in this book have given only methods of industrial biotechnology mainly to help researchers, students and chemical engineers. Since biotechnology concerns practical and diverse applications including production of new drugs, clearing up pollution etc. we have in this book given methods to control pollution in chemical industries as it has become a great health threat in India. In some cases, the damage due to environmental pollution outweighs the benefits of the product.

This book has six chapters. First chapter gives a brief description of biotechnology. Second chapter deals will proper proportion of mix of raw materials in cement industries to minimize pollution using fuzzy control theory. Chapter three gives the method of determination of temperature set point for crude oil in oil refineries. Chapter four studies the flow rates in chemical industries using fuzzy neutral networks. Chapter five gives the method of minimization of waste gas flow in chemical industries using fuzzy linear programming. The final chapter suggests when in these studies indeterminancy is an attribute or concept involved, the notion of neutrosophic methods can be adopted. The authors feel that the reader should be well versed with fuzzy models like neural networks, fuzzy relational equations, fuzzy control theory, fuzzy linear programming and neutrosophic fuzzy models like NRE together with a knowledge of the technical functioning of chemical industries.

The authors are deeply indebted to Dr. Kandasamy, Kama and Meena for their sustained cooperation.

W.B.VASANTHA KANDASAMY
FLORENTIN SMARANDACHE

**Chapter One**

# INTRODUCTION

In keeping with the definition that “biotechnology is really no more than a name given to a set of techniques and processes”, the authors apply some set of fuzzy techniques to chemical industry problems such as finding the proper proportion of raw mix to control pollution, to study flow rates, to find out the better quality of products. We use fuzzy control theory, fuzzy neural networks, fuzzy relational equations, genetic algorithms to these problems for solutions.

When the solution to the problem can have certain concepts or attributes as indeterminate, the only model that can tackle such a situation is the neutrosophic model. The authors have also used these models in this book to study the use of biotechnology in chemical industries.

The new biotechnology revolution began in the 1970s and early 1980s when scientists learned to precisely alter the genetic constitution of living organisms by processes out with traditional breeding practices. This “genetic engineering” has had a profound impact on almost all areas of traditional biotechnology and further permitted breakthroughs in medicine and agriculture, in particular those that would be impossible by traditional breeding approaches.

There are evidences to show that historically biotechnology was an art rather than a science, exemplified in the manufacture of wines, beers, cheeses etc. It is well comprehended by one and all that biotechnology is highly multi disciplinary, it has its foundations in many fields including biology, microbiology, biochemistry, molecular biology, genetics, chemistry and chemical and process engineering. It is further asserted that biotechnology will be the major technology of the twenty first century.

The newly acquired biological knowledge has already made very important contributions to health and welfare of human kind.

Biotechnology is not by itself a product or range of products; it should be regarded as a range of enabling technologies that will find significant application in many industrial sectors.

Traditional biotechnology has established a huge and expanding world market and in monetary terms, represents a major part of all biotechnology financial profits. 'New' aspects of biotechnology founded in recent advances in molecular biology genetic engineering and fermentation process technology are now increasingly finding wide industrial application.

In many ways, biotechnology is a series of embryonic technologies and will require much skilful control of its development but the potentials are vast and diverse and undoubtedly will play an increasingly important part in many future industrial processes.

It is no doubt an interaction between biology and engineering. The developments of biotechnology are proceeding at a speed similar to that of micro-electronics in the mid 1970s. Although the analogy is tempting any expectations that biotechnology will develop commercially at the same spectacular rate should be tempered with considerable caution. While the potential of new biotechnology cannot be doubted a meaningful commercial realization is now slowly occurring and will accelerate as we approach the end of the century. New biotechnology will have a considerable impact across all industrial uses of the life sciences. In each case the relative merits of competing means of production will influence the economics of a biotechnological route. There is no doubt that biotechnology will undoubtedly have great benefits in the long term in all sectors. The growth in awareness of modern biotechnology parallels the serious worldwide changes in the economic climate arising from the escalation of oil prices since 1973.

Biotechnology has been considered as one important means of restimulating the economy whether on a local, regional national or even global basis using new biotechnological methods and new raw materials. Much of modern biotechnology has been developed and utilized by large companies and corporations.

However many small and medium sized companies are realizing that biotechnology is not a science of the future but provides real benefits to their industry today. In many industries traditional technology can produce compounds causing environmental damage whereas biotechnology methods can offer a green alternative promoting a positive public image and also avoiding new environmental penalties.

Biotechnology is high technology par excellence. Science has defined the world in which we live and biotechnology in particular will become an essential and accepted activity of our culture. Biotechnology offers a great deal of hope for solving many of the problems our world faces!. As stated in the Advisory Committee on Science and Technology Report Developments in Biotechnology, public perception of biotechnology will have a major influence on the rate and direction of developments and there is growing concern about genetically modified products. Associated with genetic manipulation are diverse question of safety, ethics and welfare.

Public debate is essential for new biotechnology to grow up and undoubtedly for the foreseeable future, biotechnology will be under scrutiny. We have only given a description of the biotechnology and the new biotechnology. We have highly restricted ourselves from the technical or scientific analysis of the biotechnologies as even in the countries like USA only less than 10% of the population are scientifically literate, so the authors have only described it non-abstractly and in fact we are not in anyway concerned to debate or comment upon it as we acknowledge the deep and dramatic change the world is facing due to biotechnology and new biotechnology.

For more of these particulars please refer [1, 2, 13, 15, 17].

**Chapter Two**

# BIOTECHNOLOGY IN CHEMICAL INDUSTRIES

The chemical industries have become a great threat in India. For the problems they cause on environmental pollution is much more than the benefit derived by their product. Some of them damage other living organisms like fishes, plants and animals; some cause health hazards to people living around the industries like respiratory ailments, skin problems and damage to nervous systems. So we have chosen to illustrate the minimization of pollution by CKD in cement Industries. Most of these problems can be controlled provided one takes the proper proportion of the mix of raw materials, which would minimize the pollution.

Cement kiln dust (CKD) emits nitrogen, carbon etc., that are pollutants of the atmosphere and the waste dust affects the smooth kiln operation of the cement industry system and it reduces the production of clinker quality. Hence the minimization of waste CKD in kiln is an important one in the cement industry. The control of the waste CKD in a kiln is an uncertainty. Researchers approach this problem by mathematical methods and try to account the waste CKD in a cement kiln. But, most of their methods do not properly yield results about the waste CKD in kiln. Further, the control of the waste CKD in kiln is a major problem for this alone can lead to the minimization of atmospheric pollution by the cement industry. So in this chapter we minimize the waste CKD in kiln and account for the waste CKD in kiln using fuzzy control theory and fuzzy neural networks.

In this chapter fuzzy control theory (FCT) is used to study the cement kiln dust (CKD) problem in cement industries. Using fuzzy control method this chapter tries to minimize the cement kiln dust in cement industries. Cement industries of our country happens to be one of the major contributors of dust. The dust arising in various processing units of a cement plant varies in composition. In 1990 the national average was 9 tons of CKD generated for every 100 tons of clinker production. The control of cement kiln dust is a very important issue, because of the following reasons : 1. CKD emits nitrogen, carbon etc., which are pollutants of the atmosphere, 2. The waste dust affects the smooth kiln operation of the cement industry system and it reduces the production of clinker quality. The following creates mainly this waste dust in three ways in cement industries : (a) Cement kiln dust when not collected in time and returned into the kiln, cause air pollution, (b) Process instability and unscheduled kiln shutdowns and (c) Mixing of raw materials.

The data obtained from Graft R. Kessler [12] is used in this chapter to test the result. After using the data from Kessler [12] this chapter tries to minimize the CKD in cement factory. The minimization of CKD plays a vital role in the control of pollution in the atmosphere.

W.Kreft [21] used the interruption of material cycles method for taking account and further utilization of the waste dust in the cement factory. But this method does not properly account the waste CKD. Kesslar [12] has used volatile analysis to reduce CKD. In the volatile analysis method the alkali ratio is used to indicate the waste amount of CKD in clinker.

Kesslar [12] classifies the raw data under investigation in four ways :

I. Monitor and control of the system
II. Burning zone and fuel combustion improvements
III. CKD reprocessing
IV. Find the mix of raw materials in proper proportion.

The ratio of alkali should be lying between 0.5 to 1.5 in Kiln load material. But in this method the CKD was approximately estimated up to 40%. He has not exactly mentioned the percentage of CKD according to the alkali ratio in an online process. So this method has affected largely the kiln system.

In this chapter, in order to account for the waste CKD, the variables are expressed in terms of membership grades. This chapter considers all the four ways of waste CKD mentioned by Kesslar [12] and converts it into a fuzzy control model. This chapter consists of five sections. In section 1 we describe the cement kiln system and the nature of chemical waste dust which pollutes the atmosphere. In section 2 we adopt the fuzzy control theory to monitor and control the system and give suggestion for the improvement of burning and combustion zone. Section 3 deals with the determination of gas volume set point and temperature set point for CKD reprocessing which is vital for the determination of percentage of net CKD. The amount of waste dust depends largely on the mix of raw materials in proper proportion of raw material mix is shown in section 4. The final section deals with results and conclusion obtained from our study.

## 2.1 Description of waste CKD in cement kiln

The data available from any cement industry is used as the information and also as the knowledge about the problem. This serves as the past experience for our study for adapting the fuzzy control theory in this section. This chapter analysis the data via membership functions of fuzzy control method and minimizes the waste CKD in cement industries. Since the cement industry, emits the cement kiln dusts into the atmosphere, this waste dust pollutes the atmosphere.

This analysis not only estimates the cement kiln dust in cement industries but also gives condition to minimize the waste CKD so that the industry will get maximum profit by minimizing the waste CKD in cement industry.

CKD is particulate matter that is collected from kiln exhaust gases and consist of entrained particles of clinker, raw materials and partially calcined raw materials. The present pollution in environment is generated by CKD along with potential future liabilities of stored dust and this should make CKD reduction a high priority. Here we calculate and minimize the net CKD in kiln system. This chapter tackles the problem of minimizing waste CKD in kiln system in four stages. At the first stage we monitor and control the system. In the second stage we adopt time-to- time improved techniques in burning zone and combustion. At the third stage CKD reprocessing is carried out and in the fourth stage we optimize the mix of raw materials in proper proportion using fuzzy neutral network. The above stage-by-stage process is shown in the following figure 2.1.1. Fuzzy control theory and fuzzy neutral network (FNN) is used in this chapter for the above – described method to minimize the CKD in kiln system.

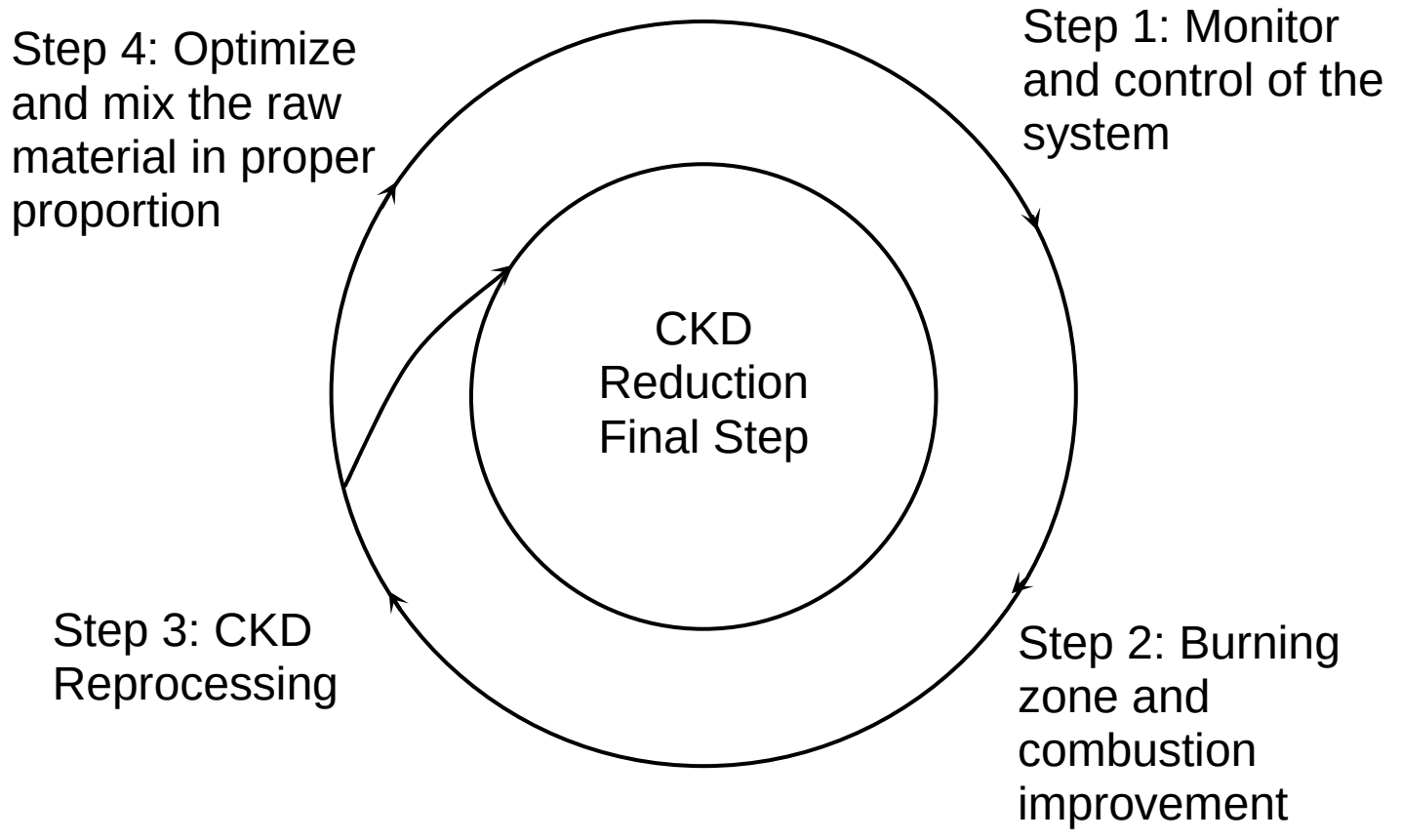

**FIGURE 1: CKD Reduction using fuzzy control**

The fuzzy controller is composed of linguistic control rule, which are conditional linguistic statements of the relationship between inputs and outputs. One of the attractive properties of fuzzy controller is its ability to emulate the behaviour of a human operator. Another important characteristic of a fuzzy controller is its applicability to systems with model uncertainty or even to unknown model systems. The use of fuzzy control applications has expanded at an increasing rate in recent years. In this chapter we use fuzzy control to monitor waste dust in cement kiln system and CKD reprocessing. The fuzzy control in kiln system is described in the figure 2.1.2. We use fuzzy neural network method and tries to find a proper proportion of material mix in cement industries.

The authors aim to achieve a desired level of lime saturation factor (LSF), silica modulus (SM) and alumina modulus (AM) of the raw mix, to produce a particular quality of the cement by controlling the mix proportions of the raw materials. To achieve an appropriate raw mix proportion is very difficult, due to the inconsistency in the chemical composition ratio given for the raw materials.

Fuzzy neural network model is used to obtain a desired quality of clinker. The raw mix as per the norms of cement industries should maintain the ranges like LSF 1.02 to 1.08, SM 2.35 to 2.55 and AM 0.95 to 1.25, which are the key factors for the burnability of clinker to obtain a good quality of cement. Fuzzy control theory method is used to minimize waste cement kiln dust. Fuzzy control theory allows varying degrees of set membership based on a membership function defined over a range of values. The membership function usually varies from 0 to 1.

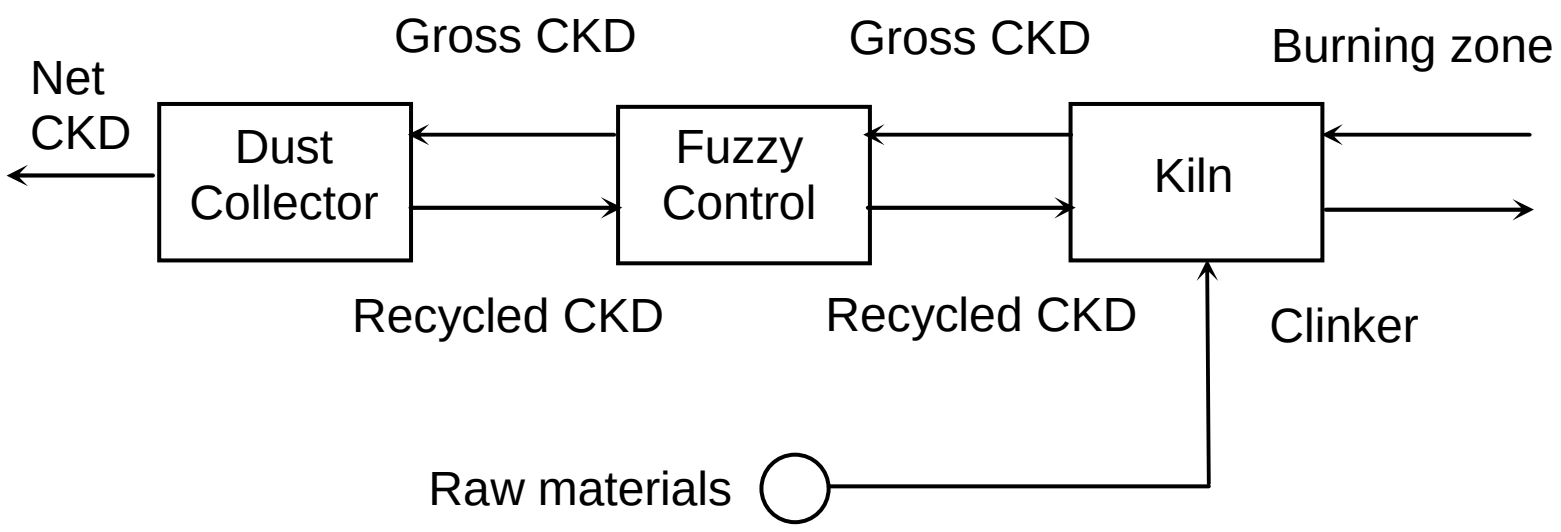

**FIGURE 2: Fuzzy control in kiln system**

## 2.2 Monitoring and control of the system using FCT and improvement of burning zone and combustion

Monitoring and control of the system is the most effective method towards CKD reduction in environment. CKD consists mainly of raw materials, which contain volatile compounds, therefore, tracking and control of the volatile compounds throughout the system often allows for the minimal CKD. The initial step in our plan towards CKD reduction is to identify the amount of the CKD. Here the indirect weighing method is applied to identify the amount of the CKD. Calculating sulphur/alkali ratio is a good indication of a possible imbalance. This ratio is calculated as the molar ratio of $SO_3/(K_2O)+Na_2O)$ in kiln load material.

CKD VOLATILE ANALYSIS

| Volatile | Molecular Weight |
|---|---|
| $Na_2O$ | 62 |
| $K_2O$ | 94.2 |
| $SO_3$ | 80 |

Ratio of alkali = $SO_3$ /$K_2O$ + $Na_2O$ = 80/156.2 = 0.512

This ratio should be between the values 0.5 to 1.5 in Kiln load material. The industry knows upto 40% of CKD exits, when the alkali ratio is between the values 0.5 to 1.5. But they cannot say exactly how much percentage of CKD waste comes from kiln by using the ratio of alkali in the online process. If industry knows this correct percentage of CKD in the online process, they can change some condition in the kiln and thus reduce the CKD in the online process. We adopt fuzzy control to estimate the percentage of CKD by using the ratio of alkali. The alkali ratio, kiln load material in tons and percentage of CKD are measured from the past happening process in kiln on a scale from 0.5 to 1.5, 5 to 25 tons and 0 to 40% respectively.

That is we assign the sulphur/alkali ratio shortly termed as alkali ratio, alkali ratio to be approximately low (L) when its value is 0.5, medium (M) when its value is 1 high (H) when its value is 1.5. In a similar way we give kiln load material ≅ {5 tons [first stage (FS)], 15 tons [second stage (SS)] and 25 tons [third stage (TS)]}. Percentage of CKD ≅ {0 [very less (VL)], 10 [less (L)], 20 [medium (M)], 30 [high (H)] and 40 [very high (VH)]}. ('≅' Denotes approximately equal). The terms of these parameters are presented in figures 2.2.1, 2.2.2 and 2.2.3.

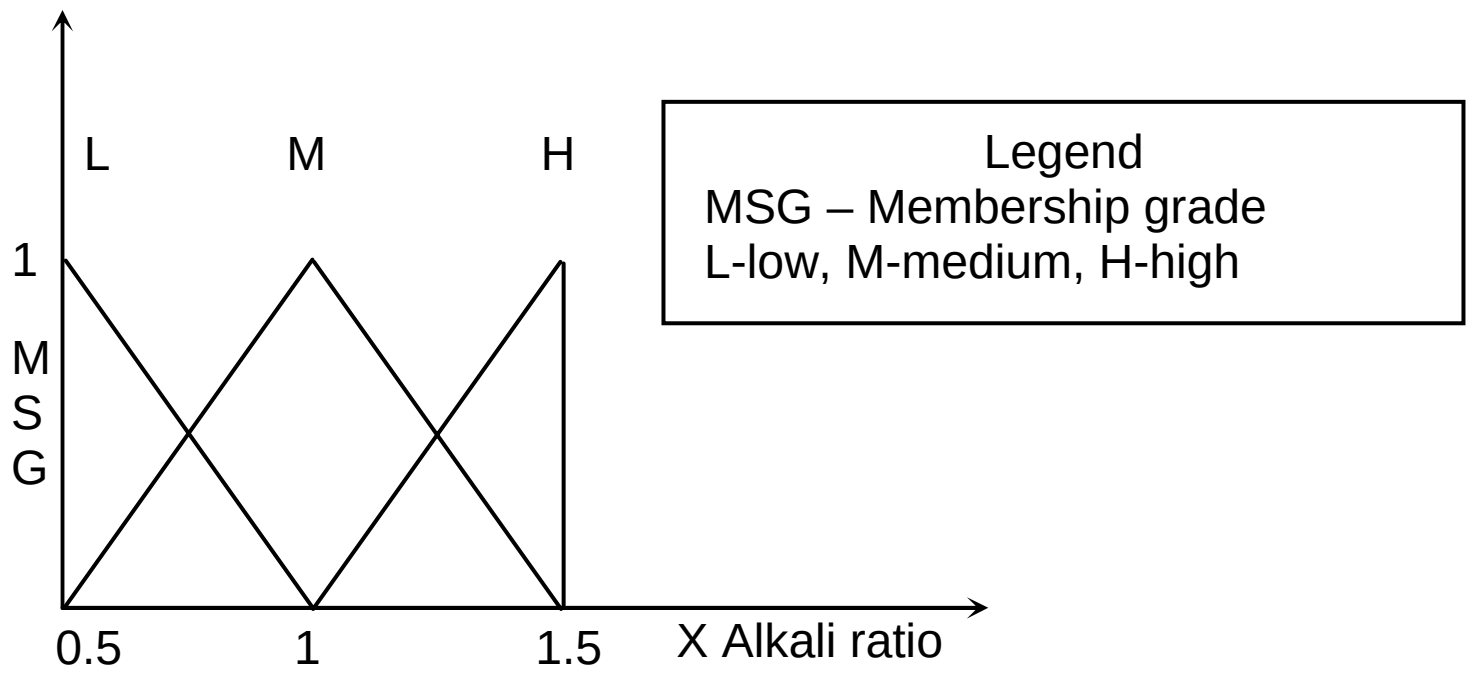

FIGURE 2.2.1: Alkali ratio- input parameter

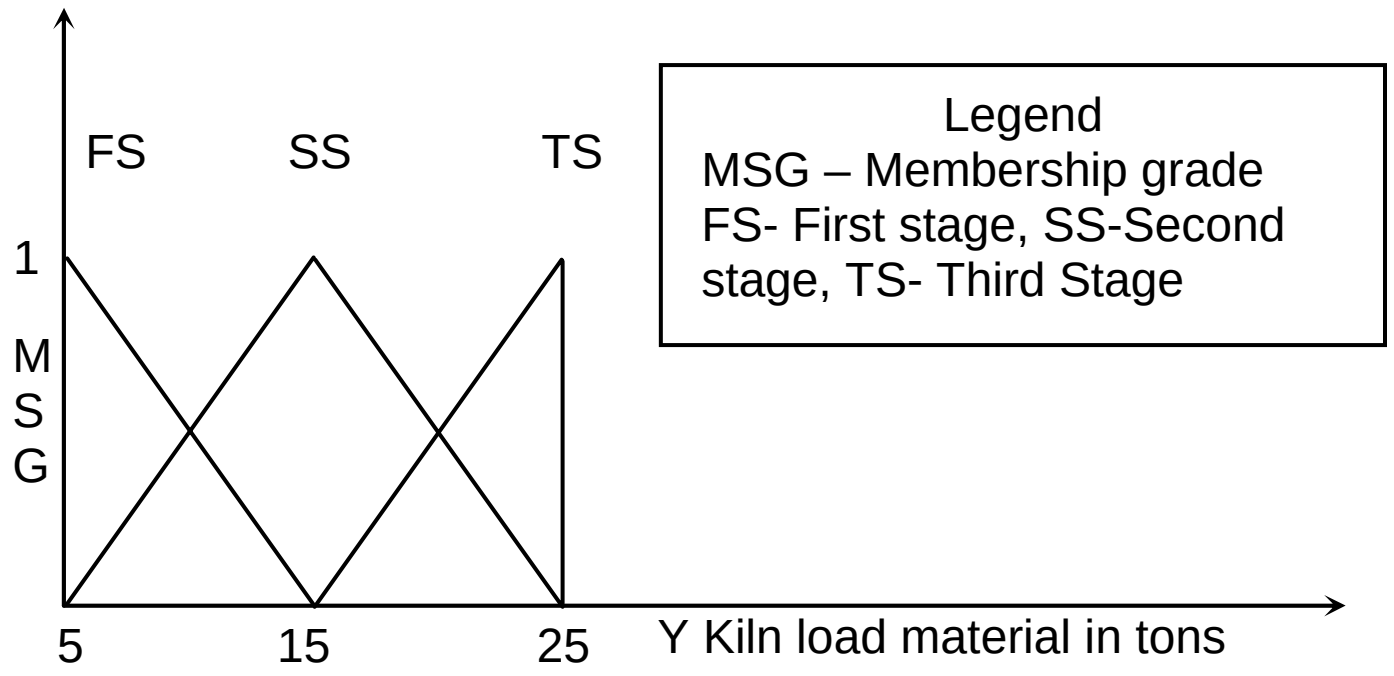

FIGURE 2.2.2: Kiln load material in tons-output parameter

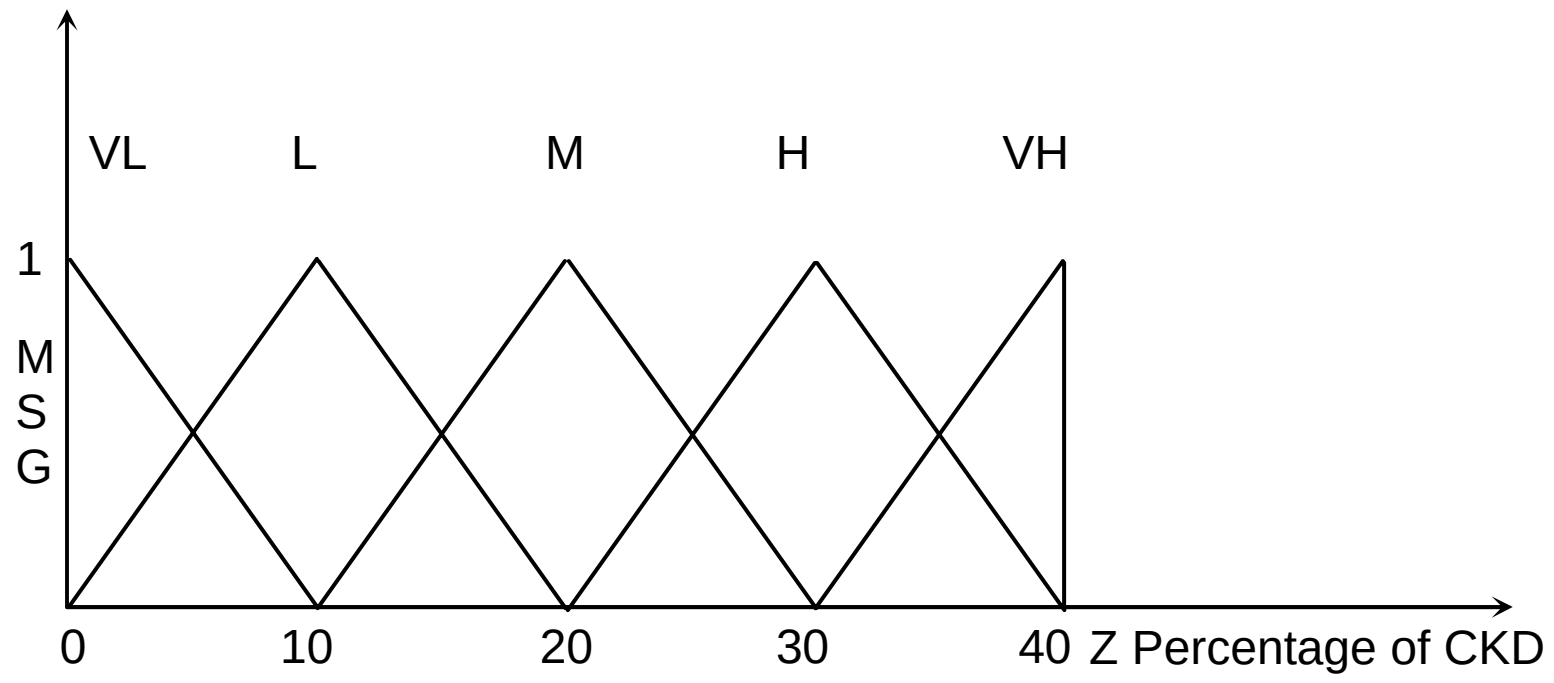

Legend
MSG – Membership grade
VL- very less, L-low, M-medium,
H-high, VH- very high

FIGURE 2.2.3: Percentage of CKD – output parameter

For the terms of alkali ratio, kiln load material in tons and percentage of CKD we give the following membership functions:

$$\mu(X)_{\text{alkali ratio}} = \begin{cases} \mu_L(X) = (1-X)/0.5 & 0.5 \le X \le 1 \\ \mu_M(X) = \begin{cases} (X-0.5)/0.5 & 0.5 \le X \le 1 \\ (1.5-X)/0.5 & 1 \le X \le 1.5 \end{cases} \\ \mu_H(X) = (X-1)/0.5 & 1 \le X \le 1.5 \end{cases} \qquad (2.2.1)$$

$$\mu(Y)_{\substack{\text{kiln} \\ \text{ratio in tons}}} = \begin{cases} \mu_{FS}(Y) = (15-Y)/10 & 5 \le Y \le 15 \\ \mu_{SS}(Y) = \begin{cases} (Y-5)/10 & 5 \le Y \le 15 \\ (25-Y)/10 & 15 \le Y \le 25 \end{cases} \\ \mu_{TS}(Y) = (Y-15)/10 & 15 \le Y \le 25 \end{cases} \qquad (2.2.2)$$

$$
\mu(Z)_{\substack{\text{percentage}\\ \text{of CKD}}} = \begin{cases} \mu_{VL}(Z) = (10-Z)/10 & 0 \le Z \le 10 \\ \mu_L(Z) = \begin{cases} Z/10 & 0 \le Z \le 10 \\ (20-Z)/10 & 10 \le Z \le 20 \end{cases} \\ \mu_M(Z) = \begin{cases} (Z-10)/10 & 10 \le Z \le 20 \\ (30-Z)/10 & 20 \le Z \le 30 \end{cases} \\ \mu_H(Z) = \begin{cases} (Z-20)/10 & 20 \le Z \le 30 \\ (40-Z)/10 & 30 \le Z \le 40 \end{cases} \\ \mu_{VH}(Z) = (Z-30)/10 & 30 \le Z \le 40 \end{cases} \quad (2.2.3)
$$

By applying the "if … and … then" rules [refer 11] to the three-membership functions $\mu(X)$, $\mu(Y)$ and $\mu(Z)$ we get the following table of rules.

The rules given in Table 2.2.1 read as follows :

Table 2.2.1

| Y / X | FS | SS | TS |
|---|---|---|---|
| L | VL | M | H |
| M | L | M | H |
| H | M | H | VH |

For example :

If alkali ratio is L and kiln load material in tons is FS then percentage of CKD is VL. If alkali ratio is H and kiln load material in tons is TS then percentage of CKD is VH; and so on.

Rules of evaluation using the membership functions defined by the equation (2.2.1) and (2.2.2), if alkali ratio is 1.2 and kiln load material is 17 tons we get the fuzzy inputs as $\mu_M(1.2) = 0.6$, $\mu_H(1.2) = 0.4$, $\mu_{SS}(17) = 0.8$ and $\mu_{TS}(17) = 0.2$. Induced decision table for percentage of CKD is as follows.

Table 2.2.2

| X \ Y | 0 | $\mu_{SS}$ (17) = 0.8 | $\mu_{TS}$(17) = 0.2 |
|---|---|---|---|
| 0 | 0 | 0 | 0 |
| $\mu_M$(1.2)=0.6 | 0 | $\mu_M$(Z) | $\mu_H$(Z) |
| $\mu_H$(1.2)=0.4 | 0 | $\mu_H$(Z) | $\mu_{VH}$(Z) |

Conflict resolutions of the four rules is as follows:

Rule 1 : If X is M and Y is SS then Z is M
Rule 2 : If X is M and Y is TS then Z is H
Rule 3 : If X is H and Y is SS then Z is H
Rule 4 : If X is H and Y is TS then Z is VH

Now, using Table 2.2.2 we calculate the strength values of the four rules as 0.6, 0.2, 0.4 and 0.2. Control output for the percentage of CKD is given in table 2.2.3.

Table 2.2.3

| X \ Y | 0 | $\mu_{SS}$ (17) = 0.8 | $\mu_{TS}$(17) = 0.2 |
|---|---|---|---|
| 0 | 0 | 0 | 0 |
| $\mu_M$(1.2)=0.6 | 0 | min{[0.6, $\mu_M$(Z)]} | min{[0.2, $\mu_H$(Z)]} |
| $\mu_H$(1.2)=0.4 | 0 | min{[0.4, $\mu_H$(Z)]} | min{[0.2, $\mu_{VH}$(Z)]} |

To find the aggregate(agg) of the control outputs, we obtain the maximum of the minimum. This is given by the following figure 2.2.4, that is $\mu_{agg}$ (Z) = max {min {[0.6, $\mu_M$(Z)] min {[0.4, $\mu_H$(Z)],)], min [0.2, $\mu_{vH}$(Z)]}. By applying the mean of maximum method for defuzzification that is the intersection points of the line $\mu$ = 0.6 with the triangular fuzzy number $\mu_M$(Z) in equation (2.2.3) we get the crisp output to be 20%.

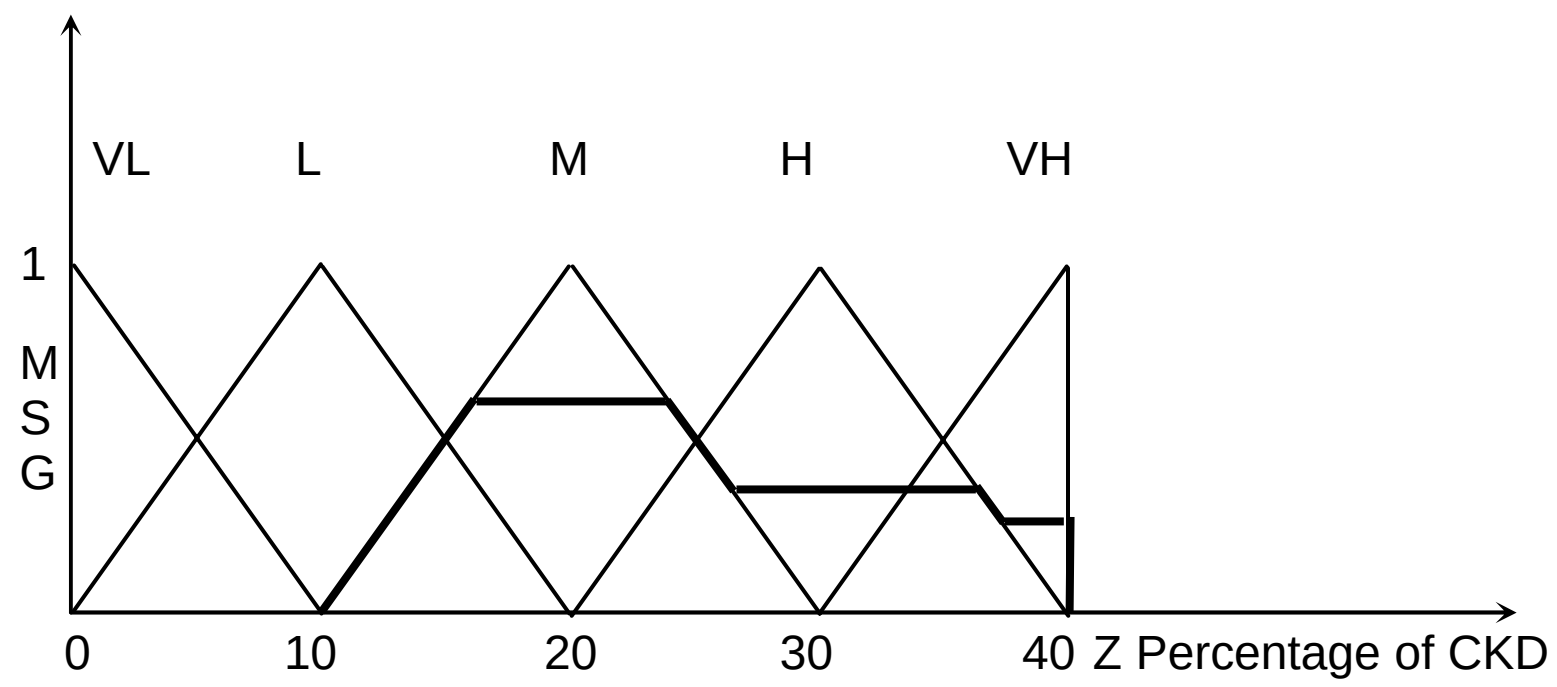

FIGURE 2.2.4: Aggregated output and defuzzificztion for the percentage of CKD

Rules of evaluation using the membership function defined by the equation (1) and (2), if alkali ratio is 0.5 and kiln load material is 5 tons we get the fuzzy inputs as $\mu_L(0.5) = 1$, $\mu_H(0.5) = 0$, $\mu_{rs}(5) = 1$ and $\mu_{ss}(5) = 0$. Induced decision table for percentage of CKD is as follows.

Table 2.2.4

| X \ Y | $\mu_{FS}(5) = 1$ | $\mu_{SS}(5) = 0$ | 0 |
|---|---|---|---|
| $\mu_L(0.5) = 1$ | $\mu_{VL}(Z)$ | $\mu_M(Z)$ | 0 |
| $\mu_M(0.5) = 0$ | $\mu_L(Z)$ | $\mu_M(Z)$ | 0 |
| 0 | 0 | 0 | 0 |

Conflict resolutions of the four rules is as follows:

Rule 1 : If X is L and Y is FS then Z is VL
Rule 2 : If X is L and Y is SS then Z is M
Rule 3 : If X is M and Y is FS then Z is L
Rule 4 : If X is M and Y is SS then Z is M.

Now, using Table 2.2.4 we calculate the strength values of the four rules as 1, 0, 0 and 0. Control output for the percentage of CKD is given in Table 2.2.5.

Table 2.2.5

| Y / X | $\mu_{FS}$ (5) = 1 | $\mu_{SS}$(5) = 0 | 0 |
|---|---|---|---|
| $\mu_L$ (0.5) = 1 | min {[1, $\mu_{VL}$(Z)]} | min {[0, $\mu_M$(Z)]} | 0 |
| $\mu_H$(0.5) = 0 | min {[0, $\mu_L$(Z)]} | min {[0, $\mu_M$(Z)]} | 0 |
| 0 | 0 | 0 | 0 |

To find the aggregate of the control outputs, we obtain the maximum of the minimum. This is given by the following figure 2.2.5 that is $\mu_{agg}$ (Z) = {min {l, $\mu_{VL}$(Z)]}, min{[0, $\mu_M$(Z)]}, min {[0, $\mu_L$ (Z)]}. By applying the mean of maximum method for defuzzification that is the intersection points of the line $\mu$ =1 with the triangular fuzzy number $\mu_{VL}$(Z) in equation (3) and get the crisp output to be 0%.

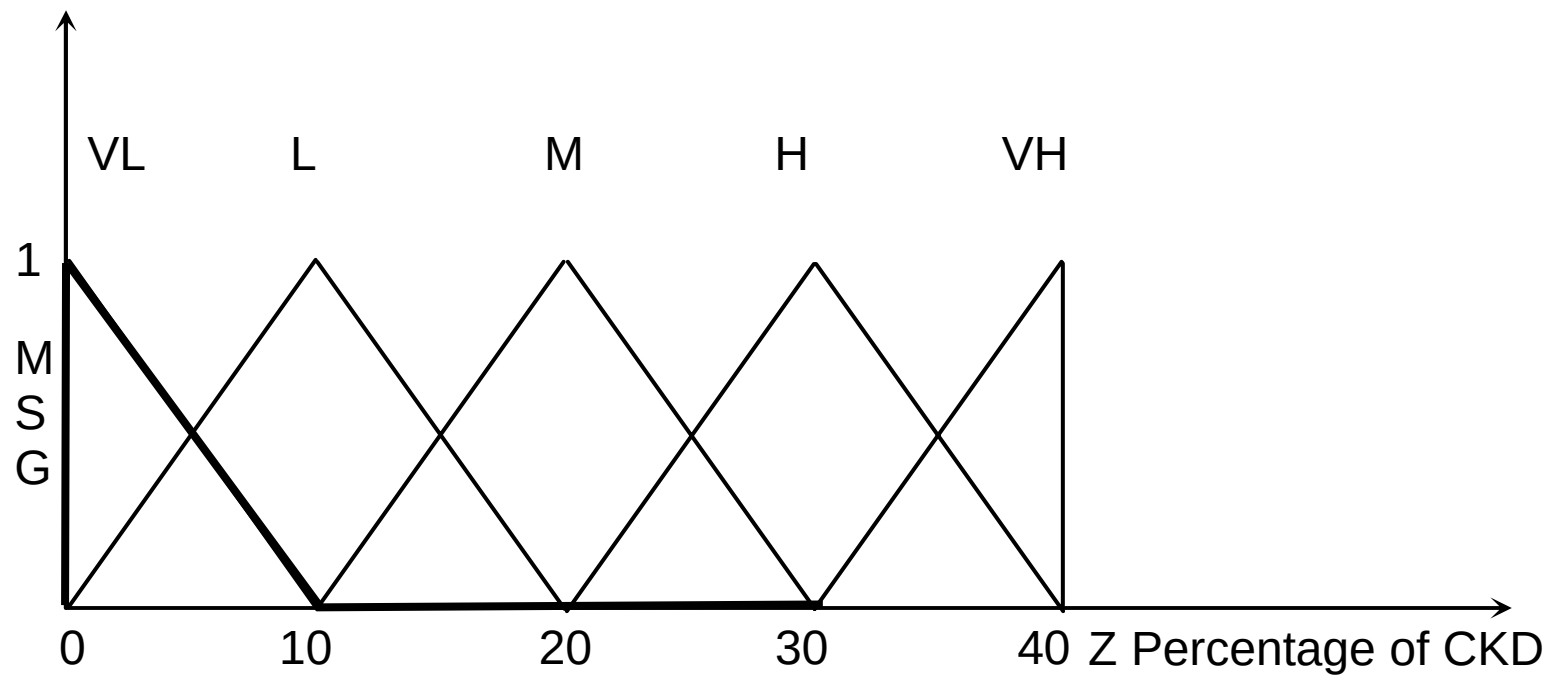

FIGURE 2.2.5: Aggregated output and defuzzificztion for the percentage of CKD

Rules of evaluation using the membership function defined by the equations (1) and (2), if alkali ratio is 1 and kiln load material is 15 tons we get the fuzzy inputs as $\mu_L$(1) = 0, $\mu_H$ (1) = 0 and $\mu_m$ (1) = 1, $\mu_{FS}$ (15) = 0, $\mu_{SS}$(15) = 1, $\mu_{TS}$ (15) = 0, Induced decision table for percentage of CKD is as follows.

Table 2.2.6

| Y / X | $\mu_{FS}(15) = 0$ | $\mu_{SS}(15) = 1$ | $\mu_{TS}(15) = 0$ |
|---|---|---|---|
| $\mu_L(1) = 0$ | $\mu_{VL}(Z)$ | $\mu_M(Z)$ | $\mu_H(Z)$ |
| $\mu_M(1) = 1$ | $\mu_L(Z)$ | $\mu_M(Z)$ | $\mu_H(Z)$ |
| $\mu_H(1) = 0$ | $\mu_M(Z)$ | $\mu_H(Z)$ | $\mu_{VH}(Z)$ |

Conflict resolutions of the nine rules is as follows :

Rule 1 : If X is L and Y is FS then Z is VL
Rule 2 : If X is L and Y is SS then Z is M
Rule 3 : If X is L and Y is TS then Z is H
Rule 4 : If X is M and Y is FS then Z is L.
Rule 5 : If X is M and Y is SS then Z is M.
Rule 6 : If X is M and Y is TS then Z is H.
Rule 7 : If X is H and Y is FS then Z is L.
Rule 8 : If X is H and Y is SS then Z is M.
Rule 9 : If X is H and Y is TS then Z is H.

Now, using Table 2.2.6 we calculate the strength values of the nine rules as 0, 0, 0, 0, 1, 0, 0, 0, 0. Control output for the percentage of CKD is given in Table 2.2.7.

Table 2.2.7

| Y / X | $\mu_{FS}(15) = 0$ | $\mu_{SS}(15) = 1$ | $\mu_{TS}(15) = 1$ |
|---|---|---|---|
| $\mu_L(1)=0$ | $\min\{[0,\mu_{VL}(Z)]\}$ | $\min\{[0,\mu_M(Z)]\}$ | $\min\{[0,\mu_H(Z)]\}$ |
| $\mu_M(1)=1$ | $\min\{[0,\mu_L(Z)]\}$ | $\min\{[0,\mu_M(Z)]\}$ | $\min\{[0,\mu_H(Z)]\}$ |
| $\mu_H(1)=0$ | $\min\{[0,\mu_M(Z)]\}$ | $\min\{[0,\mu_H(Z)]\}$ | $\min\{[0,\mu_H(Z)]\}$ |

To find the aggregate of the control outputs, we obtain the maximum of the minimum. This is given by the following figure 2.2.6, that is $\mu_{agg}(Z) = \max\{\min\{0, \mu_{VL}(Z)]\}, \min\{[0, \mu_M(Z)]\}, \min\{[0, \mu_L(Z)]\}, \{\min\{l, \mu_H(Z)]\}, \min\{[0, \mu_{VH}(Z)]\}$. By applying the mean of maximum method for defuzzification that is the intersection points of the line $\mu = 1$ with the triangular fuzzy number $\mu_{VL}(Z)$ in equation (2.2.3) and get the crisp output to 20%.

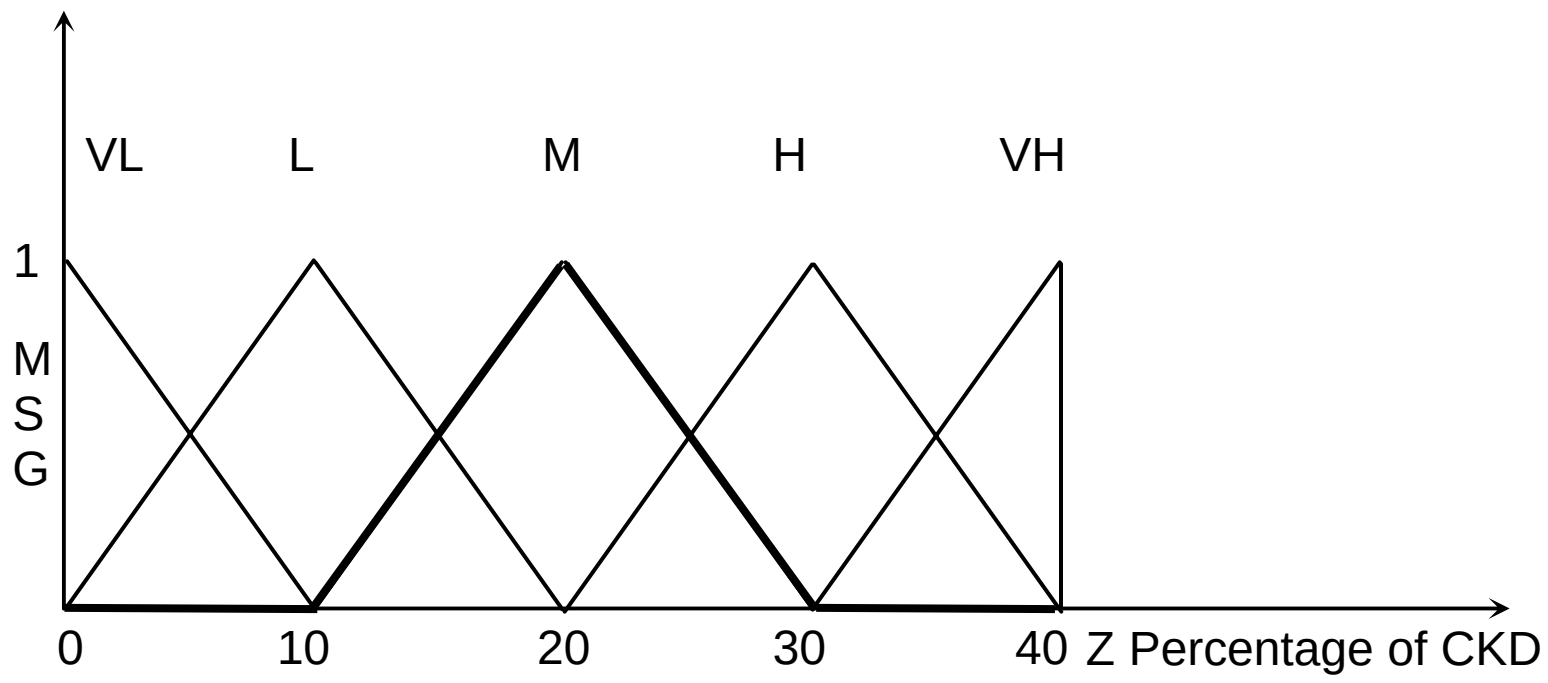

FIGURE 2.2.6: Aggregated output and defuzzification for the percentage of CKD

Rules of evaluation using the membership function defined by the equations (2.2.1) and (2.2.2), if alkali ratio is 1.5 and kiln load material is 25 tons we get the fuzzy inputs as $\mu_M$ (1.5) = 0, $\mu_H$ (1.5) = 1, $\mu_{SS}$ (25) = 0 and $\mu_{TS}$ (25) = 1. Induced decision table for percentage of CKD is as follows.

Table 2.2.8

| X \ Y | 0 | $\mu_{SS}(25) = 0$ | $\mu_{TS}(25) = 1$ |
|---|---|---|---|
| 0 | 0 | 0 | 0 |
| $\mu_M(1.5) = 0$ | 0 | $\mu_M(Z)$ | $\mu_H(Z)$ |
| $\mu_H(1.5) = 1$ | 0 | $\mu_H(Z)$ | $\mu_{VH}(Z)$ |

Conflict resolutions of the four rules is as follows :

Rule 1 : If X is M and Y is SS then Z is M
Rule 2 : If X is M and Y is TS then Z is H
Rule 3 : If X is H and Y is SS then Z is H
Rule 4 : If X is H and Y is TS then Z is VH.

Now, using Table 2.2.8 we calculate the strength values of the four rules as 0, 0, 0 and 1 Control output for the percentage of CKD is given in Table 2.2.9.

Table 2.2.9

| Y / X | 0 | $\mu_{SS}(25) = 0$ | $\mu_{TS}(25) = 1$ |
|---|---|---|---|
| 0 | 0 | 0 | 0 |
| $\mu_M(1.5) = 0$ | 0 | min [0, $\mu_M(Z)$] | min [0, $\mu_H(Z)$] |
| $\mu_H(1.5) = 1$ | 0 | min [0, $\mu_H(Z)$] | min [1, $\mu_{VH}(Z)$] |

To find the aggregate of the control outputs, we obtain the maximum of the minimum. This is given by the following figure 2.2.7, that is $\mu_{agg}$ (Z) = max {min {0, $\mu_M(Z)$]} , min {[0, $\mu_H(Z)$]}, min{[1, $\mu_{VH}$ (Z)]}. By applying the mean of maximum method for defuzzification that is the intersection points of the line $\mu$ =1 with the triangular fuzzy number $\mu_{VH}(Z)$ in equation 2.2.3 and get the crisp output to 40%.

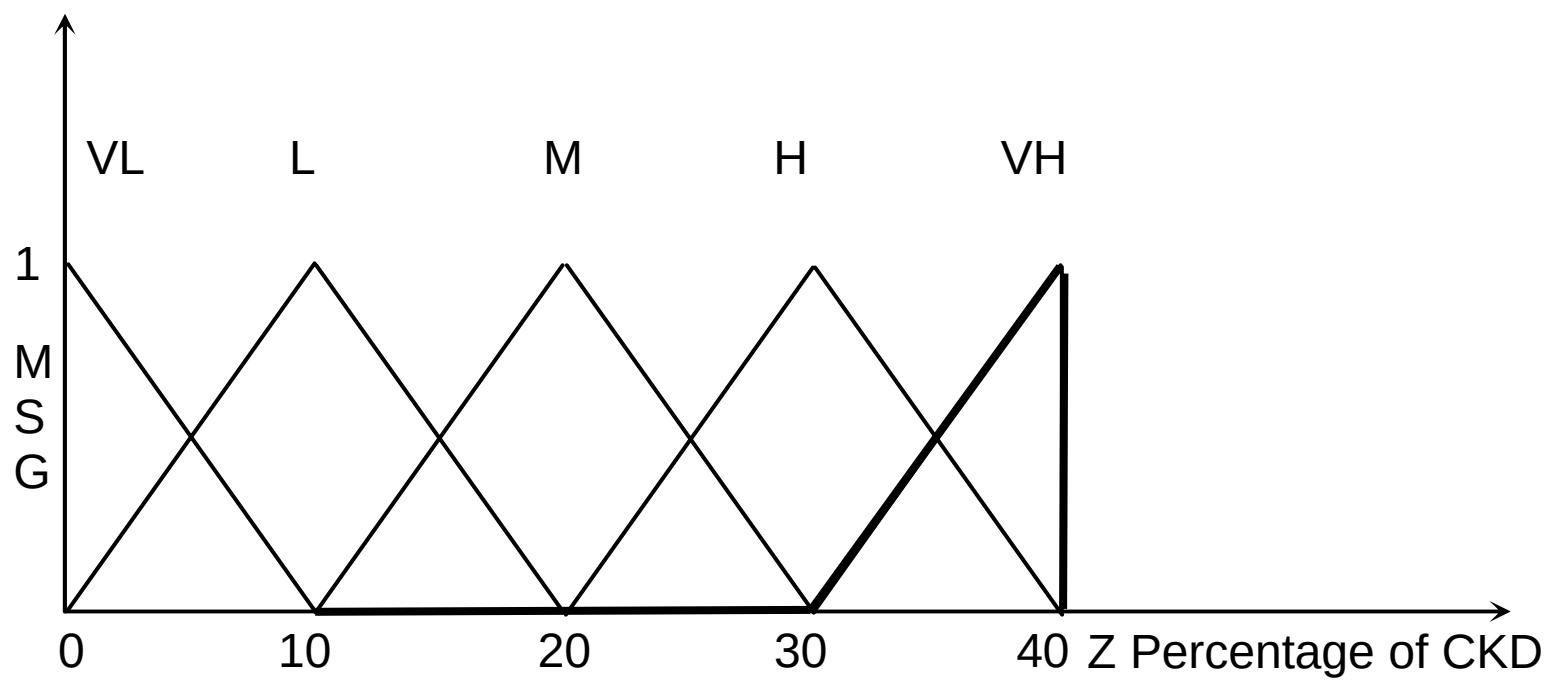

FIGURE 2.2.7: Aggregated output and defuzzification for the percentage of CKD

From our study we suggest in the online process to reduce (or) minimize the amount of CKD in the industry one should change the condition of fuel burning system and other system in kiln from time to time depending on the percentage of CKD in tons given above.

## 2.3 Determination of gas volume setpoint and temperature set point for CKD processing

The total CKD dust carried out from the kiln is again returned to the kiln as a feed (Recycled CKD). After recycled process, we get some amount of remaining CKD from kiln, which is disposed in the environment(as a waste polluting the environment). Most of the cement factory uses electrostatic precipitator(ESP) method for recycling process of CKD, as it operates by gas volume and temperature. In ESP, we mainly concentrate on gas volume in $m^3$/minute and temperature degree in celsius. The range of gas volume is varying from 11865 to 15174 $m^3$/minute and temperature is varying from $350^oC$ to $450^oC$. When in the recycle; the clinker is got from the reproduced dust to clinker by pre heater in dust collector(ESP). Generally an industry to minimize the net CKD dust upto 20% by reprocessing method randomly chooses the gas volume and temperature from the range of gas volume (11865 to 15174 $m^3$/minute) and temperature ($350^oC$ to $450^oC$) respectively. Since the gas volume and temperature are main concerns on ESP, the reprocessing directly depends on gas volume and temperature. The randomly choosing of the gas volume set point and temperature set point from the ranges of gas volume and temperature is uncertain and does not usually give the desired outcomes so, this gas volume and temperature affect the CKD reprocessing largely. In order to over come these problems we use fuzzy control to find the set point of gas volume and temperature in ESP, which is described in the following. The ranges of gas volume, temperature and percentage of net CKD are measured from the past happening data in ESP on a scale, are 11865 to 15174 $m^3$/minute, $350^oC$ to $450^oC$ and 0 to 20% respectively. Temperature ≅ {$350^oC$ [low (L)], $400^oC$ [medium (M)] and $450^oC$ [high (H)]}. Gas volume ≅ 11865 to 15174$m^3$/min [first stage (FS)], 13020 $m^3$/min [second stage (SS)], and 15174 $m^3$/min [third stage (TS)}. Percentage of net CKD ≅ {0[very less (VL)], 5 [less (L), 10[medium (M)], 15 [high (H)] and 20 [very high (VH)}. The terms of these parameters are presented in figures 2.3.1 or 2.3.3.

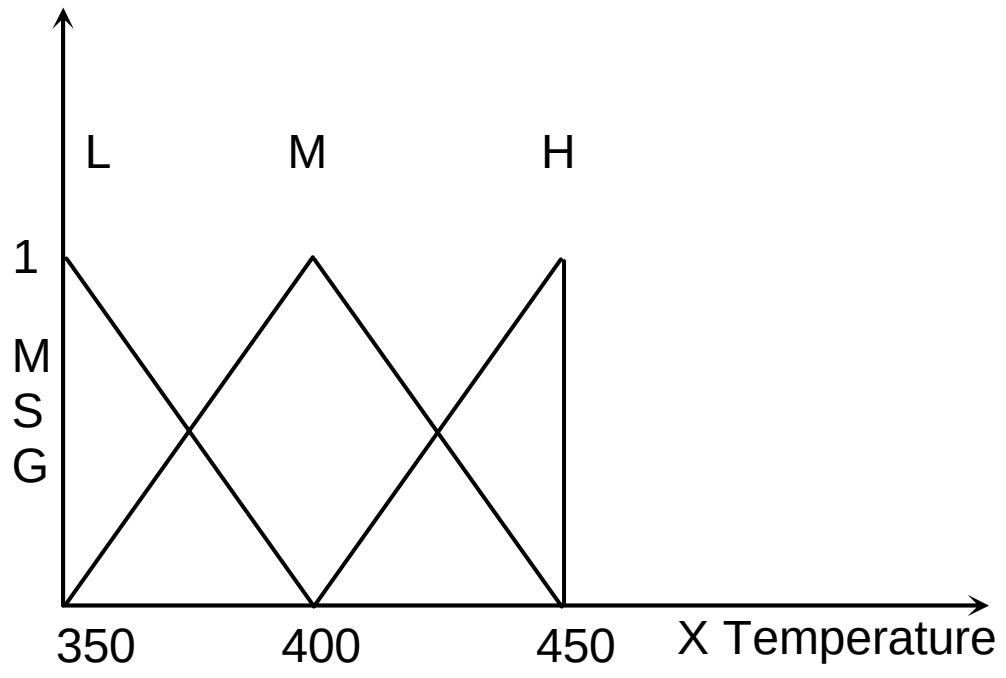

FIGURE 2.3.1: Temperature - input parameter

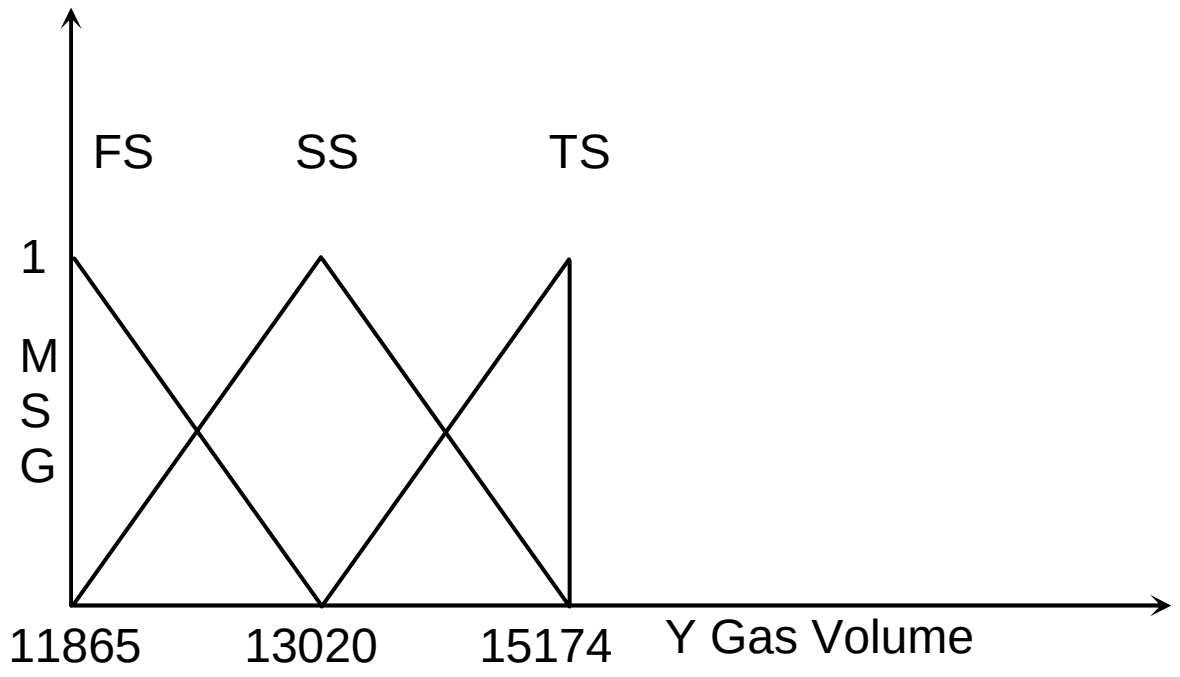

FIGURE 2.3.2: Gas Volume- input parameter

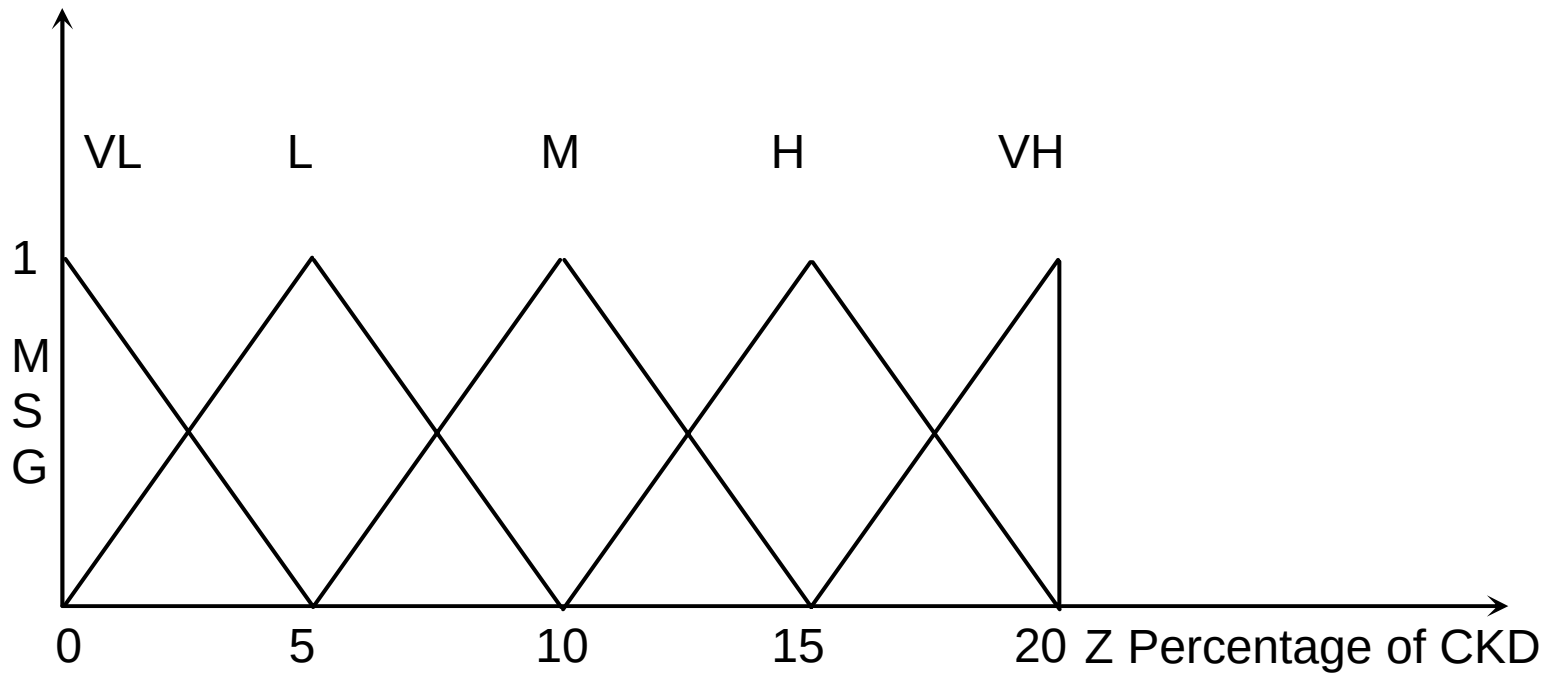

FIGURE 2.3.3: Percentage of net CKD-output parameter

For the terms of temperature, gas volume and percentage of net CKD we give the following membership functions :

$$\mu(X)_{\text{temperature}} = \begin{cases} \mu_L(X) = (400 - X)/50 & 350 \le X \le 400 \\ \mu_M(X) = \begin{cases} (X - 350)/50 & 350 \le X \le 400 \\ (450 - X)/50 & 400 \le X \le 450 \end{cases} \\ \mu_H(X) = (X - 400)/50 & 400 \le X \le 450 \end{cases} \quad (2.3.1)$$

$$\mu(Y)_{\substack{\text{gas}\\ \text{volume}}} = \begin{cases} \mu_{FS}(Y) = (13020 - Y)/1155 & 11865 \le Y \le 13020 \\ \mu_{SS}(Y) = \begin{cases} (Y - 11865)/1155 & 11865 \le Y \le 13020 \\ (15174 - Y)/2154 & 13020 \le Y \le 15174 \end{cases} \\ \mu_{TS}(Y) = (Y - 13020)/2154 & 13020 \le Y \le 15174 \end{cases}$$

(2.3.2)

$$\mu(Z)_{\substack{\text{percentage}\\ \text{of net CKD}}} = \begin{cases} \mu_{VL}(Z) = (5 - Z)/5 & 0 \le Z \le 5 \\ \mu_L(Z) = \begin{cases} Z/5 & 0 \le Z \le 5 \\ (10 - Z)/5 & 5 \le Z \le 10 \end{cases} \\ \mu_M(Z) = \begin{cases} (Z - 5)/5 & 5 \le Z \le 10 \\ (15 - Z)/5 & 10 \le Z \le 15 \end{cases} \\ \mu_H(Z) = \begin{cases} (Z - 10)/5 & 10 \le Z \le 15 \\ (20 - Z)/5 & 15 \le Z \le 20 \end{cases} \\ \mu_{VH}(Z) = (Z - 15)/5 & 15 \le Z \le 20 \end{cases} \quad (2.3.3)$$

By applying the if … and … then rules to the three-membership function μ(X), μ(Y) and μ(Z), we get the following table of rules. The rules given in Table 2.3.1 read as follows:

Table 2.3.1

| X \ Y | FS | SS | TS |
|---|---|---|---|
| L | VL | M | H |
| M | L | M | H |
| H | M | H | VH |

For example:

If temperature is L and gas volume is SS then percentage of net CKD is M.

If temperature is M and gas volume is TS then percentage of net CKD is H; and so on.

Rules of evaluation using the membership functions defined by the equation (2.3.1) and (2.3.2), if temperature is 430$^o$C and gas volume is 13080 $m^3$/min we get the fuzzy inputs as $\mu_M(430) = 0.4$, $\mu_H(430) = 0.6$, $\mu_{SS}(13080) = 0.97$ and $\mu_{TS}(3080) = 0.02$. Induced decision table for percentage of net CKD is as follows.

Table 2.3.2

| X \ Y | 0 | $\mu_{SS}(13080)= 0.97$ | $\mu_{TS}(13080)= 0.02$ |
|---|---|---|---|
| 0 | 0 | 0 | 0 |
| $\mu_M(430)=0.04$ | 0 | $\mu_M(Z)$ | $\mu_H(Z)$ |
| $\mu_H(430)=0.06$ | 0 | $\mu_H(Z)$ | $\mu_{VH}(Z)$ |

Conflict resolutions of the four rules is as follows :

Rule 1 : If X is M and Y is SS then Z is M
Rule 2 : If X is M and Y is TS then Z is H
Rule 3 : If X is H and Y is SS then Z is H
Rule 4 : If X is H and Y is TS then Z is VH.

Now, using Table 2.3.2 we calculate the strength values of the four rules as 0.4, 0.02, 0.06 and 0.02. Control output for the percentage of net CKD is given in Table 2.3.3.

Table 2.3.3

| X \ Y | 0 | $\mu_{SS}$ (13080) = 0.97 | $\mu_{TS}(13080) = 0.02$ |
|---|---|---|---|
| 0 | 0 | 0 | 0 |
| $\mu_M(430)=0.4$ | 0 | min {[0.4, $\mu_M(Z)$]} | min {[0.02, $\mu_H(Z)$]} |
| $\mu_H(430)=0.6$ | 0 | min {[0.6, $\mu_H(Z)$]} | min{[0.02, $\mu_{VH}(Z)$]} |

To find the aggregate of the control outputs, we obtain the maximum of the minimum. This is given by the following figure 2.3.4, that is $\mu_{agg}$ (Z) = max{min [0.4, $\mu_M(Z)$]}, min{[0.6, $\mu_M(Z)$]}, min {[0.02, $\mu_{VH}$ (Z)]}. We apply the mean of maximum method for defuzzification that is the intersection points of the line $\mu$ = 0.6 with the triangular fuzzy number $\mu_H(Z)$ in equation (2.3.3) and get the crisp output as 15 to 20%.

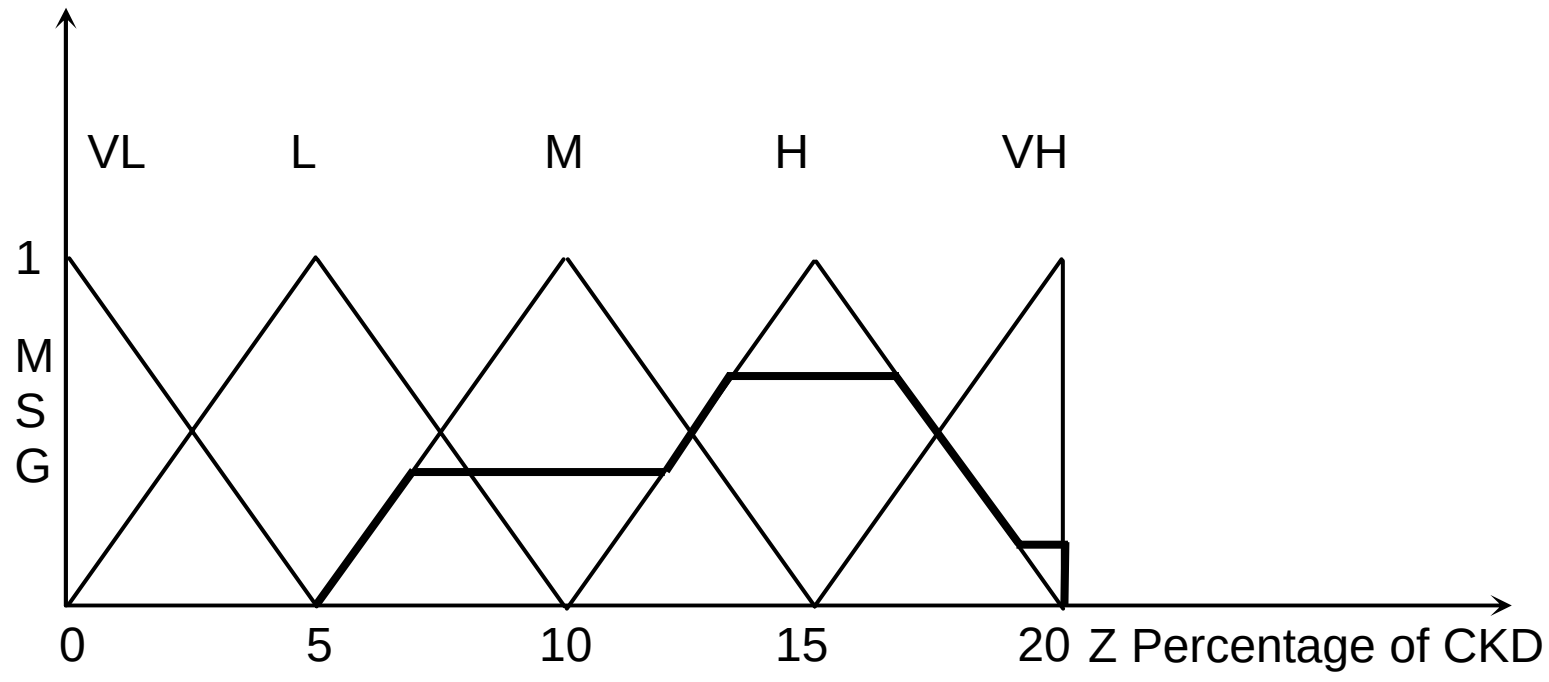

FIGURE 2.3.4: Terms of the output and defuzzification for the percentage of net CKD

Rules of evaluation using the membership functions defined by the equation (2.3.1) and (2.3.2), if temperature is 350°C and gas volume is 11865 $m^3$/min we get the fuzzy inputs as $\mu_L(350) = 1$, $\mu_H(350) = 0$, $\mu_{FS}(11865) = 1$ and $\mu_{SS}(11865) = 0$. Induced decision table for percentage of net CKD is as follows.

Table 2.3.4

| X \ Y | $\mu_{FS}$ (11865) = 1 | $\mu_{SS}(11865) = 0$ | 0 |
|---|---|---|---|
| $\mu_L(350) = 1$ | $\mu_{VL}(Z)$ | $\mu_M(Z)$ | 0 |
| $\mu_M(350) = 0$ | $\mu_L(Z)$ | $\mu_M(Z)$ | 0 |
| 0 | 0 | 0 | 0 |

Conflict resolutions of the four rules is as follows :

Rule 1 : If X is L and Y is FS then Z is VL

Rule 2 : If X is L and Y is SS then Z is M
Rule 3 : If X is M and Y is FS then Z is L
Rule 4 : If X is M and Y is SS then Z is M.

Now, using Table 2.3.4 we calculate the strength values of the four rules as 1, 0, 0 and 0. Control output for the percentage of net CKD is given in table 2.3.5.

Table 2.3.5

| X \ Y | $\mu_{FS}$ (11865) = 1 | $\mu_{SS}$(11865 ) = 0 | 0 |
|---|---|---|---|
| $\mu_L$(350)=0.4 | min {[1, $\mu_{VL}$(Z)]} | min {[0, $\mu_M$(Z)]} | 0 |
| $\mu_M$(350)= 0.6 | min {[0, $\mu_L$(Z)]} | min {[0, $\mu_M$(Z)]} | 0 |
| 0 | 0 | 0 | 0 |

To find the aggregate of the control outputs, we obtain the maximum of the minimum.

This is given by the following figure that is $\mu_{agg}$ (Z) = {min{[1, $\mu_{VL}$ (Z)]}, min{[0, $\mu_M$ (Z)]}, min {[0, $\mu_L$ (Z)]}. We apply the mean of maximum method for defuzzification that is the intersection points of the line $\mu$ = 1 with the triangular fuzzy number $\mu_{VL}$(Z) in equation (2.3.3) and get the crisp output as 0 to 5%.

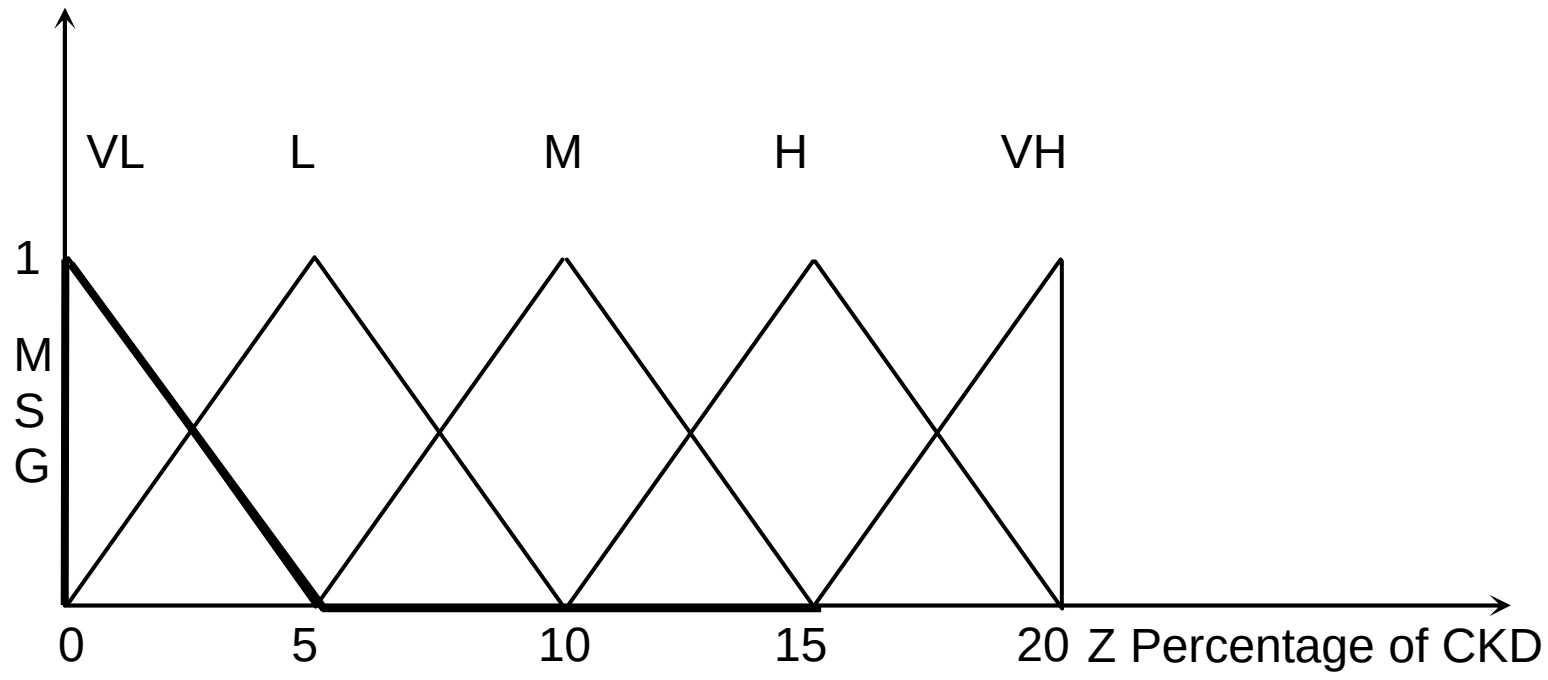

FIGURE 2.3.5: Terms of the output and defuzzification for the percentage of net CKD

Rules of evaluation using the membership functions defined by the equation (4) and (5), if temperature is 400°C and gas volume is 13020 $m^3$/min we get the fuzzy inputs as $\mu_L(400) = 0$, $\mu_M(400) = 1$, $\mu_H(400) = 0$, $\mu_{FS}(13020) = 0$, $\mu_{SS}(13020) = 1$ and $\mu_{TS}(13020) = 0$. Induced decision table 2.3.6 for percentage of net CKD is as follows.

Table 2.3.6

| X \ Y | $\mu_{FS}(13020)=0$ | $\mu_{SS}(13020)=1$ | $\mu_{TS}(13020)=0$ |
|---|---|---|---|
| $\mu_L(400)=0$ | $\mu_{VL}(Z)$ | $\mu_M(Z)$ | $\mu_H(Z)$ |
| $\mu_M(400)=1$ | $\mu_L(Z)$ | $\mu_M(Z)$ | $\mu_H(Z)$ |
| $\mu_H(400)=0$ | $\mu_M(Z)$ | $\mu_H(Z)$ | $\mu_{VH}(Z)$ |

Conflict resolutions of the nine rules is as follows:

Rule 1 : If X is L and Y is FS then Z is VL
Rule 2 : If X is L and Y is SS then Z is M
Rule 3 : If X is L and Y is TS then Z is H
Rule 4 : If X is M and Y is FS then Z is L.
Rule 5 : If X is M and Y is SS then Z is M.
Rule 6 : If X is M and Y is TS then Z is H.
Rule 7 : If X is H and Y is FS then Z is M.
Rule 8 : If X is H and Y is SS then Z is H.
Rule 9 : If X is H and Y is TS then Z is VH.

Now, using Table 2.3.6 we calculate the strength values of the nine rules as 0, 0, 0, 0, 1, 0, 0, 0, 0. Control output for the percentage of net CKD is given in Table 2.3.7.

Table 2.3.7

| X \ Y | $\mu_{FS}(13020) = 0$ | $\mu_{SS}(13020) = 1$ | $\mu_{TS}(13020) = 0$ |
|---|---|---|---|
| $\mu_L(400) =0$ | $\min\{[1,\mu_{VL}(Z)]\}$ | $\min\{[0,\mu_M(Z)]\}$ | $\min\{[0,\mu_H(Z)]\}$ |
| $\mu_M(400)=1$ | $\min\{[0, \mu_L(Z)]\}$ | $\min\{[1,\mu_M(Z)]\}$ | $\min\{[0,\mu_H(Z)]\}$ |
| $\mu_H(400)=0$ | $\min\{[0, \mu_M(Z)]\}$ | $\min\{[0,\mu_H(Z)]\}$ | $\min\{[0,\mu v_H(Z)]\}$ |

To find the aggregate(agg) of the control outputs, we obtain the maximum of the minimum.

This is given by the following figure 2.3.6, that is $\mu_{agg}(Z)$ = max {min {[0, $\mu_{VL}(Z)$]}, min {[1, $\mu_M(Z)$],)], min {[0, $\mu_L(Z)$]}, min {[0, $\mu_H(Z)$]}, min {[0, $\mu_{VH}(Z)$]}. We apply the mean of maximum method for defuzzification that is the intersection points of the line $\mu = 1$ with the triangular fuzzy number $\mu_M(Z)$ in equation (2.3.3) we get the crisp output to be 10 % to 15 %.

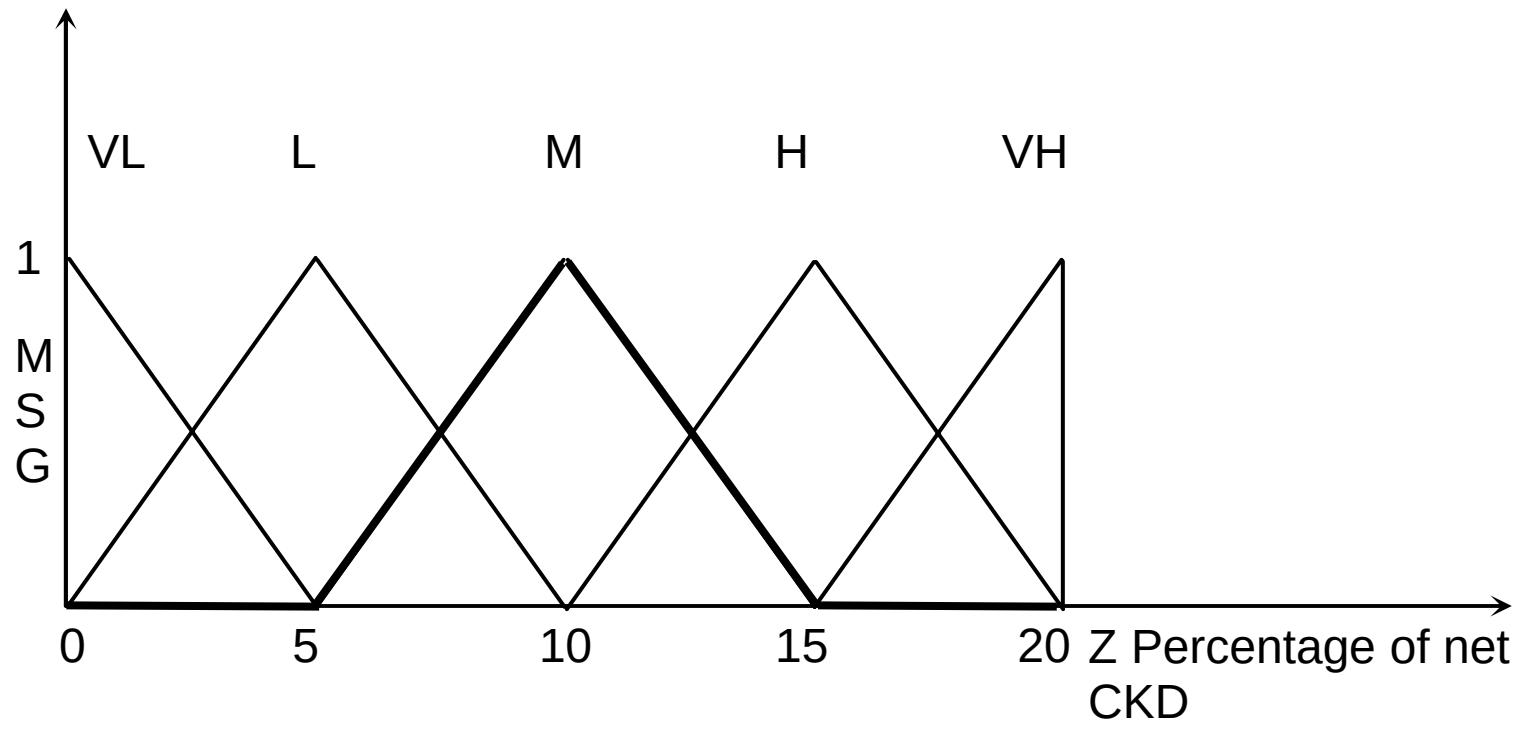

FIGURE 2.3.6: Aggregated output and defuzzification for the percentage of net CKD

Rules of evaluation using the membership functions defined by the equation (2.3.1) and (2.3.2), if temperature is 450°C and gas volume is 15174 $m^3$/min we get the fuzzy inputs as $\mu_M(450) = 0$, $\mu_H(450) = 1$, $\mu_{SS}(15174) = 0$ and $\mu_{TS}(15174) = 1$. Induced decision table for percentage of net CKD is as follows.

Table 2.3.8

| X \ Y | 0 | $\mu_{SS}(15174) = 0$ | $\mu_{TS}(15174) = 1$ |
|---|---|---|---|
| 0 | 0 | 0 | 0 |
| $\mu_M(450) = 0$ | 0 | $\mu_M(Z)$ | $\mu_H(Z)$ |
| $\mu_H(450) = 1$ | 0 | $\mu_H(Z)$ | $\mu_{VH}(Z)$ |

Conflict resolutions of the four rules is as follows:

Rule 1 : If X is M and Y is SS then Z is M
Rule 2 : If X is M and Y is TS then Z is H
Rule 3 : If X is H and Y is SS then Z is H
Rule 4 : If X is H and Y is TS then Z is VH.

Now, using Table 2.3.8 we calculate the strength values of the four rules as 0, 0, 0 and 1. Control output for the percentage of net CKD is given in Table 2.3.9.

Table 2.3.9

| X \ Y | 0 | $\mu_{SS}(15174) = 0$ | $\mu_{TS}(15174) = 0$ |
|---|---|---|---|
| 0 | 0 | 0 | 0 |
| $\mu_M(450)=0$ | 0 | $\min[0, \mu_M(Z)]$ | $\min[0,\mu_H(Z)]$ |
| $\mu_H(400)= 1$ | 0 | $\min[0, \mu_H(Z)]$ | $\min[1,\mu_{VH}(Z)]$ |

To find the aggregate(agg) of the control outputs, we obtain the maximum of the minimum.

This is given by the following figure 2.3.7, that is $\mu_{agg}$ (Z) = max{min{[0, $\mu_M$(Z)]}, min{[0, $\mu_H$(Z)],)], min{[1, $\mu_{VH}$(Z)]}. We apply the mean of maximum method for defuzzification that is the intersection points of the line $\mu$ = 1 with the triangular fuzzy number $\mu_M$(Z) in equation (2.3.3) we get the crisp output to be 20 %.

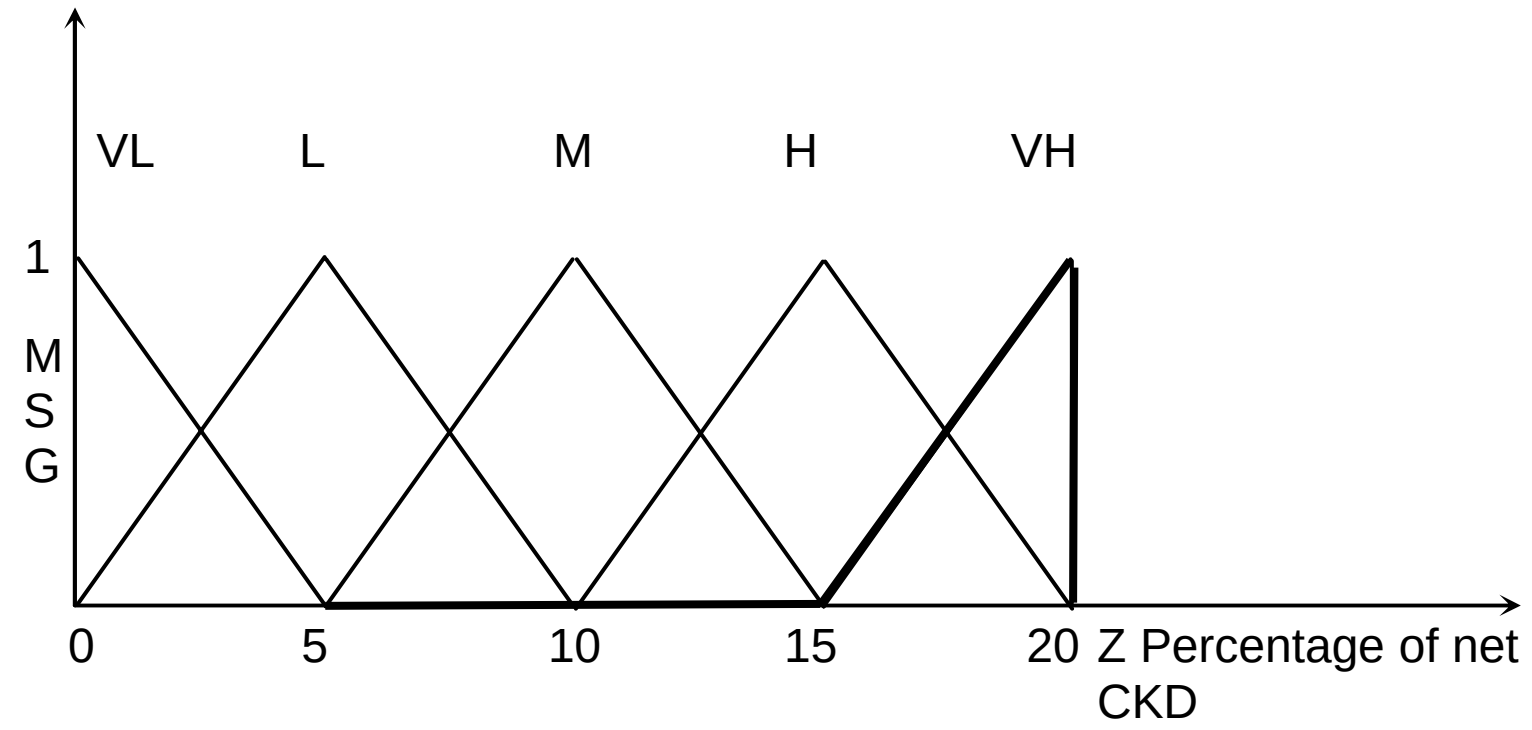

FIGURE 2.3.7: Aggregated output and defuzzification for the percentage of CKD

## 2.4 Finding the MIX of raw materials in proper proportion and minimize the waste dust using fuzzy neural network

The study of proper proportions of material mix during the clinkerization process is very difficult due to inconsistency in the chemical and mineralogical composition and the variation of these characteristic affects kiln operation, fuel consumption, clinker quality and above all the amount of CKD vent into the atmosphere. Further the raw mix should maintain a fixed range for a specific quality of cement. The problem of satisfying this range involves lot of randomness and uncertainty, which in turn speaks about the desired quality of the clinker. Chemical and mineralogical composition contains $SiO_2$, $Al_2O_3$, $Fe_2O_3$, CaO, MgO, $K_2O$ and $Na_2O$. Since all terms used to determine the proper proportions of material mix is very ambiguous, we felt it would be proper to use fuzzy theory approach to study the problem. We adopt fuzzy relational neural network method to find the correct proportion of raw mix so that the desired quality of the clinker is achieved. This is done by taking experts opinion about the proportions and then by giving fuzzy weights. This membership grades are varied a finite number of times till the error function reaches zero, which is equivalent to studying the set point values. The clinker of desired chemical composition is expected to satisfy the following modulus related to the chemical composition of the raw mix.
Lime saturation factor (LSF),

$$\text{LSF} = \frac{\text{CaO} \times 100}{2.8\ \text{SiO}_2 + 1.2\ \text{Al}_2\text{O}_3 + 0.65\text{Fe}_2\text{O}_3} \qquad (2.4.1)$$

A high LSF requires high heat consumption for clinker burning inside the kiln and this gives more strength to the cement.
Silica Modulus (SM)

$$\text{SM} = \frac{\text{SiO}_2}{\text{Al}_2\text{O}_3 + \text{Fe}_2\text{O}_3} \qquad (2.4.2)$$

A higher SM decreases the liquid phase content, which impairs the burnability of the clinker and reduces the cement setting time.

Alumina Modulus (AM)

$$AM = \frac{Al_2O_3}{Fe_2O_3} \qquad (2.4.3)$$

The value of AM determines the composition of liquid phase in the clinker.

Here we describe the working of the block schematic of raw mill processing. The raw mill grinder receives raw materials such as limestone, silica, iron and bauxite for the production of cement in separate feeders, called weigh feeders. All the raw materials are ground in a raw mill grinder to a powder form. A sample of this ground raw mix is collected at the output of the raw mill grinder by an auto sampler, and a sample is prepared after being finely ground by vibration mill and pressed by hydraulic press and then is analysed in the laboratory by an X-ray sequential spectrometer. The results of X-ray analysis, which are obtained through sampling and analysing the equipment, are fed to the computer through a data communication line, for the required control action. The entire process is illustrated in figure 2.4.1.

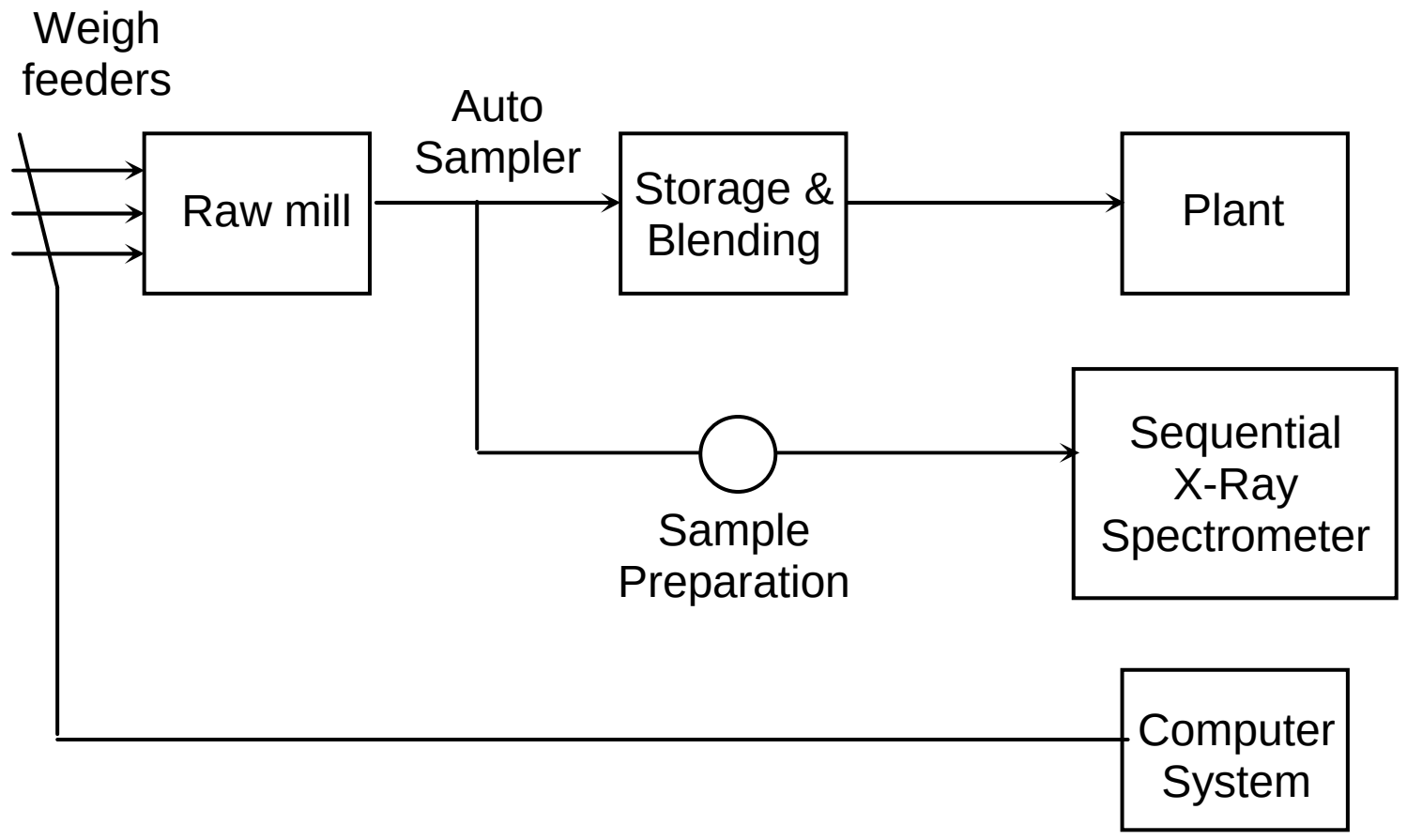

FIGURE 2.4.1 Block schematic of raw mill processing steps

The past researchers developed a control algorithm for raw mixing proportion based of singular value decomposition(SVD) methods. The purpose of this algorithm is to calculate the change in raw materials in each of the weigh feeders to achieve the raw mixing that is LSF, SM and AM.

Singular value decomposition(SVD) is one of the most basic and important tools in the analysis and solution of the problems in numerical linear algebra, and are finding increasing applications in control and digital signal processing. The potential of SVD technique is first exploited in the domain of linear algebra, where it provides a reliable determination of the rank of the matrix, thereby leading to accurate solutions of linear equations.

Here we adopt raw mix proportion control algorithm to our problem. The purpose of this algorithm is to calculate the change in raw materials in each of the weigh feeders to achieve the target value of the chemical composition ratio (or) module of LSF, SM and AM.

Suppose at any instant the action of the control system gives rise to the composition change as dLSF', dSM' and dAM' in response to the required composition change as dLSF, dSM and dAM respectively then the total mean square error at that instant will be

$$E = (dLSF - dLSF')^2 + (dSM - dSM')^2 + (dAM - dAM')^2 \quad (2.4.4)$$

The problem now is to minimize E with respect to the change in the feeder content(dw;: i = 1, 2, …, n). Differentiating equation (2.4.4) with respect to dw and equating to zero, we will have dLSF' = dLSF, dSM' = dSM and dAM' = dAM. As mentioned earlier, the values of LSF, SM and AM of the raw material, change constantly.

Our objective is to keep the values of LSF, SM and AM of the raw mix at the raw mill outlet fixed by changing the quantity of the raw material in the weigh feeders. So the module LSF, SM and AM are functions of the change in the raw material in different feeders. This can be represented as

$$dSLF' = dLSF = \sum_{i=1}^{n} \frac{\partial LSF}{\partial w_i} dw_i \qquad (2.4.5)$$

$$dSM' = dSM = \sum_{i=1}^{n} \frac{\partial SM}{\partial w_i} dw_i \qquad (2.4.6)$$

$$dAM' = dAM = \sum_{i=1}^{n} \frac{\partial AM}{\partial w_i} dw_i \qquad (2.4.7)$$

$$\sum_{i=1}^{n} dw_i = 0 \qquad (2.4.8)$$

$$LL_i \leq dw_i \leq HL_i \qquad (2.4.9)$$

where $w_i$ is the mix ratio of raw material in the feeder, $LL_i$ and $HL_i$ are the lower limit and the higher limit respectively of the raw material change possible for the $i^{th}$ feeder(i = 1, 2, …, n). The composition change, for example in LSF is given by

$$dLSF' = LSF_{sp} - LSF_{meas} \qquad (2.4.10)$$

Here 'sp' stands for set point that is the desired value and 'meas' stands for the measured value that is the value achieved.

Now consider the solution of equation (2.3.5) to (2.3.8). The number of unknowns is the same as the number of weigh feeders. If there are four unknowns then there are four weigh feeders, we have the following set of equations with four unknowns.

$$\frac{\partial LSF}{\partial w_1} dw_1 + \frac{\partial LSF}{\partial w_2} dw_2 + \frac{\partial LSF}{\partial w_3} dw_3 + \frac{\partial LSF}{\partial w_4} dw_4 = dLSF' \qquad (2.4.11)$$

$$\frac{\partial SM}{\partial w_1}dw_1 + \frac{\partial SM}{\partial w_2}dw_2 + \frac{\partial SM}{\partial w_3}dw_3 + \frac{\partial SM}{\partial w_4}dw_4 = dSM' \tag{2.4.12}$$

$$\frac{\partial AM}{\partial w_1}dw_1 + \frac{\partial AM}{\partial w_2}dw_2 + \frac{\partial AM}{\partial w_3}dw_3 + \frac{\partial AM}{\partial w_4}dw_4 = dAM' \tag{2.4.13}$$

$$dw_1 + dw_2 + dw_3 + dw_4 = 0 \tag{2.4.14}$$

Rearranging equations (2.4.11) to (2.4.14) in matrix form yields

$$\begin{bmatrix} \frac{\partial LSF}{\partial w_1} & \frac{\partial LSF}{\partial w_2} & \frac{\partial LSF}{\partial w_3} & \frac{\partial LSF}{\partial w_4} \\ \frac{\partial SM}{\partial w_1} & \frac{\partial SM}{\partial w_2} & \frac{\partial SM}{\partial w_3} & \frac{\partial SM}{\partial w_4} \\ \frac{\partial AM}{\partial w_1} & \frac{\partial AM}{\partial w_2} & \frac{\partial AM}{\partial w_3} & \frac{\partial AM}{\partial w_4} \\ 1 & 1 & 1 & 1 \end{bmatrix} \begin{bmatrix} dw_1 \\ dw_2 \\ dw_3 \\ dw_4 \end{bmatrix} = \begin{bmatrix} dLS\,F' \\ dS\,M'' \\ dA\,M' \\ 0 \end{bmatrix} \tag{2.4.15}$$

If the stacker reclaimed, a macline that feeds limestone of constant chemical composition to the weigh feeders is available, then LSF value will more(or) less remain constant; So in this case, one must give importance to achieving desired value for SM and AM. To cope with this situation in SVD method one can simply ignore equation(2.4.11). Also this method can be used in the event of feeder failure, or the addition of a feeder. In these cases, the number of feeder is simply changed and the corresponding equations, similar to equation (2.4.11) are added (or) deleted as appropriate.

The value is the amount of change for that modulus with unit change in raw material mix proportion sent into the grinder. This can be obtained from the calculation of the composition of the raw materials, but in cement production process the composition of the raw materials fed into the mill changes constantly. So it is not possible to get fixed values for these differential factors. Raw materials from a particular quarry have the composition varying over very narrow ranges for your purpose we have chosen a typical composition of raw material with its values as the average value of the material received from the quarry. The raw materials in each feeder consist of CaO, $SiO_2$, $Al_2O_3$, and $Fe_2O_3$ thus affecting all the three moduli such as LSF, SM and AM as given in equations (2.4.1, 2.4.2 and 2.4.3) so these moduli can now be redefined as

$$LSF = \frac{\sum_{i=1}^{n} CaO_1 \cdot w_i}{\sum_{i=1}^{n} \left[ 2.8\left(SiO_2\right) \cdot w_i + 1.2\left(Al_2O_3\right)_i \cdot w_i + 0.65\left(Fe_2O_3\right) \cdot w_i \right]} \tag{2.4.16}$$

$$SM = \frac{\sum_{i=1}^{n} \left(Si_2O_3\right)_i \cdot w_i}{\sum_{i=1}^{n} \left[ \left(Al_2O_3\right)_i \cdot w_i + \left(Fe_2O_3\right)_i \cdot w_i \right]} \tag{2.4.17}$$

$$AM = \frac{\sum_{i=1}^{n} \left(Al_2O_3\right)_i \cdot w_i}{\sum_{i=1}^{n} \left[ \left(Fe_2O_3\right)_i \cdot w_i \right]} \tag{2.4.18}$$

where n is the number of feeders. Now the differential coefficients of equation (2.4.11), (2.4.12) and (2.4.13) can be obtained by differentiating the equation (2.4.16), (2.4.17) and (2.4.18) with respect to $w_i$.

Adaptation of fuzzy neural network to raw mix proportion control algorithm :

Let P represent the coefficient of raw mixing ratio that is $\frac{\partial LSF}{\partial w_i}$, $\frac{\partial SM}{\partial w_i}$, $\frac{\partial AM}{\partial w_i}$ where i = 1, 2, 3, 4, Q represents the unknown quantities for four weigh feeders that is $dw_1$, $dw_2$, $dw_3$ and $dw_4$ and R represents the known values that is dLSF', dSM' and dAM'. Generally researchers used some other non-fuzzy method to estimate the unknowns $dw_1$, $dw_2$, $dw_3$ and $dw_4$ but since one is not always certain of solving these equations, fuzzy neural network model is adopted. By this method one is always guaranteed of a solution.

The problem is tackled in two stages according as if a solution exist using fuzzy relation equations $P \circ Q = R$ then all the quantities for four weigh feeders are determined. If $P \circ Q = R$ does not give solution the fuzzy neural network method is adapted to the fuzzy relation equation as the second stage. By adopting fuzzy neural network method to the fuzzy relation equation, unknown quantities for four weigh feeders are determined.

We get the matrices according to $P \circ Q = R$.

$$\begin{bmatrix} \frac{\partial LSF}{\partial w_1} & \frac{\partial LSF}{\partial w_2} & \frac{\partial LSF}{\partial w_3} & \frac{\partial LSF}{\partial w_4} \\ \frac{\partial SM}{\partial w_1} & \frac{\partial SM}{\partial w_2} & \frac{\partial SM}{\partial w_3} & \frac{\partial SM}{\partial w_4} \\ \frac{\partial AM}{\partial w_1} & \frac{\partial AM}{\partial w_2} & \frac{\partial AM}{\partial w_3} & \frac{\partial AM}{\partial w_4} \\ 1 & 1 & 1 & 1 \end{bmatrix} \begin{bmatrix} dw_1 \\ dw_2 \\ dw_3 \\ dw_4 \end{bmatrix} = \begin{bmatrix} dLSF' \\ dSM'' \\ dAM' \\ 0 \end{bmatrix} \qquad (2.4.19)$$

where

$$P = \begin{bmatrix} \frac{\partial LSF}{\partial w_1} & \frac{\partial LSF}{\partial w_2} & \frac{\partial LSF}{\partial w_3} & \frac{\partial LSF}{\partial w_4} \\ \frac{\partial SM}{\partial w_1} & \frac{\partial SM}{\partial w_2} & \frac{\partial SM}{\partial w_3} & \frac{\partial SM}{\partial w_4} \\ \frac{\partial AM}{\partial w_1} & \frac{\partial AM}{\partial w_2} & \frac{\partial AM}{\partial w_3} & \frac{\partial AM}{\partial w_4} \\ 1 & 1 & 1 & 1 \end{bmatrix}$$

$$Q = \begin{bmatrix} dw_1 \\ dw_2 \\ dw_3 \\ dw_4 \end{bmatrix} \text{ and } R = \begin{bmatrix} dLS\,F' \\ dS\,M' \\ dA\,M' \\ 0 \end{bmatrix}$$

we in this problem minimize the error between the rise to the composition change and required composition change. The membership value $p_{ij} \in [0,1]$ are given by experts.

Equation (2.4.19) can be rewritten as

$$\begin{bmatrix} p_{11} & p_{12} & p_{13} & p_{14} \\ p_{21} & p_{22} & p_{23} & p_{24} \\ p_{31} & p_{32} & p_{33} & p_{34} \\ p_{41} & p_{42} & p_{43} & p_{44} \end{bmatrix} \begin{bmatrix} dw_1 \\ dw_2 \\ dw_3 \\ dw_4 \end{bmatrix} = \begin{bmatrix} dLSF' \\ dSM' \\ dAM' \\ 0 \end{bmatrix}.$$

It can partitioned into 4 equations.

$$\begin{bmatrix} p_{11} & p_{12} & p_{13} & p_{14} \end{bmatrix} \begin{bmatrix} dw_1 \\ dw_2 \\ dw_3 \\ dw_4 \end{bmatrix} = \begin{bmatrix} dLS\,F' \\ dS\,M' \\ dA\,M' \\ 0 \end{bmatrix},$$

$$\begin{bmatrix} p_{21} & p_{22} & p_{23} & p_{24} \end{bmatrix} \begin{bmatrix} dw_1 \\ dw_2 \\ dw_3 \\ dw_4 \end{bmatrix} = \begin{bmatrix} dLS\,F' \\ dS\,M' \\ dA\,M' \\ 0 \end{bmatrix},$$

$$\begin{bmatrix} p_{31} & p_{32} & p_{33} & p_{34} \end{bmatrix} \begin{bmatrix} dw_1 \\ dw_2 \\ dw_3 \\ dw_4 \end{bmatrix} = \begin{bmatrix} dLSF' \\ dS\,M' \\ dA\,M' \\ 0 \end{bmatrix}$$

and

$$\begin{bmatrix} p_{41} & p_{42} & p_{43} & p_{44} \end{bmatrix} \begin{bmatrix} dw_1 \\ dw_2 \\ dw_3 \\ dw_4 \end{bmatrix} = \begin{bmatrix} dLSF' \\ dSM' \\ dAM' \\ 0 \end{bmatrix}.$$

If the above partitioned equation do not satisfy the condition max $q_{ik} < r_{ik}$ (where $q_{ik}$ are unknown quantities for weigh feeders and $r_{ik}$ are known values that is dLSF', dSM' and dAM') then the system of equations has final solution. If the above partitioned equation satisfy this condition max $q_{ik} < r_{ik}$ where $q_{ik}$ are unknown quantities for weigh feeders and $r_{ik}$ are known values that is dLSF', dSM' and dAM' then the system of equations has no solution. In this case fuzzy neural network method is adopted for fuzzy relation equation as the second stage.

$$P_1 \begin{bmatrix} dw_1 \\ dw_2 \\ dw_3 \\ dw_4 \end{bmatrix} = \begin{bmatrix} dLSF' \\ dSM' \\ dAM' \\ 0 \end{bmatrix}, \; P_2 \begin{bmatrix} dw_1 \\ dw_2 \\ dw_3 \\ dw_4 \end{bmatrix} = \begin{bmatrix} dLSF' \\ dSM' \\ dAM' \\ 0 \end{bmatrix},$$

$$P_3 \begin{bmatrix} dw_1 \\ dw_2 \\ dw_3 \\ dw_4 \end{bmatrix} = \begin{bmatrix} dLSF' \\ dSM' \\ dAM' \\ 0 \end{bmatrix}, \; P_4 \begin{bmatrix} dw_1 \\ dw_2 \\ dw_3 \\ dw_4 \end{bmatrix} = \begin{bmatrix} dLSF' \\ dSM' \\ dAM' \\ 0 \end{bmatrix}$$

The linear activation function f defined earlier gives the output $y_i=f[\max (w_{ij}\, x_j)]$, $i \in N_n$. First calculate $w_{11}x_1$, $w_{12}x_2$, $w_{13}x_3$, and $w_{14}x_4$ then find $y_1= f\,[\max(w_{1j}\, x_j)]$ which gives $dw_1$. Similarly calculate $w_{21}x_1$, $w_{22}x_2$, $w_{23}x_3$, and $w_{24}x_4$ to find $y_2=f\,[\max (w_{2j}\, x_j)]$ which gives $dw_2$, calculate $w_{31}x_1$, $w_{32}x_2$, $w_{33}x_3$, and $w_{34}x_4$ to find $y_3=f\,[\max (w_{3j}\, x_j)]$ which gives $dw_3$ and calculate $w_{41}x_1$, $w_{42}x_2$, $w_{43}x_3$, and $w_{44}x_4$ to find $y_4=f\,[\max (w_{4j}\, x_j)]$ which gives $dw_4$. Then we can find out using equation (2.4.4) that is whether the error function reaches 0 or not, suppose the error function does not reach 0, we change the weights that is the membership grades till the error reaches zero, which is explained in figure(2.4.8). Thus to achieve the value of error function to be zero, we give different membership grades to the weigh feeders (finite number of times) and make the value of required composition change to be equal to the raise in composition change.

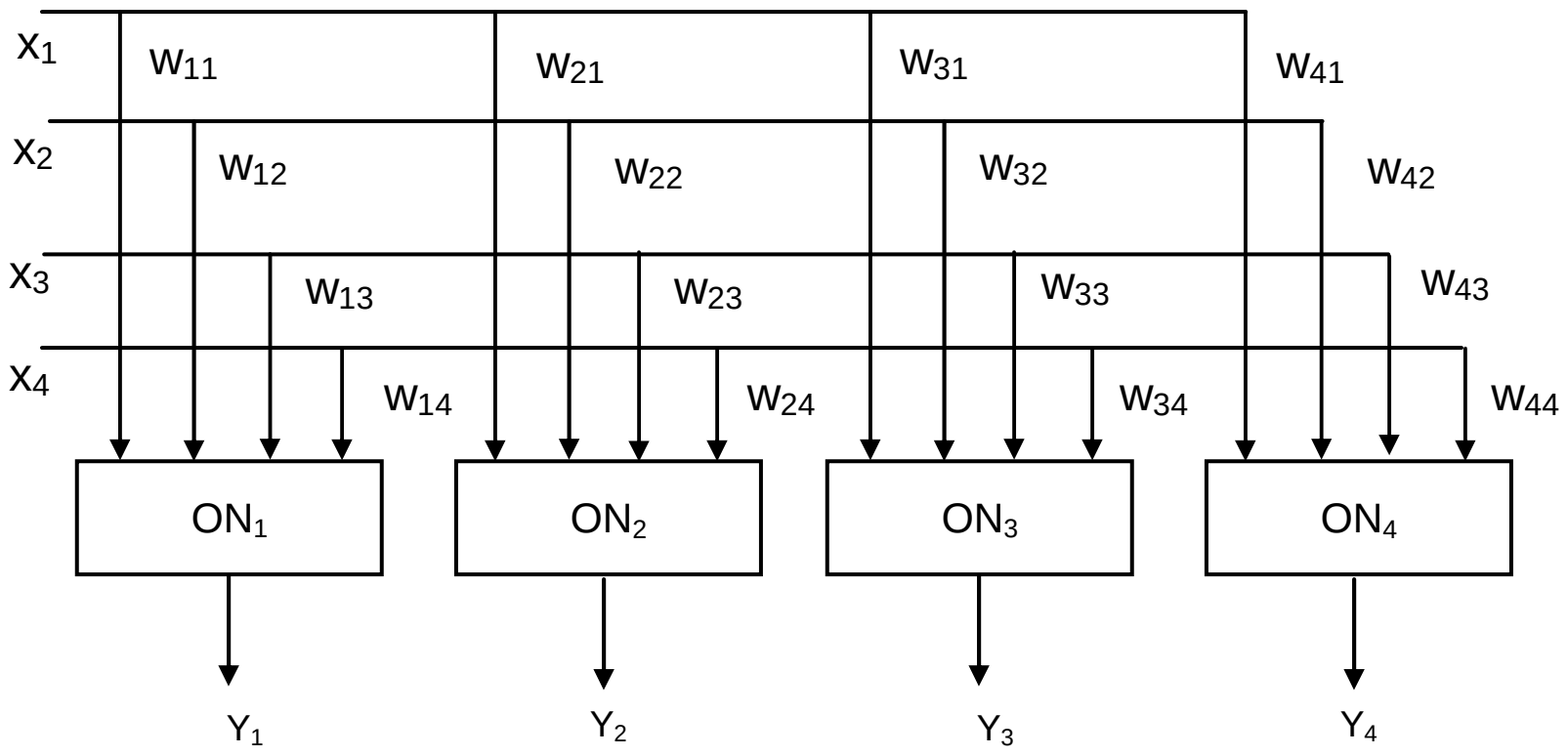

Figure 2.4.2: The feed forward neural network

## 2.5 Conclusions

The fuzzy control method described and defined in this chapter has the following problem:

1) Monitoring and control of the system
2) CKD reprocessing

***Monitoring and control :***

We have analyzed the alkali ratio (0.5 to 1.5) in kiln load material. The alkali ratio and the kiln load material in tons are considered as the two input parameter of fuzzy control. The output, is the percentage of CKD. The estimated results are as follows :

1) When the alkali ratio is 0.5 in the 5 tons of kiln load material, the estimated percentage of CKD is 0.
2) When the alkali ratio is 1 in the 15 tons of kiln load material, the estimated percentage of CKD is 20.
3) When the alkali ratio is 1.2 in the 17 tons of kiln load material, the estimated percentage of CKD is 20.
4) When the alkali ratio is 1.5 in the 25 tons of kiln load material, the estimated percentage of CKD is 40.

From our study we suggest in the online process to reduce or minimize the amount of CKD in the industry one should change the condition of fuel burning system and other system in kiln from time to time depending on the percentage of CKD in tons.

***CKD Reprocessing :***

CKD reprocessing, mainly concentrates on gas volume and temperature set point using fuzzy control. The fuzzy control method suggests the following results to minimize the CKD in reprocessing.

1) The suggested gas volume is 11865 $m^3$/min and temperature set point 350$^oC$ for reprocessing of CKD. At the time of reprocessing with suggested gas volume and temperature set point, the percentage of net CKD occurs from 0 to 5.
2) The suggested gas volume is 13020 $m^3$/min and temperature set point 400$^oC$ for reprocessing of CKD. At the time of reprocessing with suggested gas volume and temperature set point, the percentage of net CKD occurs from 10 to 15.
3) The suggested gas volume is 13080 $m^3$/min and temperature set point 430$^oC$ for reprocessing of CKD. At the time of reprocessing with suggested gas volume and temperature set point, the percentage of net CKD is 20.
4) The suggested gas volume is 15174 $m^3$/min and temperature set point 450$^oC$ for reprocessing of CKD. At the time of reprocessing with suggested gas volume and temperature set point, the percentage of net CKD is 20.

***Raw material mix using fuzzy network :***

The fuzzy neural network method defined has the following merits. The solution exists for all unknown weigh feeders and made the error between raise to the composition change and required composition change of the raw material close to zero. That means the change in raw materials in each of weigh feeders $dw_1$, $dw_2$, $dw_3$, $dw_4$ is achieved by the membership grade. This is very important one in cement industries to produce a desired quality of clinker.

***The merits of fuzzy control method :***

1) In earlier method the cement industry estimated the percentage of CKD approximately upto 40%. The industry did not know how much percentage of CKD occurs in each process. Using fuzzy control method, the estimated percentage of CKD in each process. By using this, the industry can change some internal condition of kiln and minimize the CKD in the online process.
2) The earlier methods adopted by the cement industries, choose the temperature set point and gas volume randomly for electrostatic precipitator to minimize the net CKD in reprocessing. But the random choice did not in general give the desired out comes. Using fuzzy control method, gives exact temperature set point and set point of gas volume from the range of temperature set point and gas volume for electrostatic precipitator. By using this temperature set point and gas volume set point, the industry will get the desired outcomes.

***The merits of the fuzzy neural networks :***

1) Solution exists to all unknown weigh feeders.
2) The target value of the chemical composition is achieved by minimizing the error between raise to the composition change and required composition change.

**Chapter Three**

# DETERMINATION OF TEMPERATURE SET POINTS FOR CRUDE OIL

## 3.1 Introduction

Study of Temperature set point in Chemical industries happen to be an important feature. Here we give an illustration how fuzzy control method is adopted for finding precise temperature set point to distil different crude in an oil refinery. Oil that comes from the ground is called the "crude oil". By cooking, the crude is converted to useful oil. Here the temperature set point plays a vital role at the time of cooking the crude oil. Since the quality and quantity of the crude is dependent on the temperature set point, the crude oil refinery has different temperature, set points to distil different crudes. Here, this chapter tries to determine a precise temperature set point for the crude oil refinery to maximise the distillation of the crude and the quantity of the crude for long hours. This chapter six sections.

This study is significant because most of the crude oil refineries have common type of operating systems. The analysis of this study is focussed on Kalundborg Refinery [Ebbesen (1992)]. Here we approach the problem of finding the precise temperature set point for different crudes using fuzzy control theory. The data is taken from Kalundborg Refinery [Ebbesen (1992)].

In 1995, Friedman developed a Mass and Enthalpy balance method and used it to improve the quality of crudes. In 1992, Ebbesen studied about the crude operating in Kalundborg Refinery. He made some derivations from the theory of Friedman. Finally he gave a range of temperature set points for the distillation of different crudes. However at the end of his study he made it clear that in the case of kerosene, 90% stayed within $1^0$C of its set point of temperature, in the case of naphtha 95% distillation stayed within $1^0$C of its set point of temperature. After two hours, the quality during crude switches was different indicating a lower quality.

Here we establish the result using the data taken from Ebbesen(1992). This gives temperature set point for the distillation of kerosens, naphtha and gasoil. The range of temperature set points and the various percentage of distillation are converted into the fuzzy control theory. Here membership grade is assigned to each temperature set point and percentage of distillation and fuzzy control rule is used to each temperature set point and percentage of distillation. Finally center max-min rule is used, to find the precise temperature set point for kerosene, naphtha and gasoil.

## 3.2 Description of Crude Oil Refineries

Using the data available from Ebbesen(1992) of the Kalundborg oil refinery, we analyse the data via fuzzy rules and membership grades of fuzzy control theory method and find the precise temperature set points for different crudes to maximize the quality and distillation of crude for long hours. Crude oil refinery selects temperature set points randomly from the range of temperature set points for the distillation of the crude; as a result, the quality and quantity of the processed crude are maintained only for very few hours. Thus, to be more precise the aim to find the precise temperature set points for kerosene, naphtha and gasoil using the data of Ebbesen (1992) from the Kalundborg oil refinery.

Kalundborg oil refinery operates with different crudes on a regular basis. Here this crude oil refinery distils kerosene, naphtha and gasoil. A schematic diagram of various streams is shown in the following figure [Ebbessen(1992)].

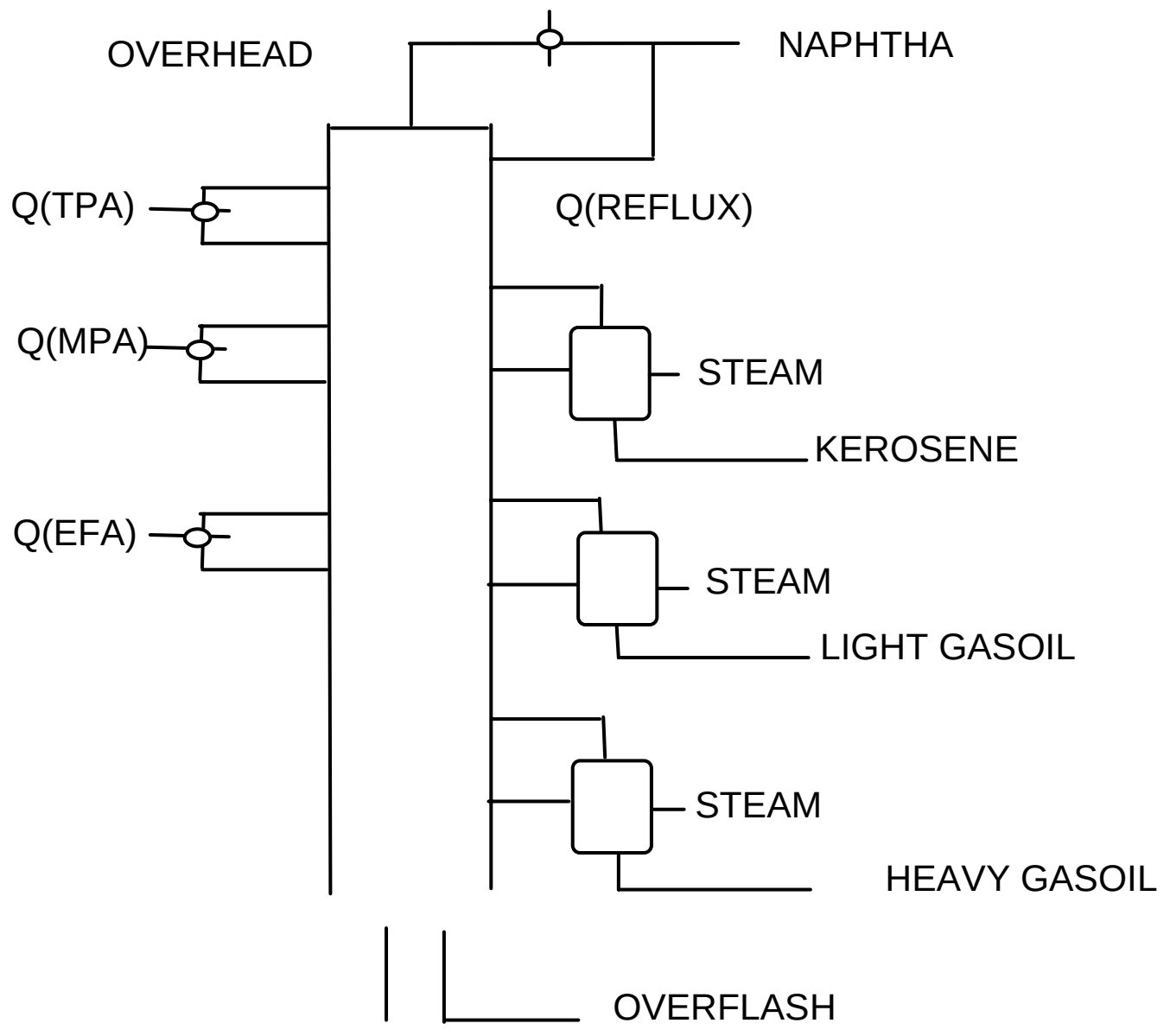

FIGURE 3.2.1: CRUDE REFINERY

where TPA – denotes the top pump-around, MPA – denotes mid pump-around and the BPA-bottom pump-around respectively. Q-denotes the heat removed. These are mainly used for controlling the temperature.

The random choice of temperature set points for different crudes with distillation taken from Ebbesen(1992) are described for kerosene, naphtha and gasoil.

### Range of temperature set points for kerosene with the percentage of distillation

The crude oil refinery gives $227^0$ C-temperature set point for distillation of kerosene. The crude oil refinery selects randomly this $227^0$ C temperature set point from the given range of set points {$220^0$ C, $221^0$ C, $222^0$ C, $223^0$ C, $224^0$ C, $225^0$ C, $226^0$ C, $227^0$ C, $228^0$ C, $229^0$ C, $230^0$ C}. This temperature set point $227^0$ C gives 90% distillation and it says within $1^0$ C of its set point of temperature. The temperature graph is given below.

Graph 3.2.1: Graph depicting the 90% distillation of Kerosene

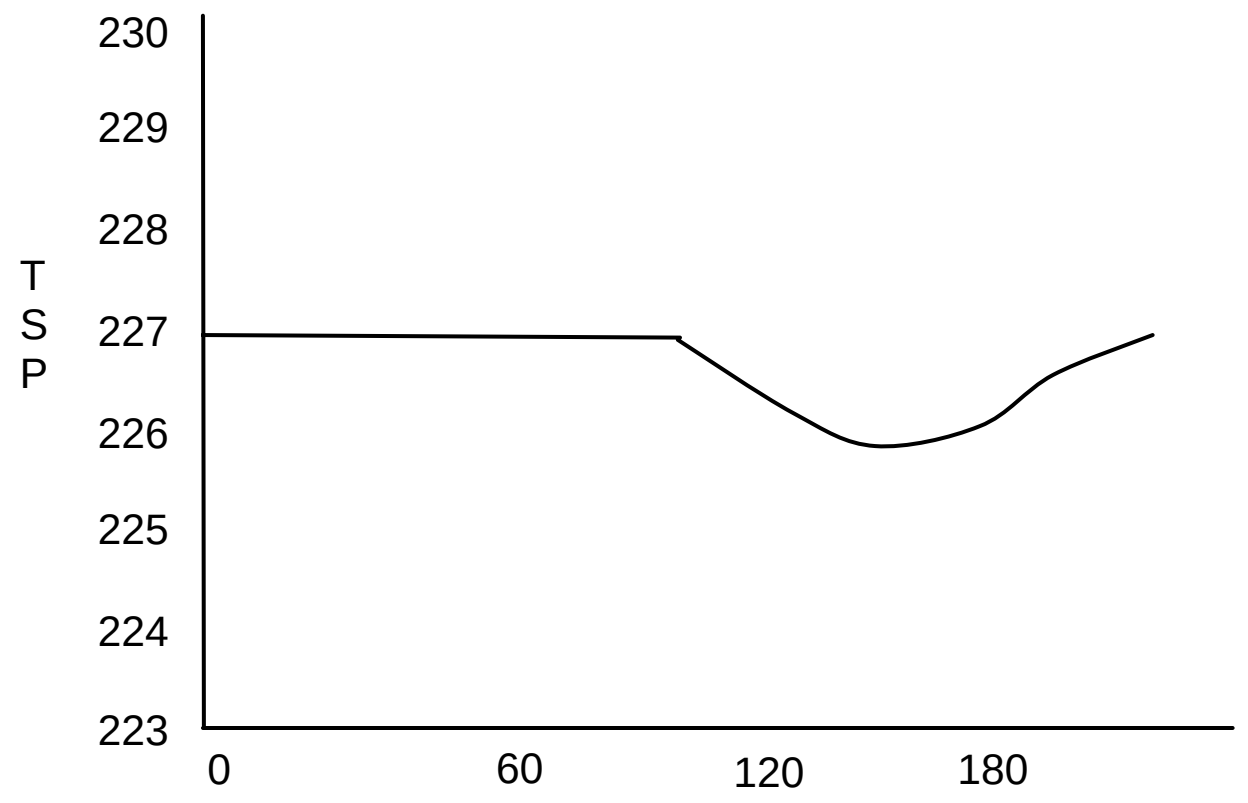

X axis– time in minute and Y axis set point of temperature for Kerosene. Set point was $227^oC$

Legend
TSP: Temperature set points

From this Ebbesen (1992) concludes that for the set point $227^0$ C the distillation of kerosene was 90%.

### Range of temperature set point of naphtha with percentage of distillation

The crude oil refinery gives $160^0$ C-temperature set point for the distillation of naphtha. The crude oil refinery selects randomly this $160^o$ C temperature set point from the given range of set points {$155^o$ C, $156^o$ C, $157^o$ C, $158^o$ C, $160^o$ C, $161^o$ C, $162^o$ C, $163^o$ C, $164^o$ C, $165^o$ C}. This temperature set point $160^o$ C gives 95% distillation and it stays within $1^o$ C of its set point of temperature.

The temperature graph is given below.

Graph 3.2.2: Graph depicting the 95% distillation of Naphtha

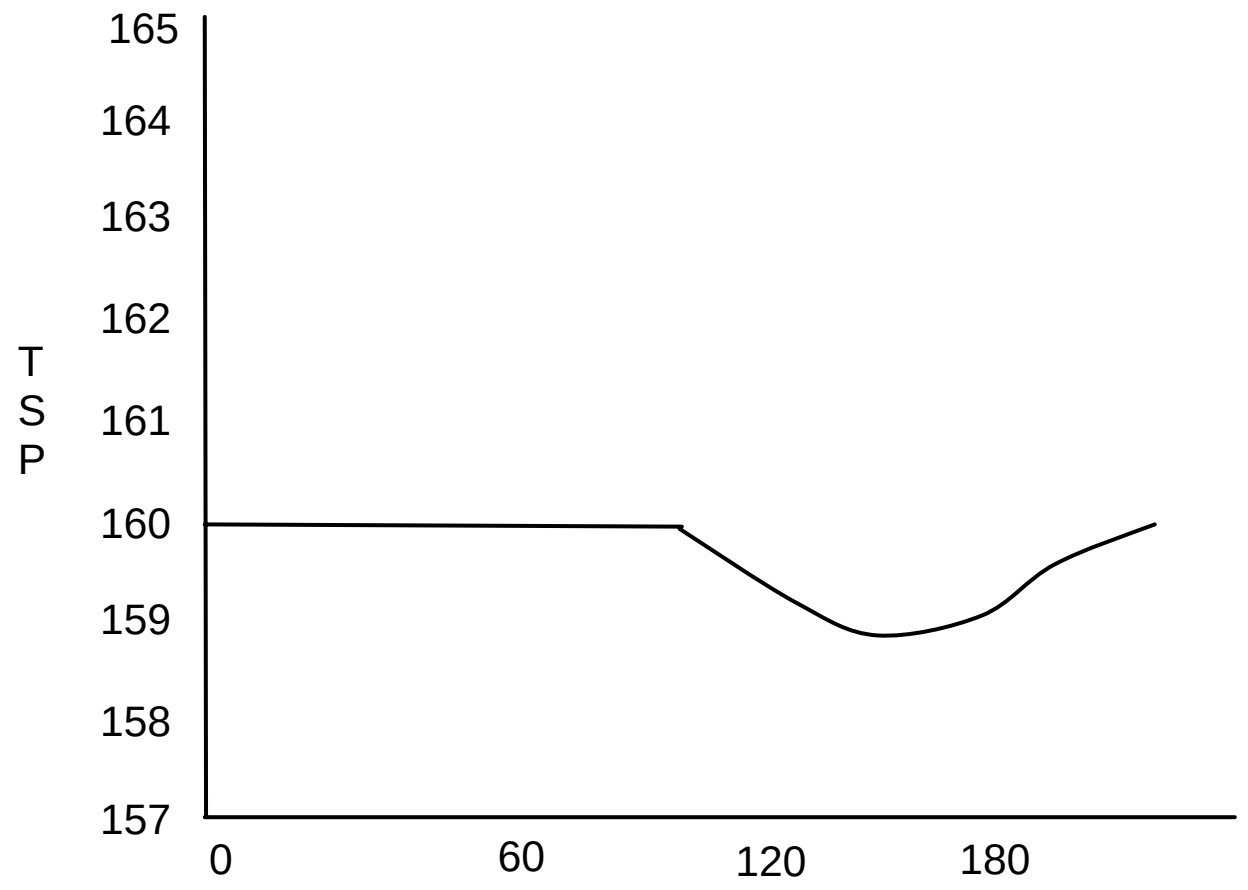

X axis– time in minute and Y axis set point of temperature for naphtha. Set point was $160^oC$

From this Ebbesen(1992) concludes that for the set point $160^0$ C the distillation of naphtha was 95%.

**Range of temperature set point of gasoil with the percentage of distillation**

The temperature set point $-4.5^0$ C gives 95% distillation and it stays within $1^0$ C of its set point temperature.

The temperature graph is given below.

Graph 3.2.3: Graph depicting the 95% distillation of gasoil

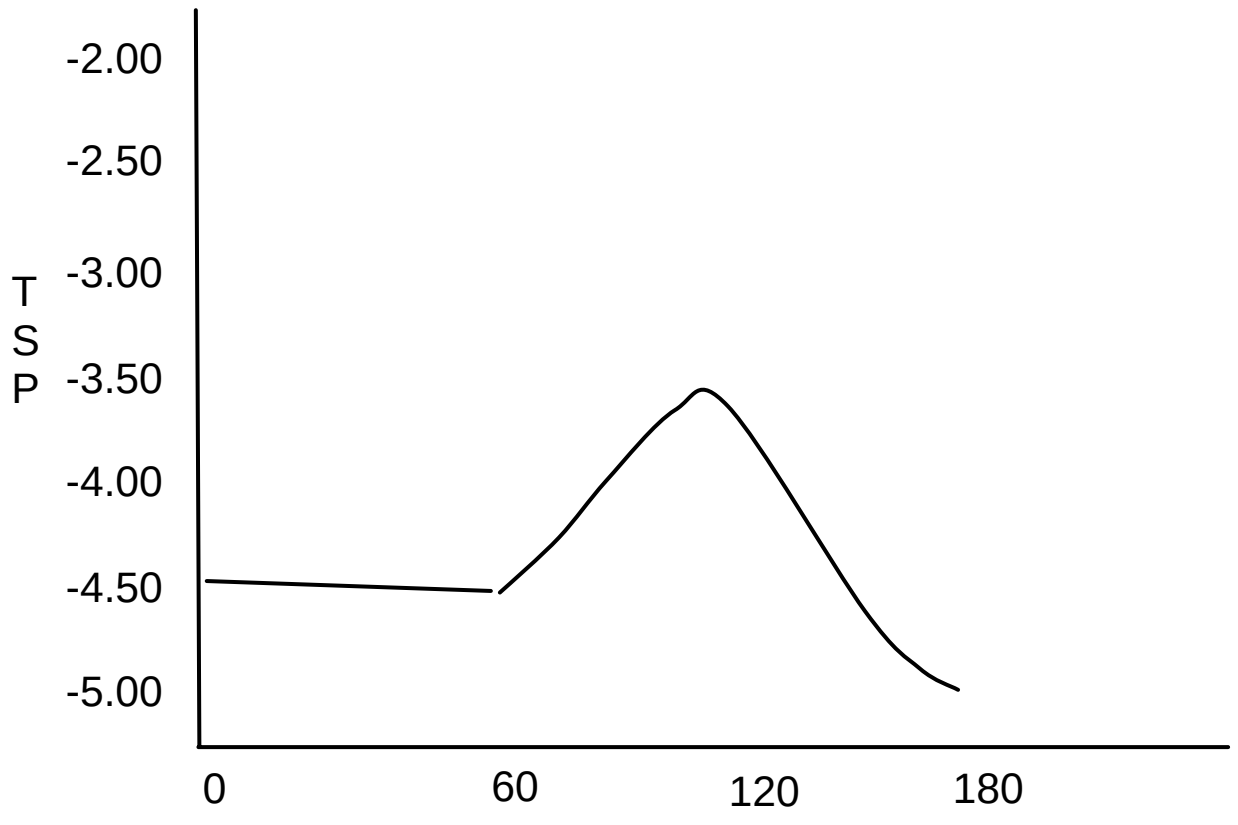

X axis– time in minute and Y axis set point of temperature for naphtha. Set point was $-4.50^{o}C$

From this Ebbesen(1992) concludes that for the set point $-4.5^{0}$ C the distillation of gasoil was 95%.

### 3.3 Determination of Temperature Set-Point of Kerosene Resulting in Better Distillation Using Fuzzy Control Theory

The given possible ranges of temperature set points are {$220^{0}$ C, $221^{0}$ C, $222^{0}$ C, $223^{0}$ C, $224^{0}$ C, $225^{0}$ C, $226^{0}$ C, $227^{0}$ C, $228^{0}$ C, $229^{0}$ C, $230^{0}$ C} and possible percentages of distillation are (88%, 89%, 90%, 91%, 92%} in case of kerosene as observed by Ebbesen(1992). As fuzzy control theory is the tool adaptable only when the past performance data is available, now this chapter considers the given possible range of temperature set points and distillation as the inputs of fuzzy control theory. To identify the precise temperature set points from the possible range of temperature set points, this chapter assigns membership grades to each input of fuzzy control theory. Here, the fuzzy control theory is used to find a precise temperature set point for kerosene.

In the procedure developed here membership grade is from the interval [0, 1] to the input of each temperature set point and each percentage of distillation. After membership grades are assigned to each input of temperature set points, the following graph results representing the membership grades of temperatures set point.

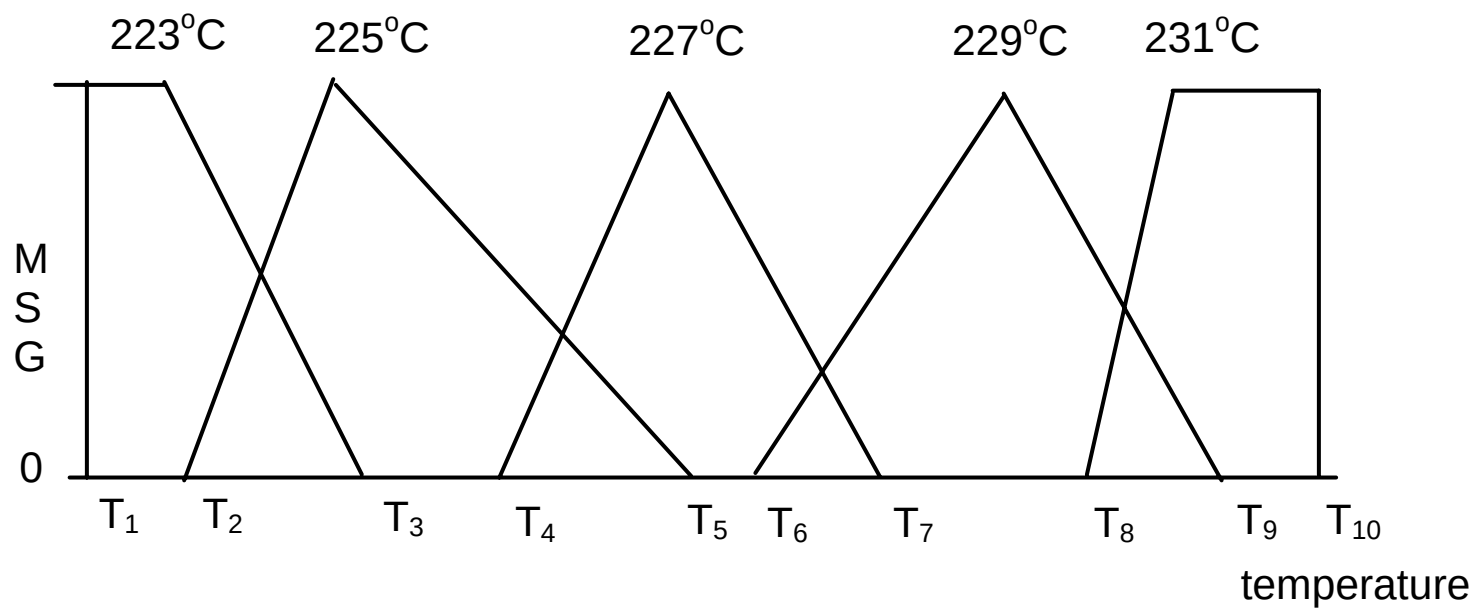

FIGURE 3.3.1: Membership function of temperature set point

After membership grades are assigned to each input of distillation, the following graph results representing the membership grades of distillation.

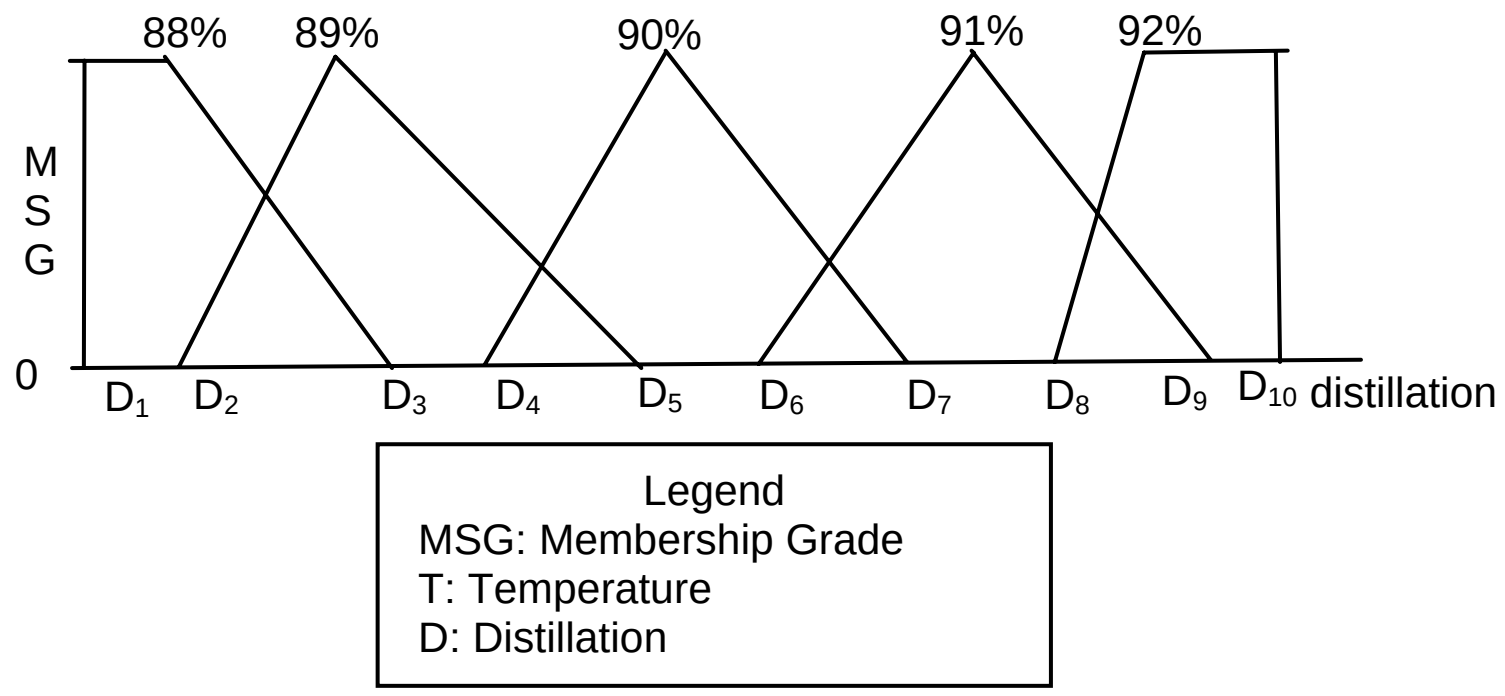

FIGURE 3.3.2: Membership function of percentage of distillation

The membership grade varies from 0 to 1.

For getting precise temperature set point for kerosene the throttle variables(The grade of membership) are qualified into five subsets. Here fuzzy rules are used to find the possible percentage of distillation for each temperature set point and the Center Max-Min rule is used to find a throttle membership grade for the existing fuzzy rules. To get the grade of membership to each existing fuzzy rule throttle variables are qualified into five subsets as follows:

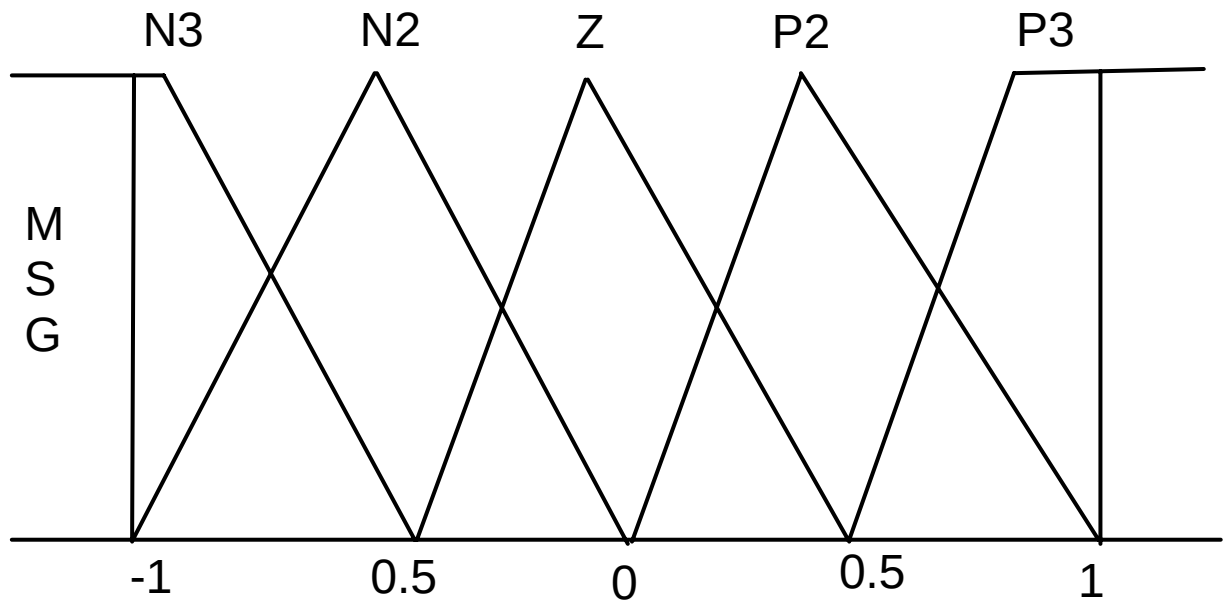

N3: Very Big Negative
N2: Big Negative
Z : Normal
P2: Big Positive
P3: very Big Positive

FIGURE 3.3.3: Throttle values

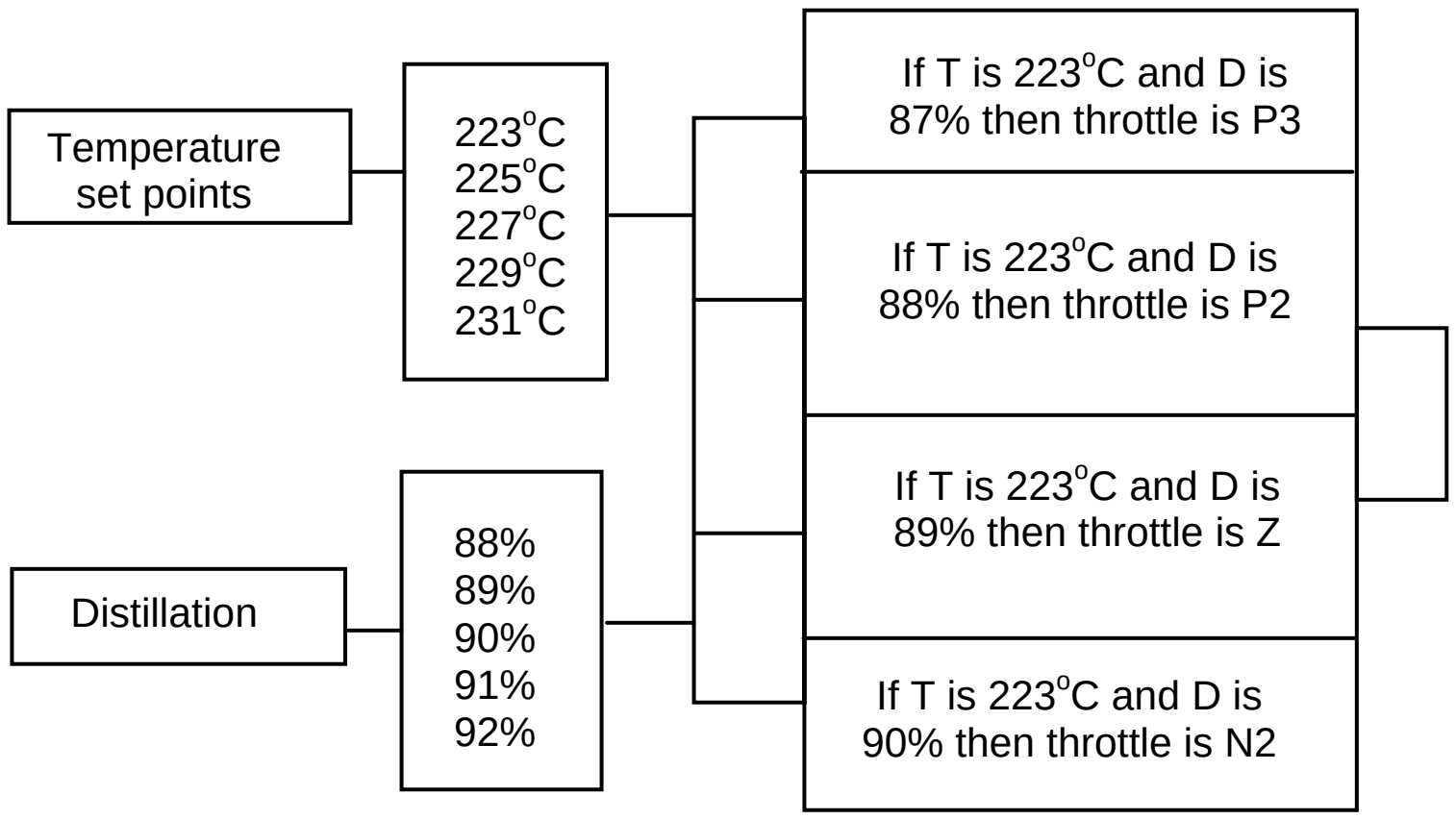

T: Temperature set point D: Percentage of Distillation

FIGURE 3.3.4:Fuzzy rules for the temperature set point 223°C

The following are the fuzzy rules :

Rule -1 : If T is $223^0$ C temperature set point AND D is 87% THEN throttle is P3.
Rule -2 : If T is $223^0$ C temperature set point AND D is 88% THEN throttle is P2.
Rule -3 : If T is $223^0$ C temperature set point AND D is 89% THEN throttle is Z.
Rule -4 : If T is $223^0$ C temperature set point AND D is 90% THEN throttle is N2.

We conclude the throttle value to the temperature set point for $223^0$ C, by the above stated rules, only rule -2, and rule-3 are applicable that is distillation is 88% and 89% respectively.

**Rule-2**
The throttle value to the temperature set point $223^0$ C for 88% distillation is calculated using figures 4.2 and 4.3.

$$\text{Throttle} = (0.41+0.38)/2 = 0.395$$

The graphical representation of the membership grade of the temperature set point $223^0$ C for 88 percentage of distillation is as follows.

Graph 3.3.1: Graphical representation of Rule 2

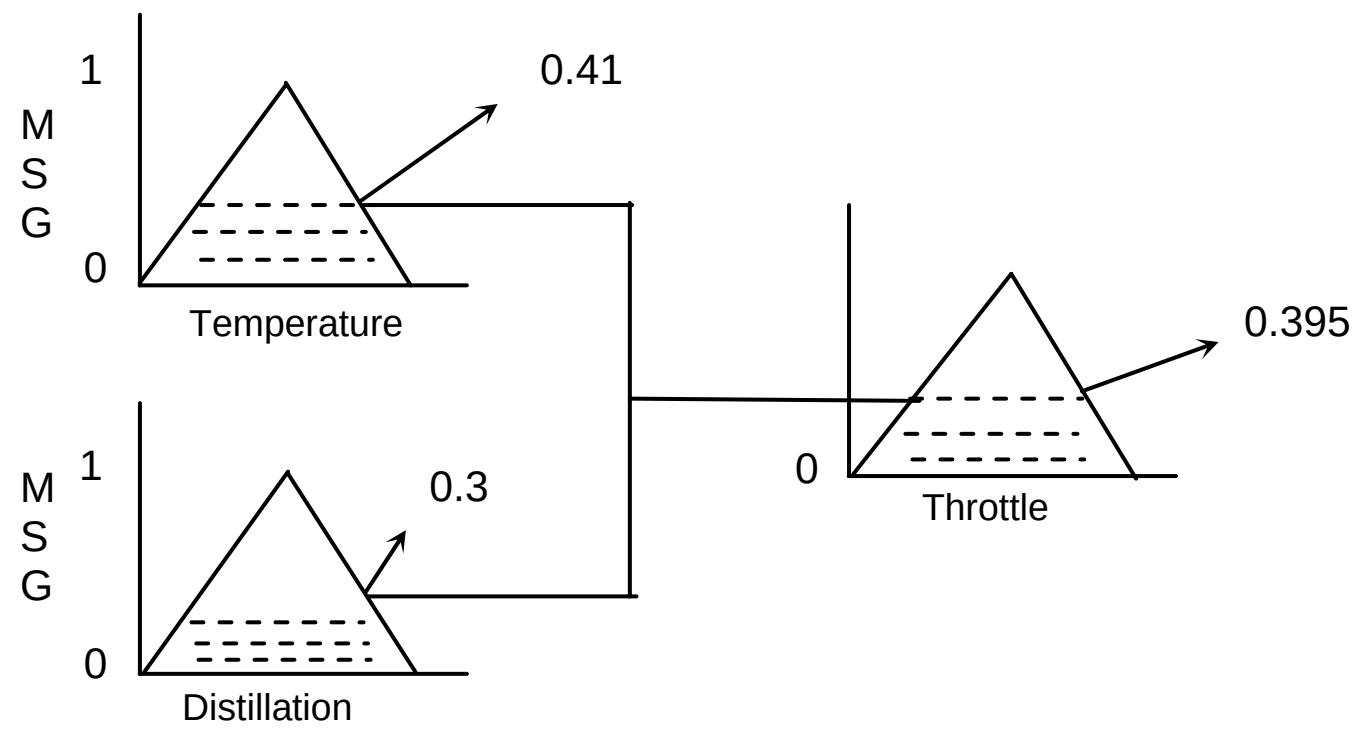

Graph 3.3.2: The two outputs are then defuzzified
by center max-min rule

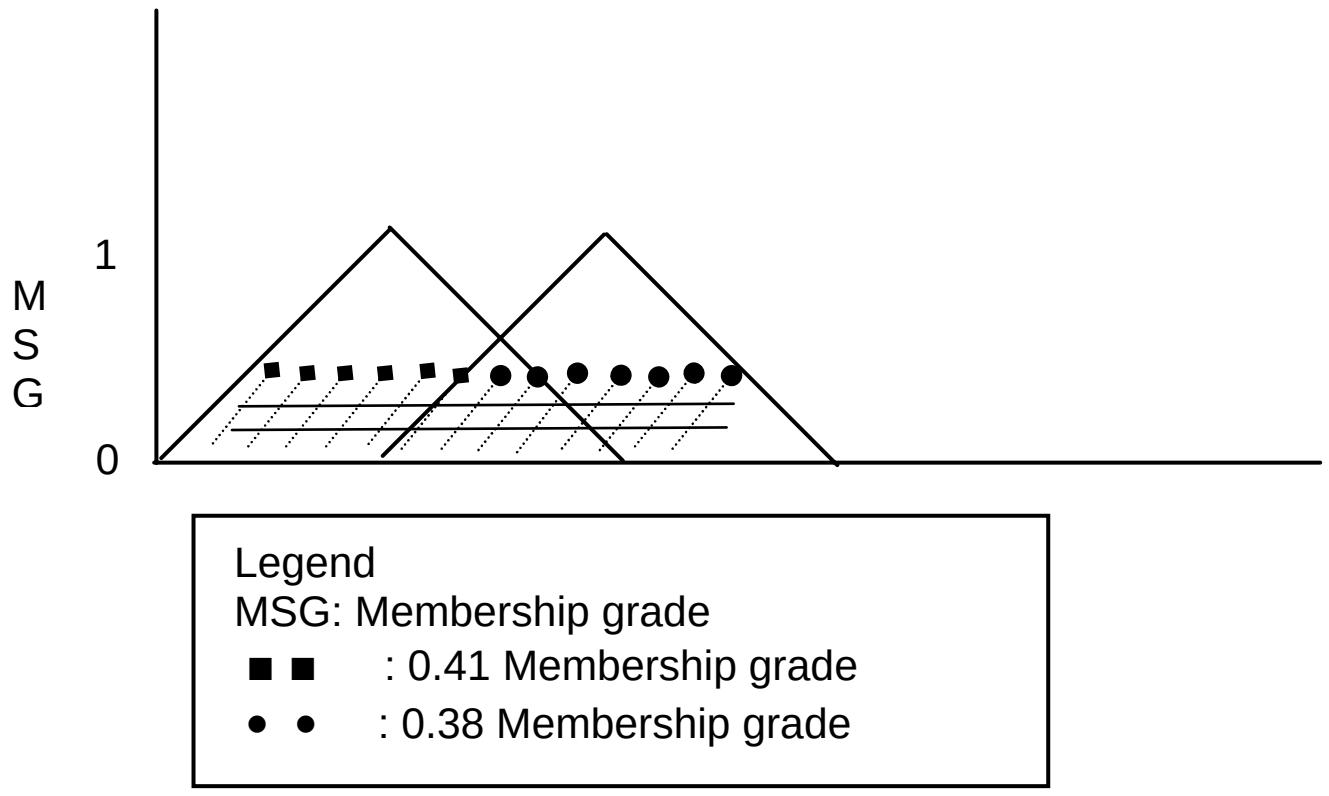

**Rule -3**

The throttle value to the temperature set point $223^0$ C for 89% distillation from figures 4.2 and 4.3 is as follows

Throttle=(0.47+.45)/2=0.46

The graphical representation of the membership grade of the temperature set point $223^0$ C for 88% percentage of distillation is as follows:

Graph 3.3.3: Graphical representation of Rule 3

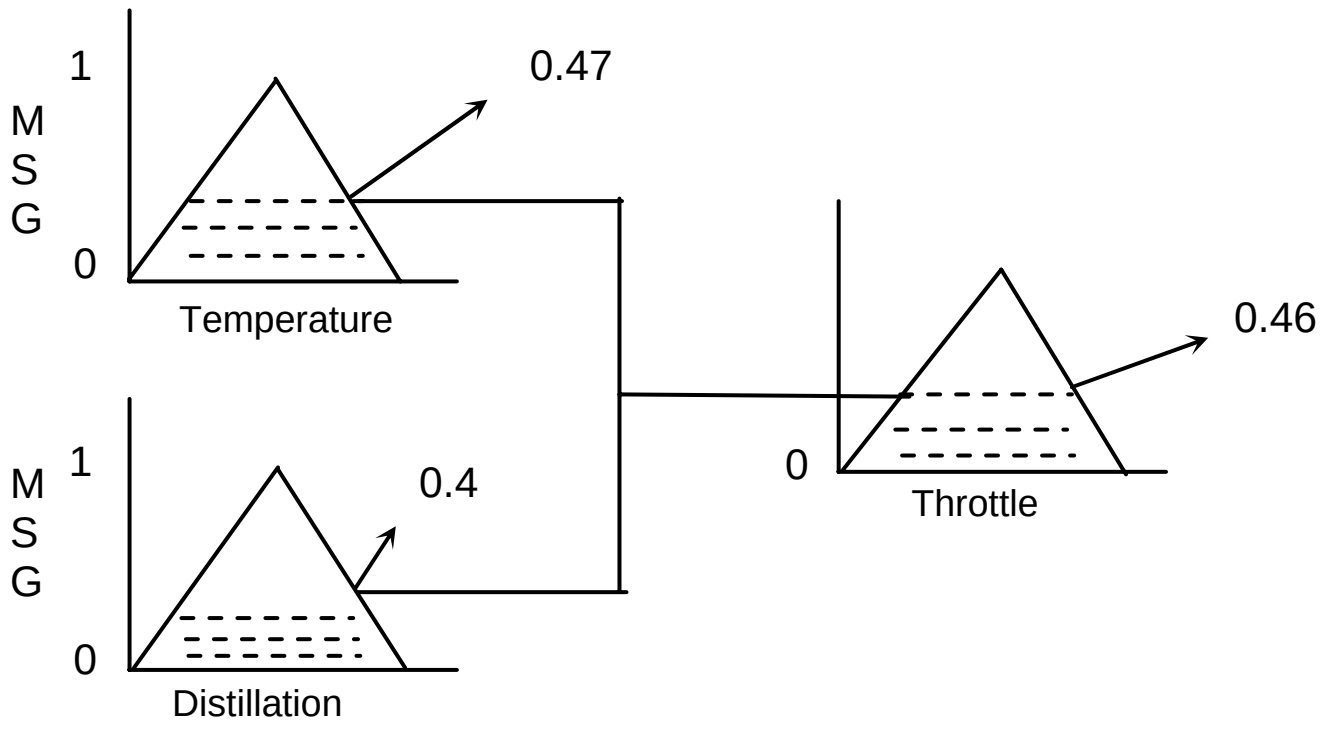

Graph 3.3.4: The two outputs are then defuzzified by center max-min rule

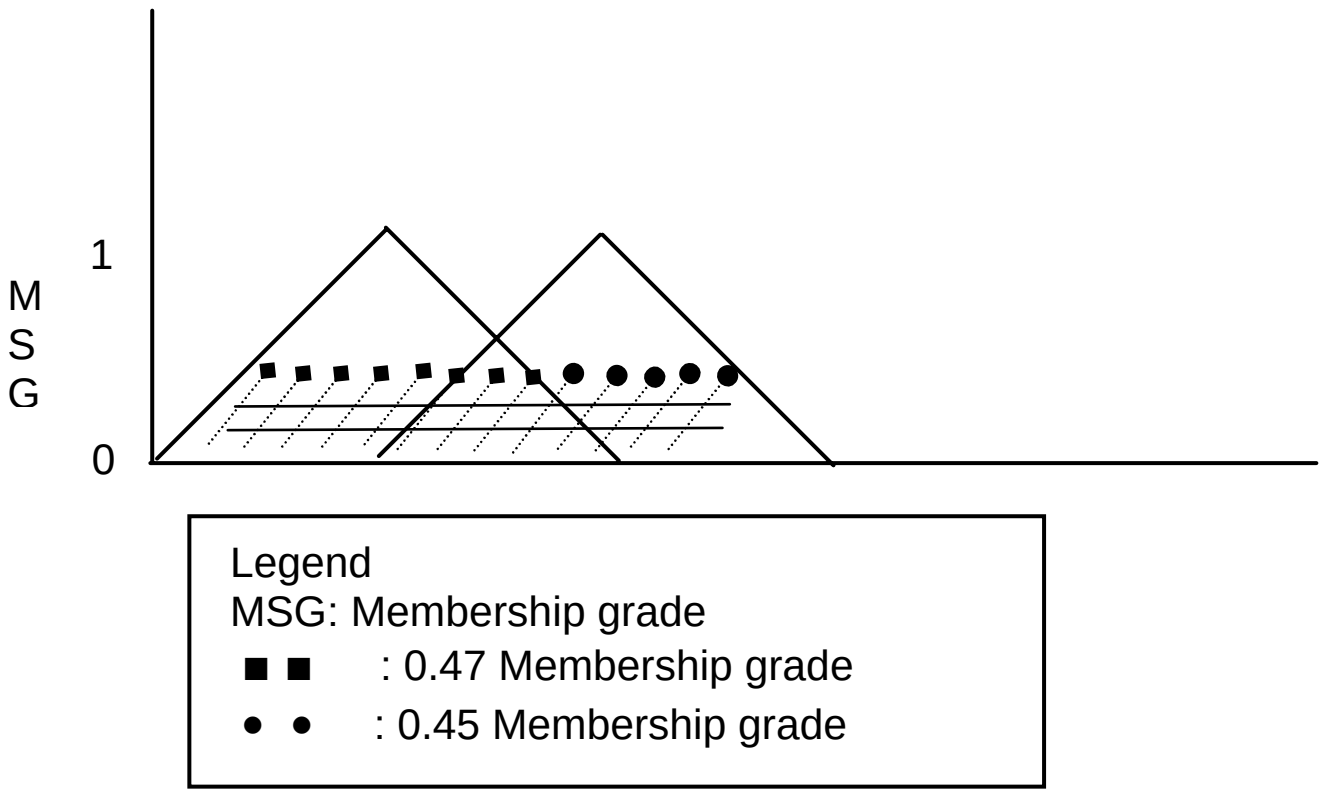

Here, the Center Max-Min rule is used to find a precise temperature set point.

**Using Center Max-Min rule to find precise temperature set point for kerosene**

Throttle(grade of membership) = m(P3) × Location(P2) + m(Z) + Location(N2) = 0.427. Graphs for the other rules have not been given explicitly but after calculations, the values are given as the same procedure is adopted.

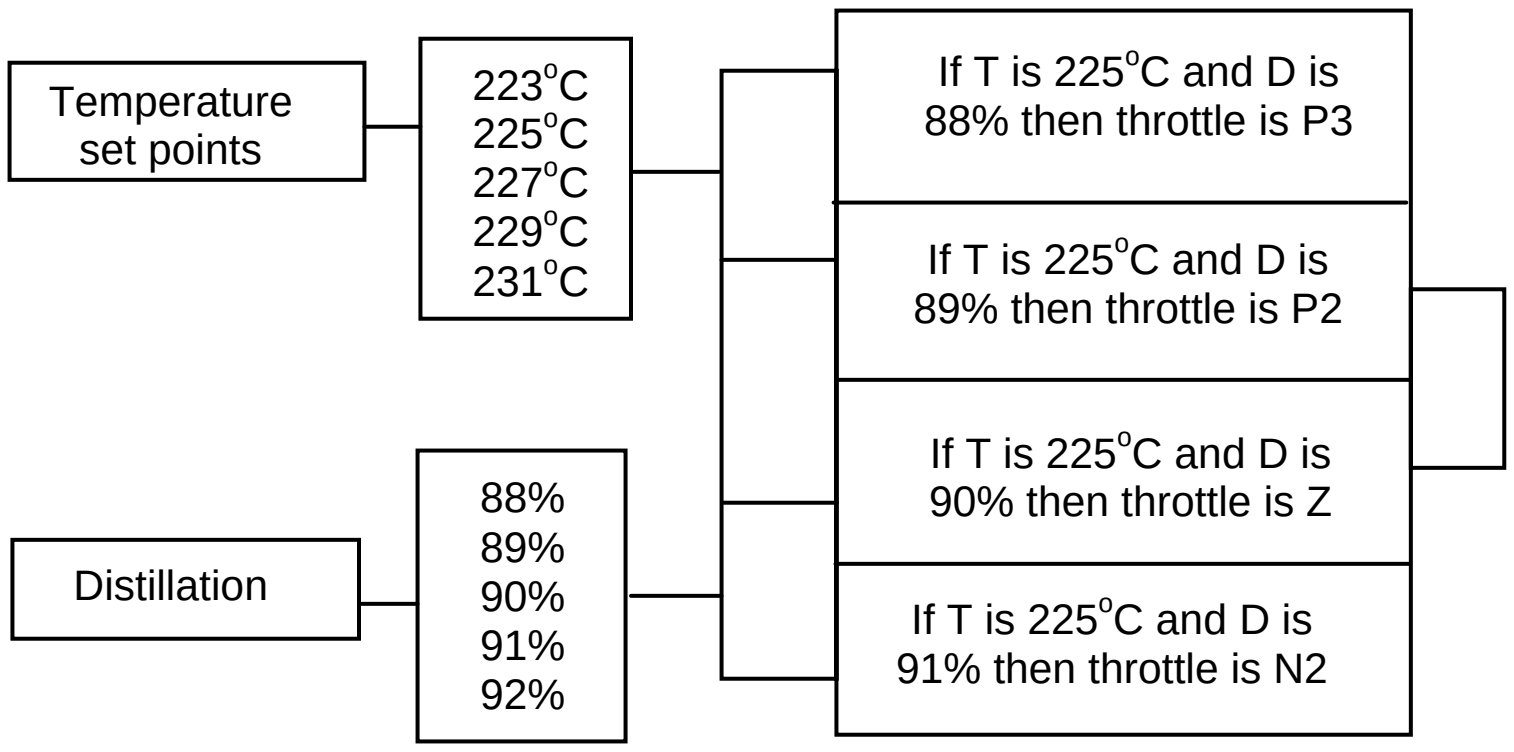

T: Temperature D: Distillation

FIGURE 3.3.5:Fuzzy rules for the temperature set point 225$^o$C

The fuzzy rule for distillation of kerosene:

Rule -1 : If T is $225^0$C temperature set point AND D is 88% THEN throttle is P3.
Rule -2 : If T is $225^0$C temperature set point AND D is 89% THEN throttle is P2.
Rule - 3 : If T is $225^0$C temperature set point AND D is 90% THEN throttle is Z.
Rule - 4 : If T is $225^0$C temperature set point AND D is 91% THEN throttle is N2.

**Rule -2**

The throttle value to the temperature set point $225^0$C for 89% distillation is calculated as follows.

Throttle= (.47+.45)/2=0.46,

**Rule -3**

The throttle value to the temperature set point $225^0$C for 90% distillation is calculated as follows.

Throttle= (.51+.23)/2=0.37,

**Using Center Max-Min rule to find precise temperature set point for kerosene**

Throttle(grade of membership) = m(P3) × Location(P2) + m(Z) + Location(N2)
= 0.415,

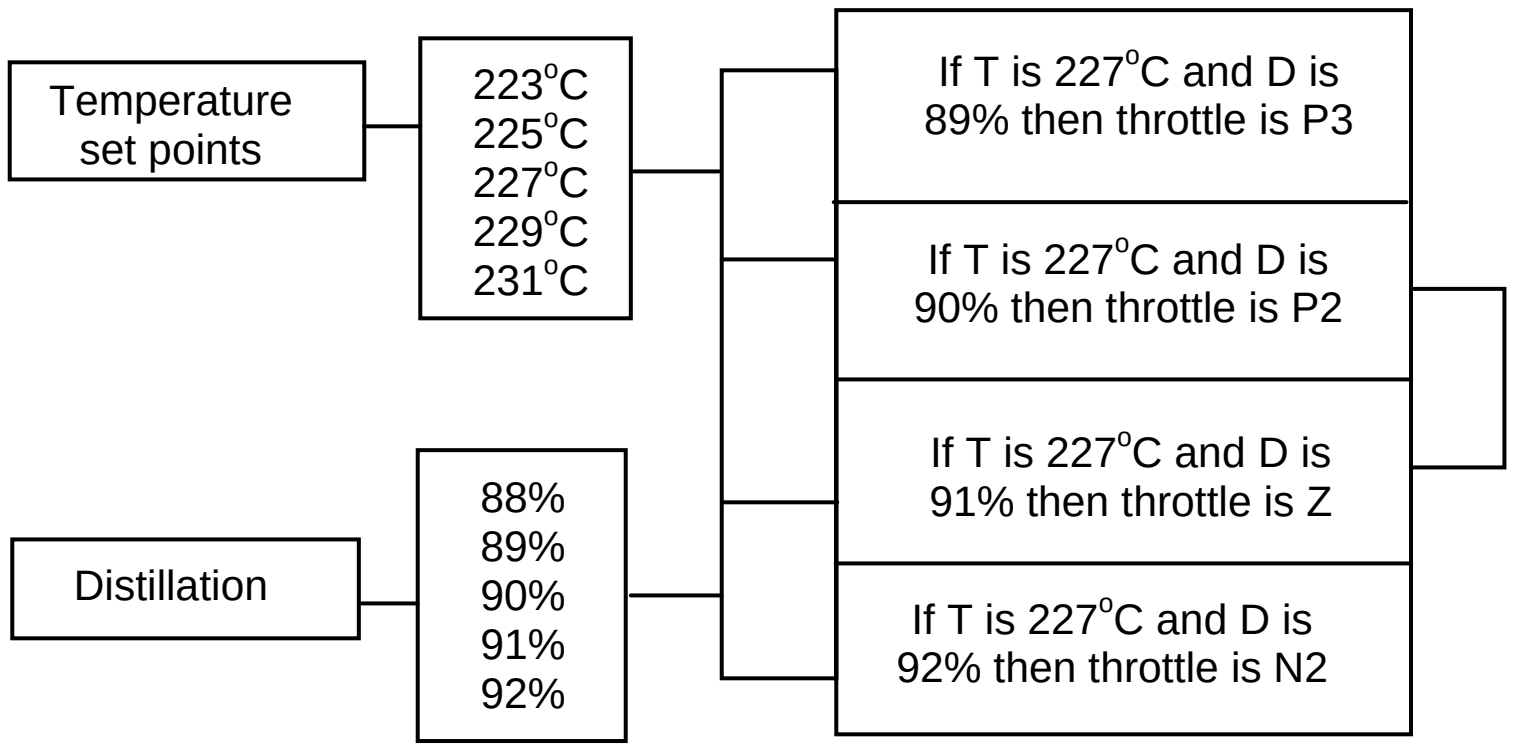

FIGURE 3.3.6:Fuzzy rules for the temperature set point $227^{o}$C

The fuzzy rule for distillation of kerosene:

Rule -1 : If T is $227^0$C temperature set point AND D is 89% THEN throttle is P3.
Rule -2 : If T is $227^0$C temperature set point AND D is 90% THEN throttle is P2.
Rule - 3 : If T is $227^0$C temperature set point AND D is 91% THEN throttle is Z.
Rule - 4 : If T is $227^0$C temperature set point AND D is 92% THEN throttle is N2.

Here we calculate the throttle value to the temperature set point for $227^0$ C, by the above stated rules, only rule -2, and rule-3 are applicable that is distillation is 90% and 91% respectively.

**Rule-2**

The throttle value to the temperature set point $227^0$ C for 90% distillation is calculated as follows:

Throttle = (0.51+0.23)/2
= 0.37,

**Rule-3**

The throttle value to the temperature set point $227^0$ C for 91% distillation is calculated as follows.

Throttle = (0.17+0.59)/2
= 0.38,

**Using Center Max-Min rule to find precise temperature set point for kerosene**

Throttle(grade of membership) = m(P3) x Location(P2) + m(Z) + Location(N2)
= 0.375,

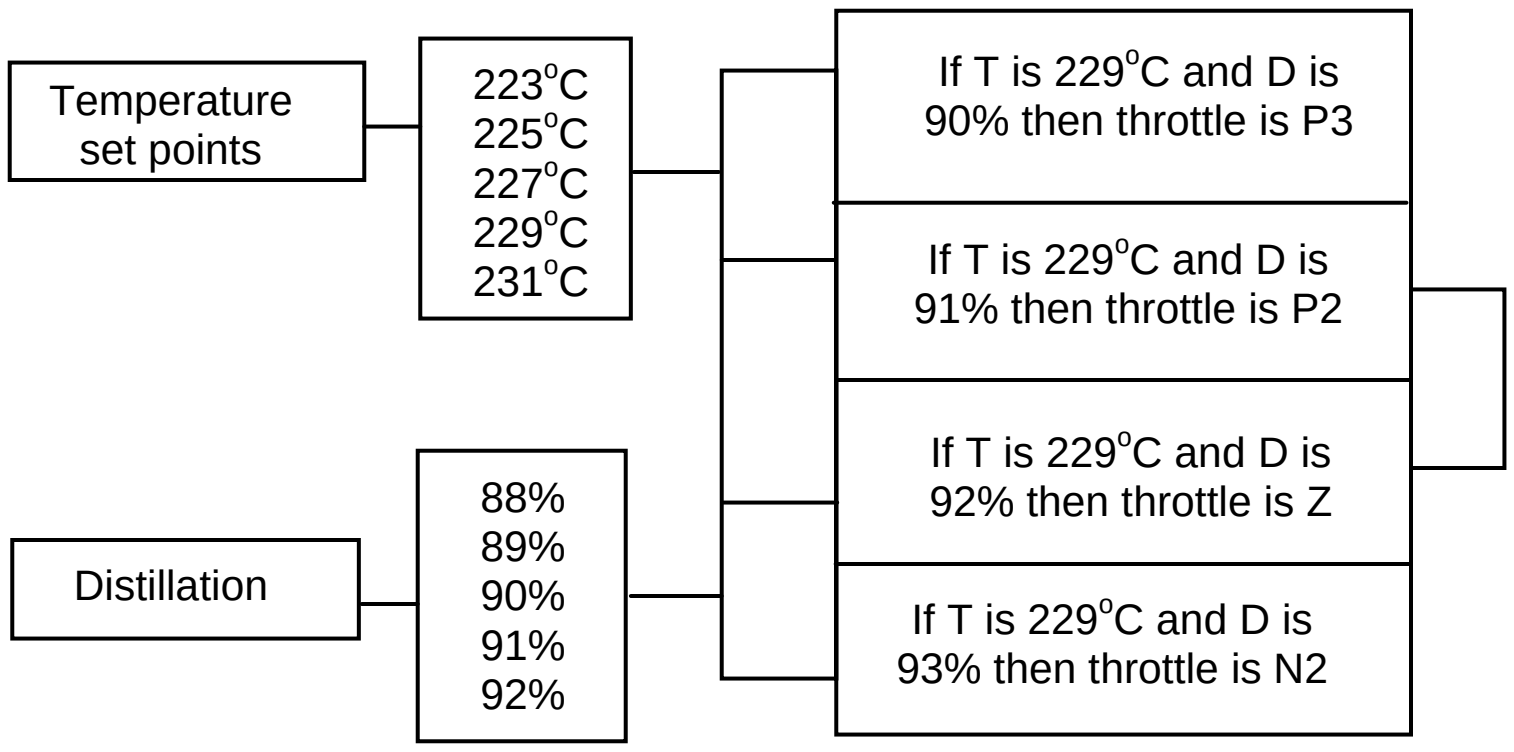

FIGURE 3.3.7:Fuzzy rules for the temperature set point $229^{o}$C

Here we obtain the throttle value to the temperature set point for $229^{0}$ C, by the above stated rules, only rule -2, and rule-3 are applicable that is distillation is 91% and 92% respectively.

**Rule-2**

The throttle value to the temperature set point $229^{0}$ C for 91% distillation is calculated as follows.

$$\text{Throttle} = (0.\ 17 + 0.59)/2 = 0.38,$$

**Rule-3**

The throttle value to the temperature set point $229^{0}$ C for 92% distillation is calculated as follows.

$$\text{Throttle} = (0.4 + 0.4)/2 = 0.4,$$

**Using Center Max-Min rule to find precise temperature set point for kerosene**

$$\begin{aligned}\text{Throttle(grade of membership)} &= m(P3) \times \text{Location}(P2) + \\ &\quad m(Z) + \text{Location}(N2) \\ &= 0.420,\end{aligned}$$

We have taken the range of temperature set points for distillation of kerosene from the crude oil refinery [Ebbesen

(1992)] to find a precise temperature set point. This data is analysed with rules of fuzzy control theory. The fuzzy rules expressed in terms of degree of membership grade to each temperature set point. Finally the ultimate membership grade was obtained using Center Max-Min rule for the distillation of kerosene.

It has been observed that the highest membership grade using Center Max-Min rule was given to the temperature set point $223^0$C.

### 3.4 Determination of Temperature Set Point of Naphtha Resulting in Better Distillation using Fuzzy Control Theory

The given possible ranges of temperature set points are {$155^0$C, $156^0$C, $157^0$C, $159^0$C, $160^0$C, $161^0$C, $162^0$C, $163^0$C, $164^0$C, $165^0$C} and possible distillation are {93%, 94%, 95%, 96%, 97%} in the case of naphtha as observed by Ebbesen(1992). Using this data as inputs of fuzzy control theory, we identify the precise temperature set points from possible range of temperature set points. Now membership grade is assigned to the input of each temperature set point and each percentage of distillation. The following graph is represents the membership grades of temperature set point.

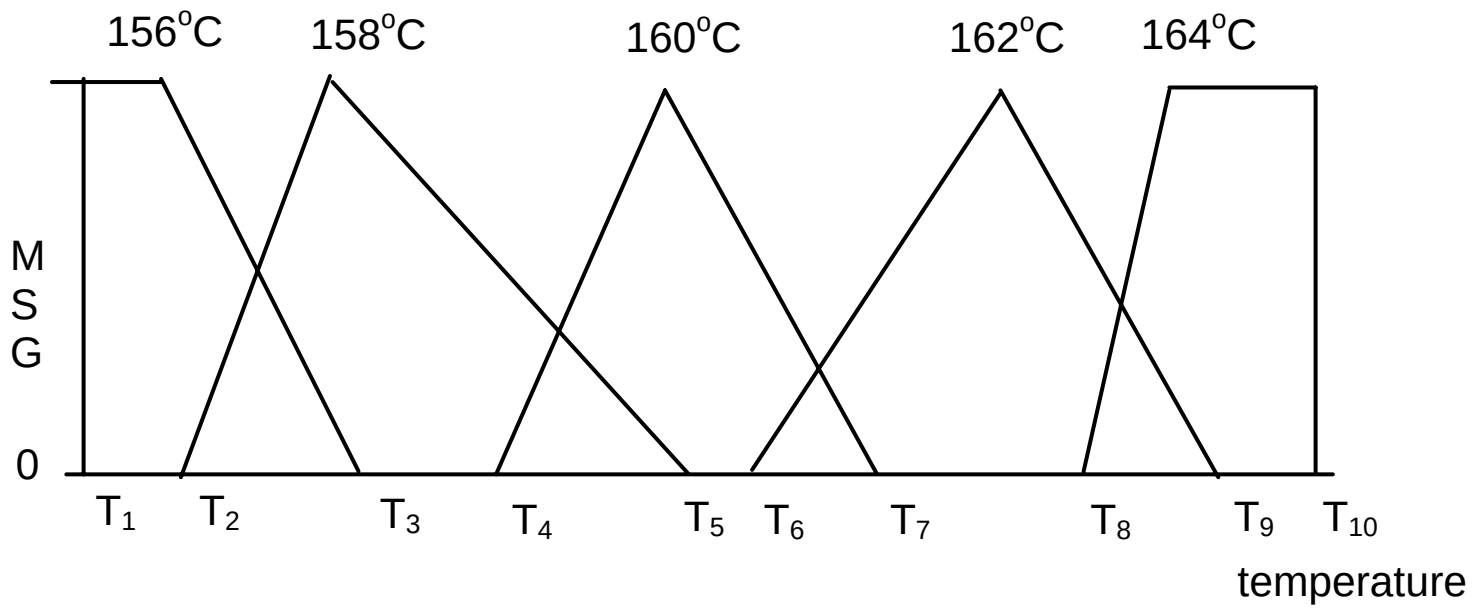

FIGURE 3.4.1: Membership grade of temperature set points

After assigning membership grades in the interval [0,1] to each input of the percentage of distillation the following graph is obtained representing the membership grades of percentage of distillation.

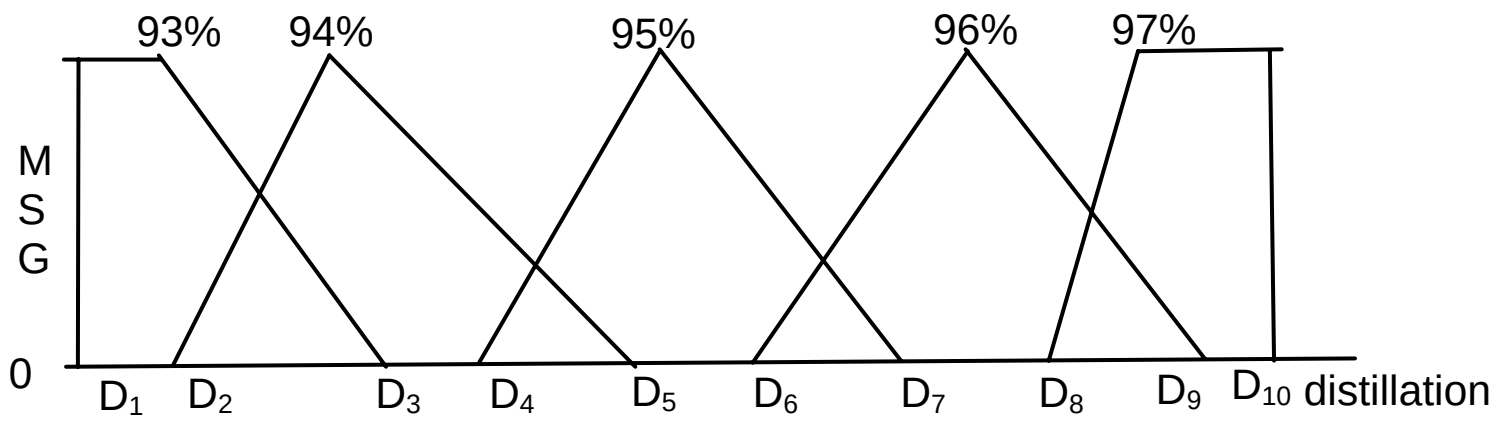

Legend
MSG: Membership Grade
T: Temperature set point
D: Distillation

FIGURE 3.4.2: Membership function of percentage of distillation

The membership grade varies from 0 to 1.

For getting precise temperature set point for naphtha the throttle variables(the grade of membership) quantified into five subsets. Here fuzzy rules are used to find the possible percentage of distillation for each temperature set point and the Center Max-Min rule is used to find a throttle membership grade for the existing fuzzy rules. To get the grade of membership to each existing fuzzy rule, throttle variables are quantified into five subsets as follows:

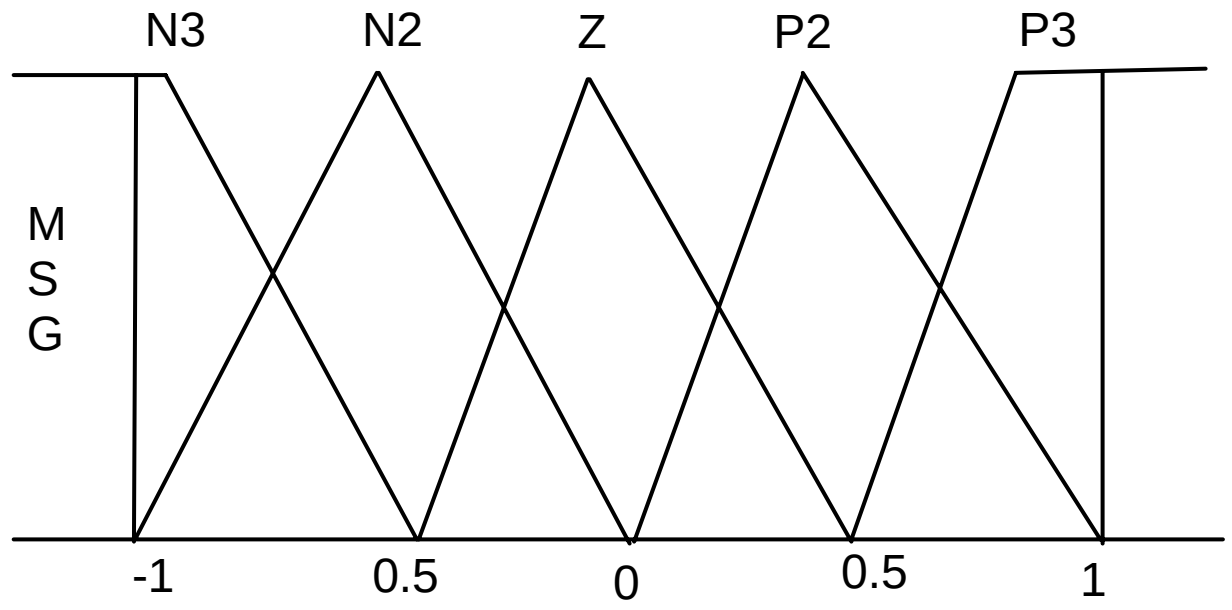

N3: Very Big Negative
N2: Big Negative
Z : Nominal
P2: Big Positive
P3: very Big Positive

FIGURE 3.4.3: Throttle values

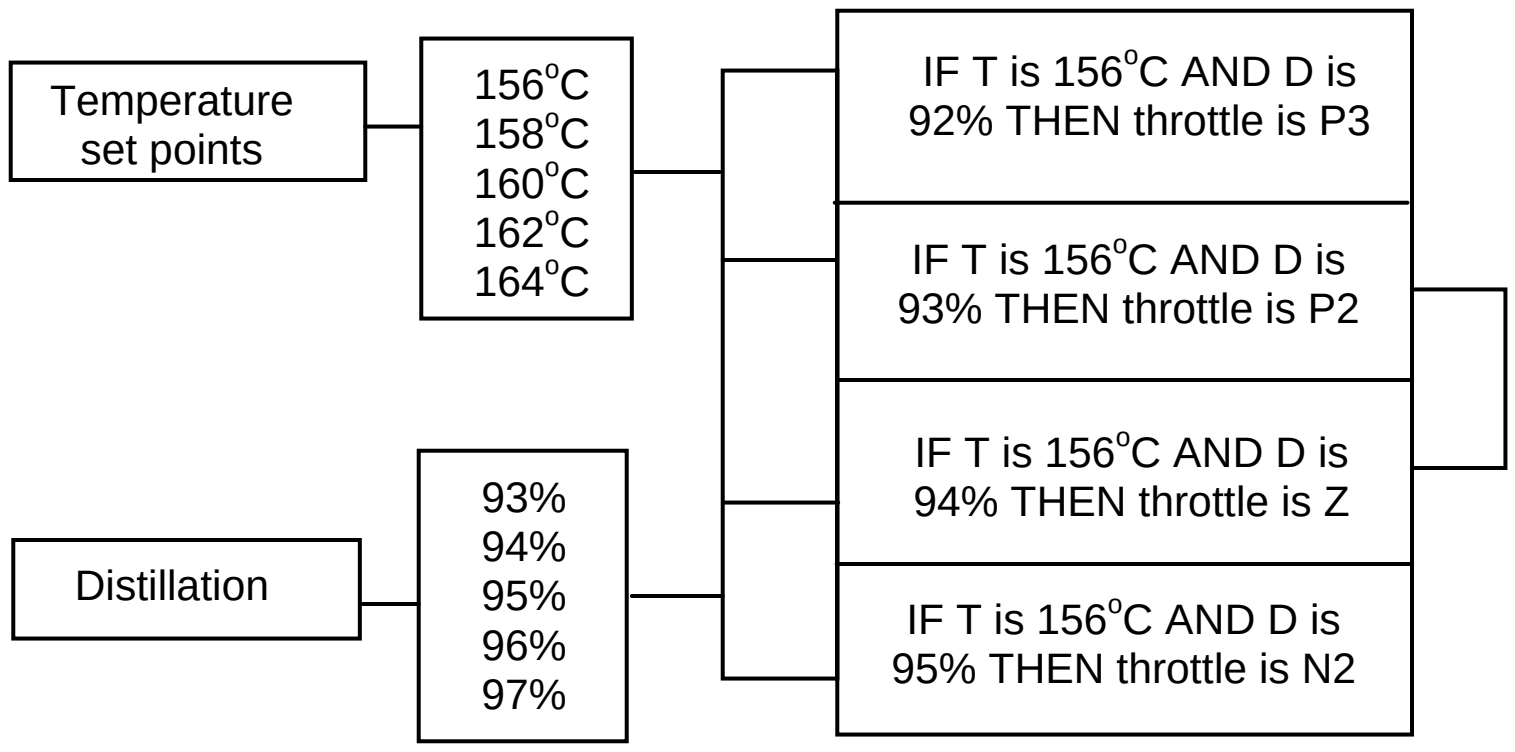

FIGURE 3.4.4: Fuzzy rules for the temperature set point $156^o$C

**The fuzzy rule for distillation of naphtha :**

Rule -1 : If T is $156^0$C temperature set point AND D is 92% THEN throttle is P3.
Rule -2 : If T is $156^0$C temperature set point AND D is 93% THEN throttle is P2.
Rule - 3 : If T is $156^0$C temperature set point AND D is 94% THEN throttle is Z.
Rule - 4 : If T is $156^0$C temperature set point AND D is 95% THEN throttle is N2.

We calculate the throttle value to the temperature set point for $156^0$C. In the above stated rules, only rule -2, and rule-3 are applicable that is only we get 93% and 94% of distillation respectively.

**Rule-2**
The throttle value to the temperature set point $156^0$C for 93% distillation is calculated as follows.

Throttle = (0.34 + 0.57)/2 = 0.455,

The graphical representation of the membership grade of the temperature set point $156^0$C for 93 percentage of distillation is as follows.

Graph 3.4.1: Graphical representation of Rule 2

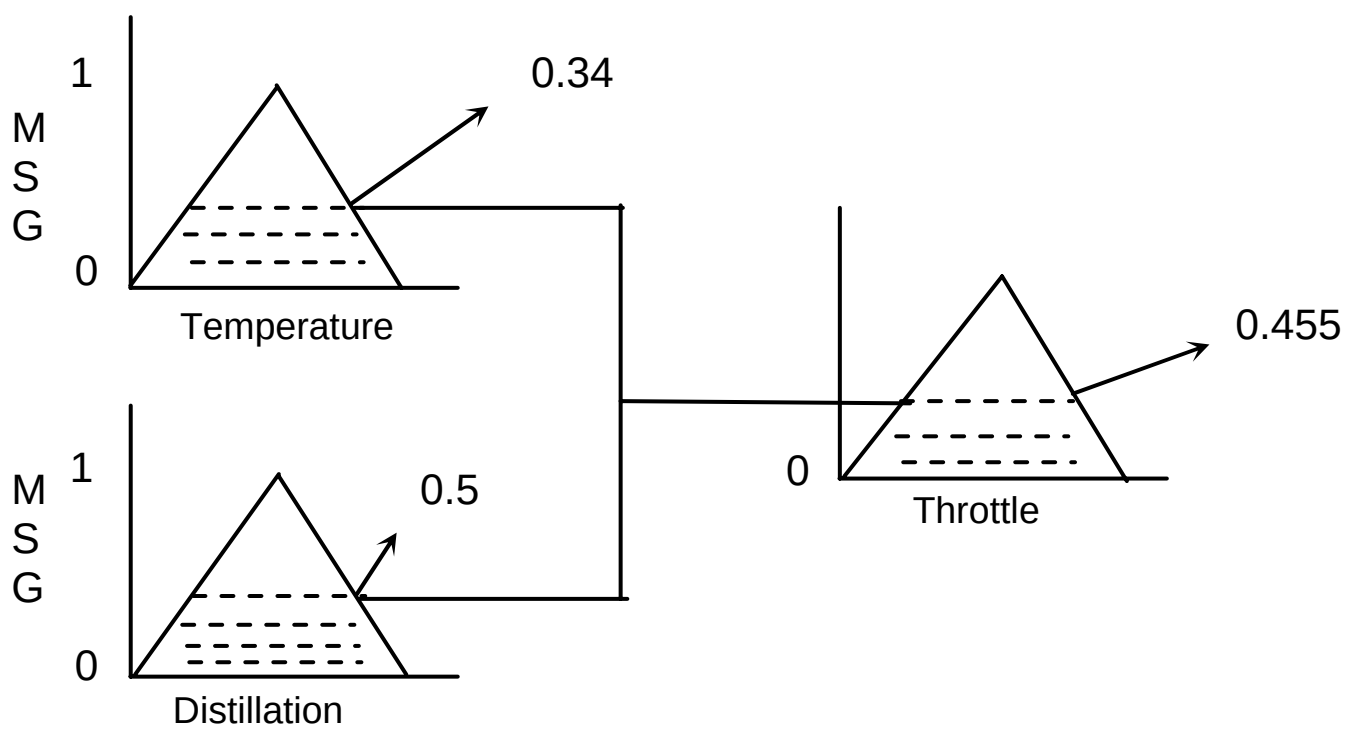

Graph 3.4.2: The two outputs are then defuzzified
by center max-min rule

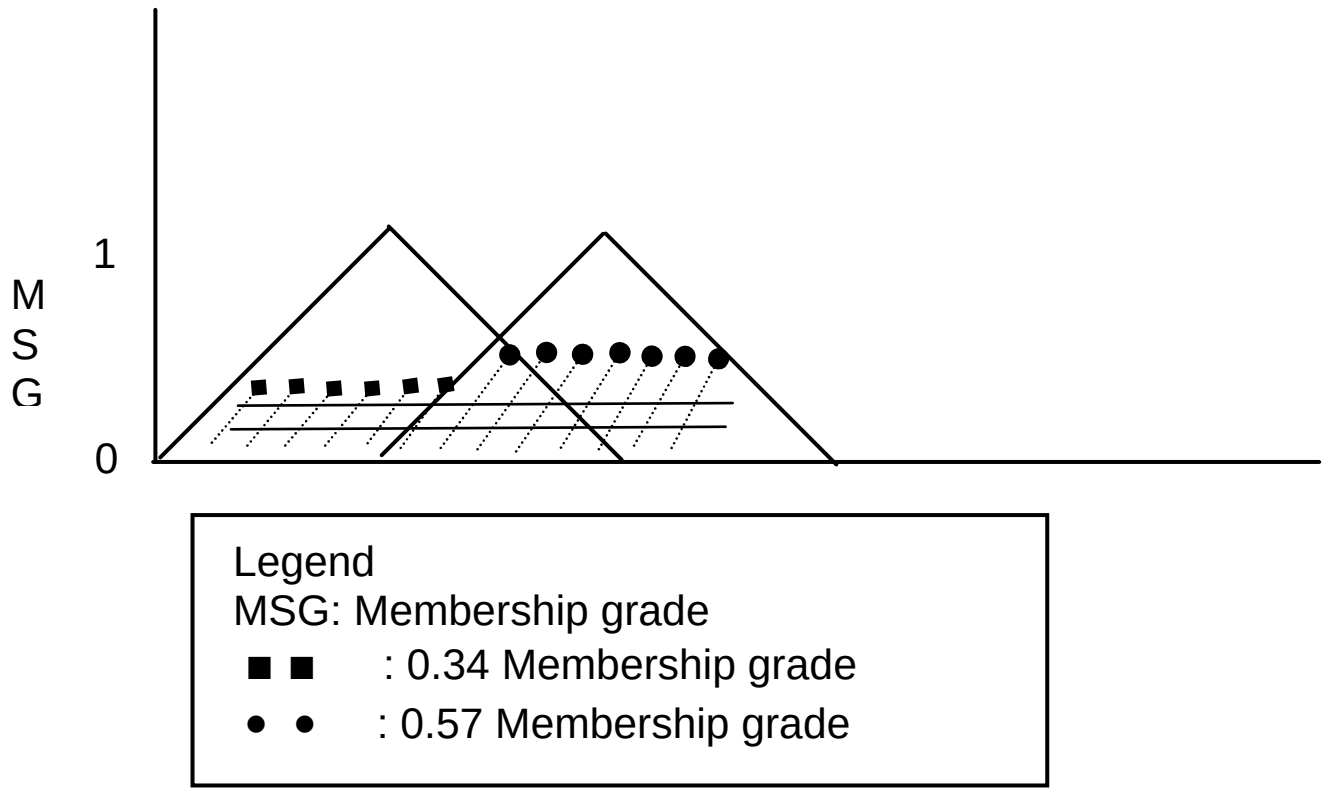

## Rule-3

The throttle value to the temperature set point 156$^0$C for 94% distillation is calculated as follows.

Throttle = (0.59+0.42)/2 = 0.505,

The graphical representation of the membership grade of the temperature set point $156^0C$ for 94 percentage of distillation is as follows.

Graph 3.4.3: Graphical representation of Rule 3

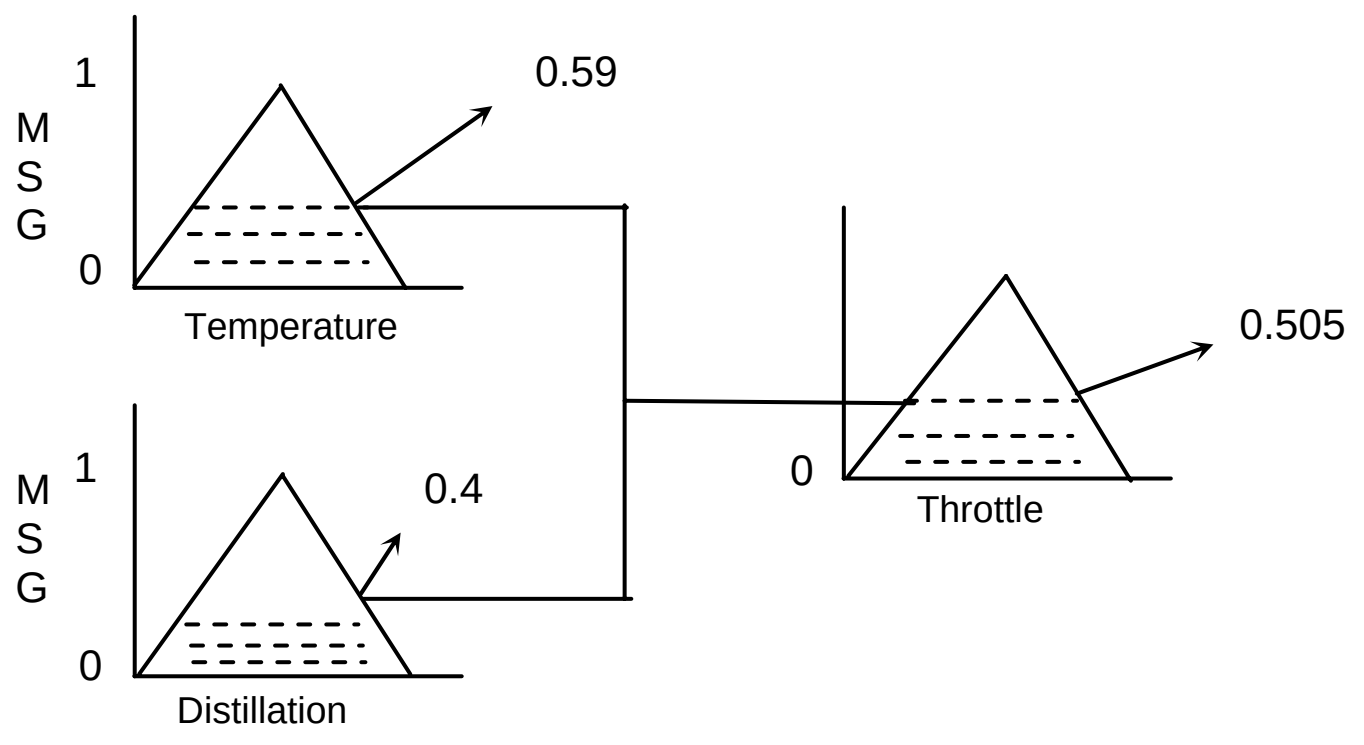

Graph 3.4.4: The two outputs are then defuzzified by center max-min rule

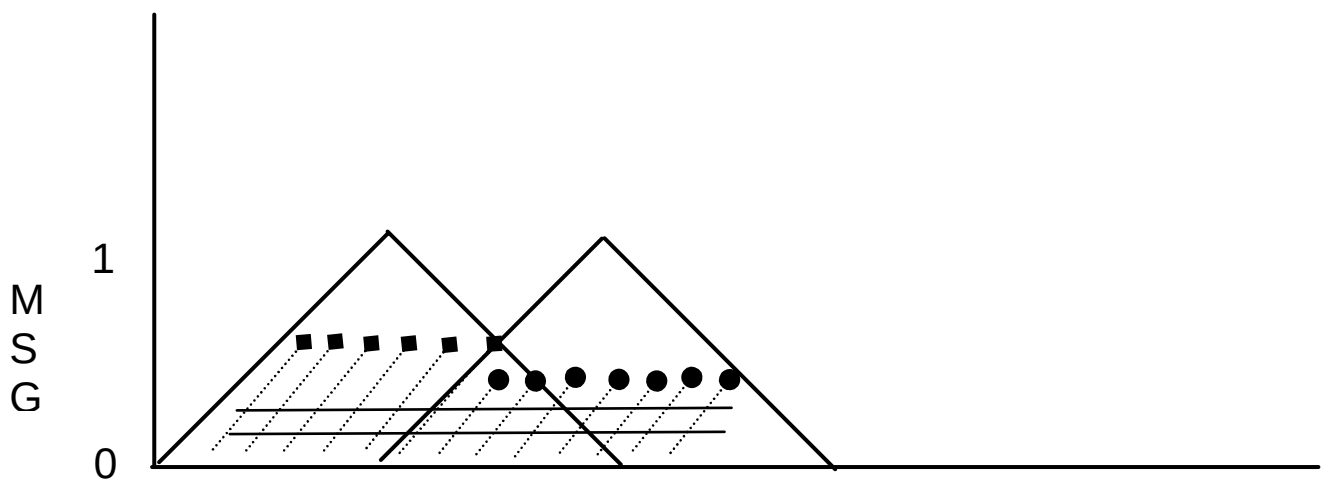

Legend
MSG: Membership grade
■ ■ : 0.59 Membership grade
• • : 0.42 Membership grade

Here, the Center Max-Min rule is used to find a precise temperature set point by grade f membership(throttle)value.

Here, the Center Max-Min rule is used to find a precise temperature set point.

**Using Center Max-Min rule to find precise temperature set point for naphtha**

Throttle(grade of membership) = m(P3) $\times$ Location(P2) + m(Z) + Location(N2)
= 0.730.

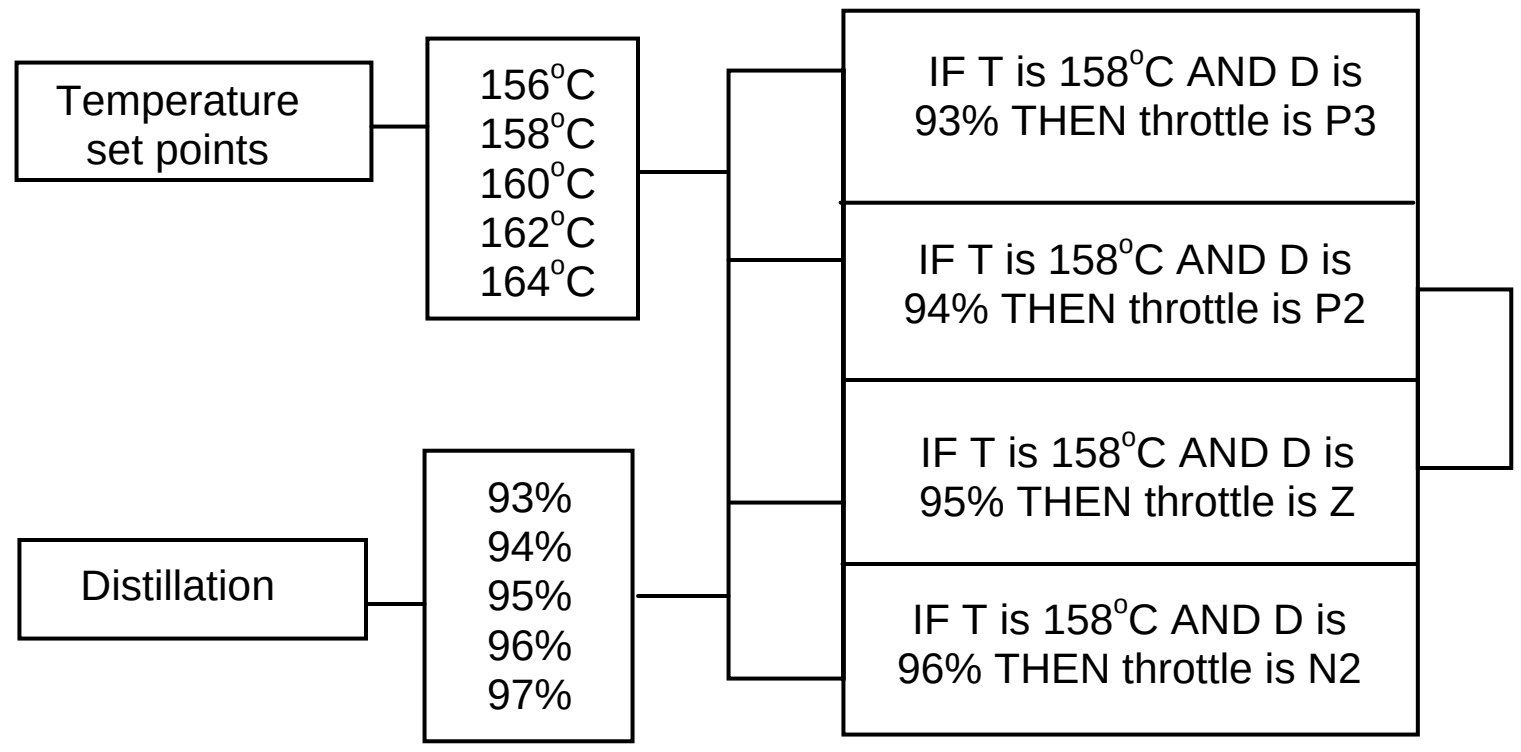

FIGURE 3.4.5: Fuzzy rules for the temperature set point $158^{o}$C

**The fuzzy rule for distillation of naphtha :**

Rule -1 : If T is $158^0$C temperature set point AND D is 93% THEN throttle is P3.
Rule -2 : If T is $158^0$C temperature set point AND D is 94% THEN throttle is P2.
Rule - 3 : If T is $158^0$C temperature set point AND D is 95% THEN throttle is Z.
Rule - 4 : If T is $158^0$C temperature set point AND D is 96% THEN throttle is N2.

We calculate the throttle value to the temperature set point for $158^0$C. In the above stated rules, only rule -2, and rule-3 are applicable that is distillation is 94% and 95% respectively.

**Rule-2**
The throttle value to the temperature set point $158^0$C for 94% is calculated as follows:

Throttle = (0.59+0.42)/2 = 0.505,

**Rule-3**

The throttle value to the temperature set point $158^0$C for 95% distillation is calculated as follows.

Throttle = (0.5+0.35)/2 = 0.425,

**Using Center Max-Min rule to find precise temperature set point for Naphtha**

Throttle(grade of membership) = m(P3) $\times$ Location(P2) + m(Z) + Location(N2)
= 0.717,

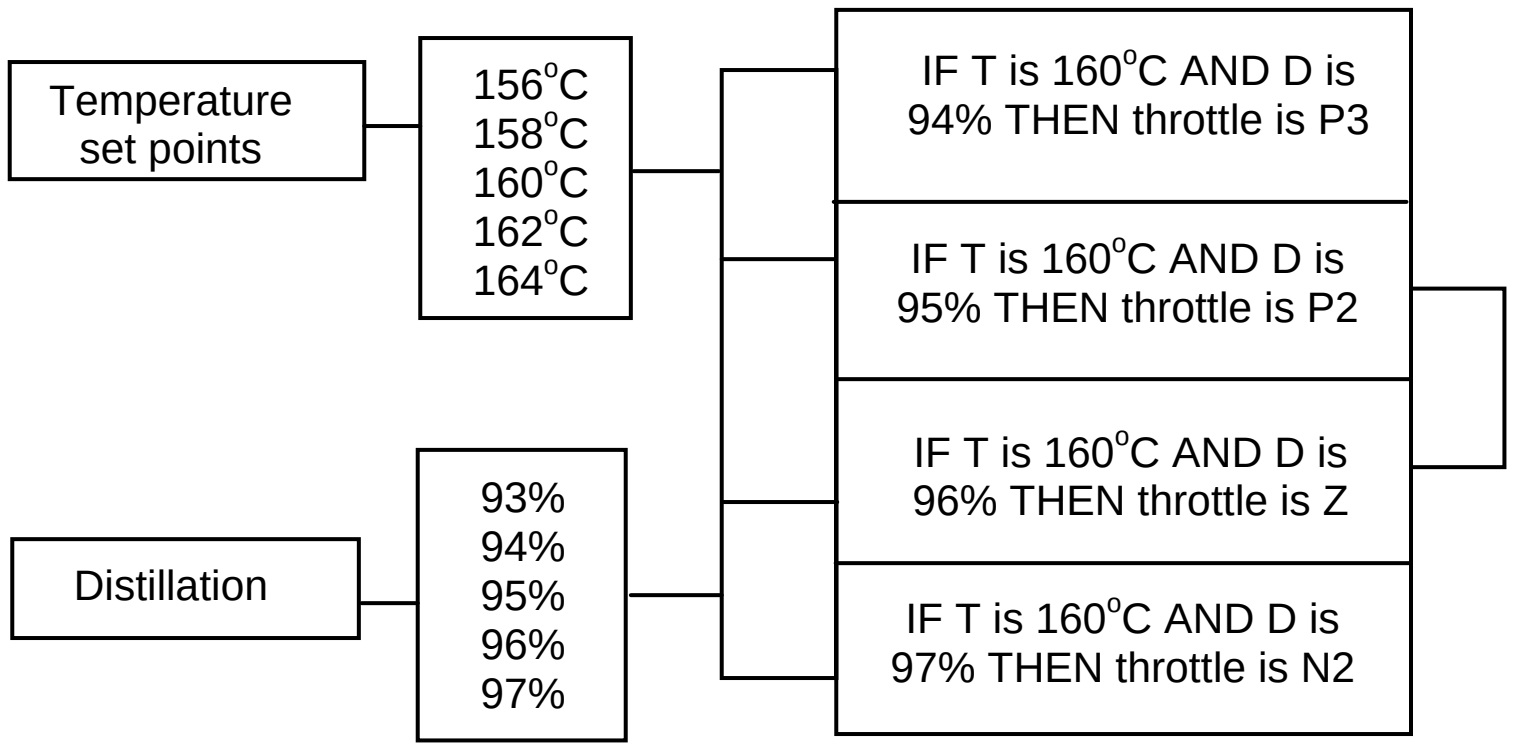

FIGURE 3.4.6: Fuzzy rules for the temperature set point $160^o$C

**The fuzzy rule for distillation of naphtha :**
Rule -1 : If T is $160^0$C temperature set point AND D is 94% THEN throttle is P3.
Rule -2 : If T is $160^0$C temperature set point AND D is 95% THEN throttle is P2.
Rule - 3 : If T is $160^0$C temperature set point AND D is 96% THEN throttle is Z.
Rule - 4 : If T is $160^0$C temperature set point AND D is 97% THEN throttle is N2.

We obtain the throttle value to the temperature set point for $160^0$C. In the above stated rules, only rule -2, and rule-3 are applicable that is distillation is 95% and 96% respectively.

**Rule-2**

The throttle value to the temperature set point $160^0$C for 95% distillation is calculated as follows.

$$\text{Throttle} = (0.5+0.35)/2 = 0.425,$$

**Rule-3**

The throttle value to the temperature set point $160^0$C for 96% distillation is as follows.

$$\text{Throttle} = (0.59+0.5)/2 = 0.55,$$

**Using Center Max-Min rule to find precise temperature set point for naphtha**

$$\begin{aligned} \text{Throttle(grade of membership)} &= m(P3) \times \text{Location}(P2) + m(Z) + \text{Location}(N2) \\ &= 0.7625, \end{aligned}$$

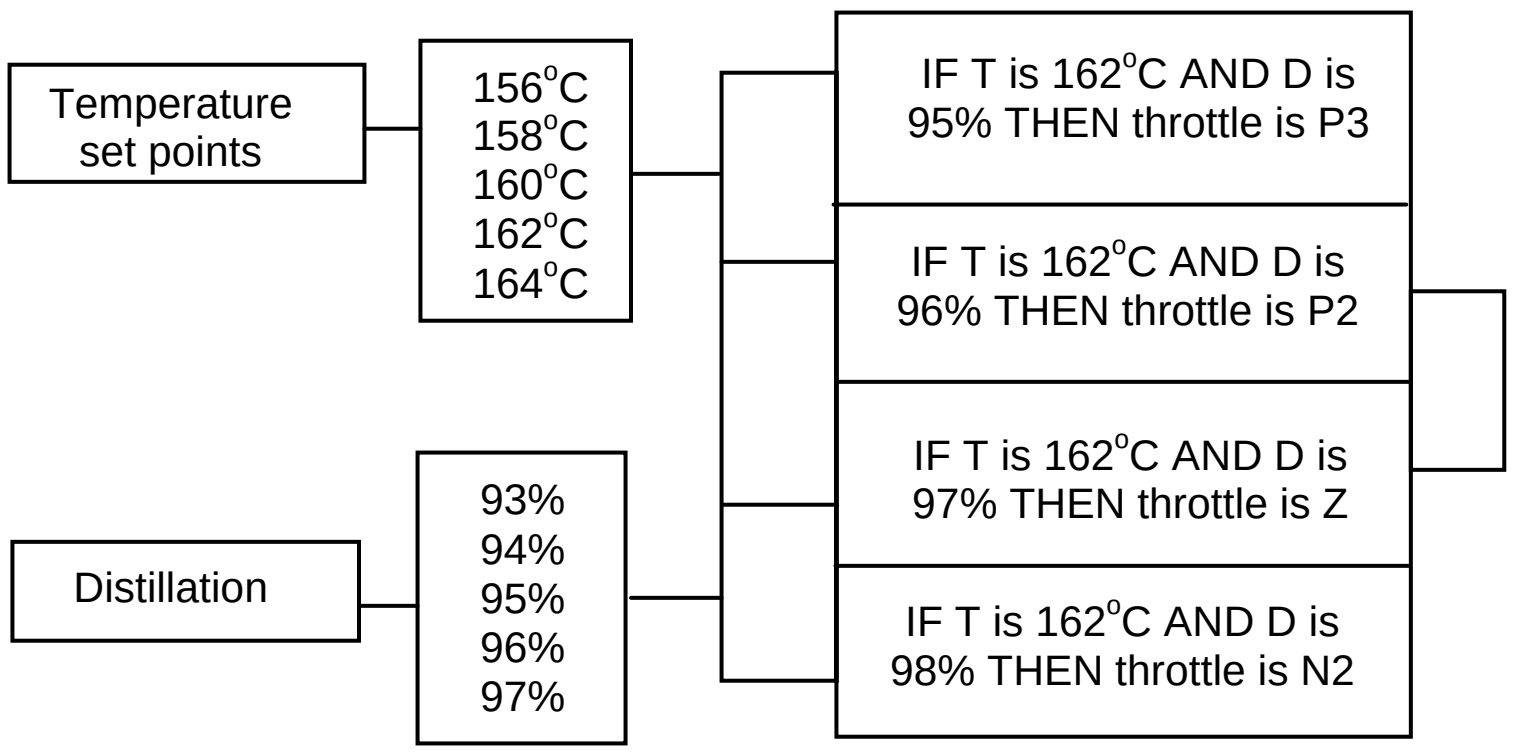

FIGURE 3.4.9: Fuzzy rules for the temperature set point $162^o$C

We calculate the throttle value to the temperature set point for $162^0$ C. In the above stated rule only rule -2, and rule-3 are applicable that is distillation is 96% and 97% respectively.

**Rule-2**

The throttle value to the temperature set point $162^0$ C for 96% distillation is calculated as follows.

Throttle = (0.5+0.6)/2 = 0.55

**Rule -3**

The throttle value to the temperature set point $162^0$ C for 97% distillation is calculated as follows

Throttle = (0.45+.55)/2 = 0.48.

**Using Center Max-Min rule to find precise temperature set point for kerosene**

Throttle(grade of membership) = m(P3) $\times$ Location(P2) + m(Z)+Location(N2)

Throttle = 0.79,

We have taken a range of temperature set points for the distillation of naphtha from the crude oil refinery to find a precise temperature set point. This data analyzed with rules of fuzzy control theory. The fuzzy rules are expressed in terms of the degree of membership grade to each temperature set point. Finally the ultimate membership grade was obtained using the centre max-min rule for the distillation of naphtha.

We observe that the highest membership grade using Center Max-Min rule result in maximum distillation of naphtha for the temperature set point $162^0$C.

### 3.5 Determination of Temperature Set-Point of Gasoil Resulting in Better Distillation Using Fuzzy Control Theory

The given possible ranges of temperature set points are {–5.50$^0$C, –5.00$^0$C, –4.50$^0$ C, –4.00$^0$ C, –3.50$^0$C} and possible percentages of distillation are (93%, 94%, 95%, 96%, 97%} in the case of gasoil as observed by Ebbesen(1992). Using this data as input we use fuzzy control theory to find a precise temperature set point for gasoil.

Now membership grade is assigned to the input of each temperature set point and each distillation. After membership grades are assigned to each input of temperature set points, the following graph results representing the membership grades of temperatures set point.

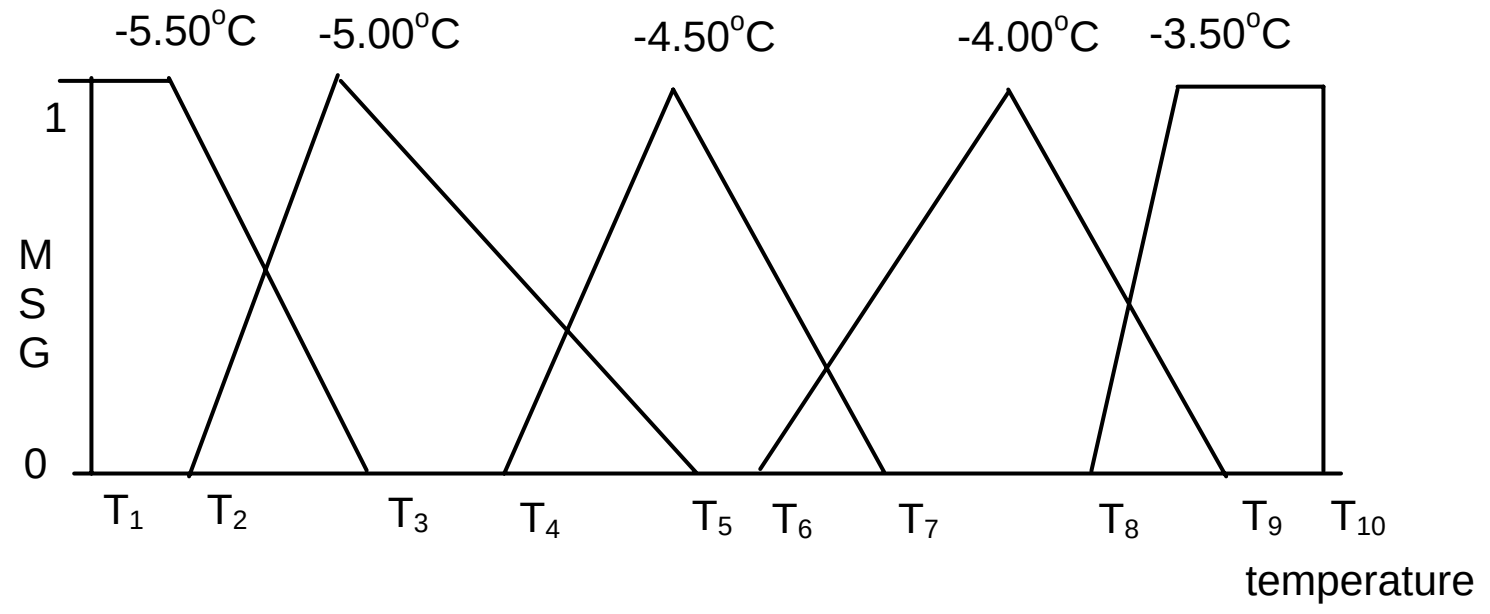

FIGURE 3.5.1: Membership grade of temperature set points

After membership grades are assigned to each input of temperature set points, the following graph results representing the membership grades of temperatures set point.

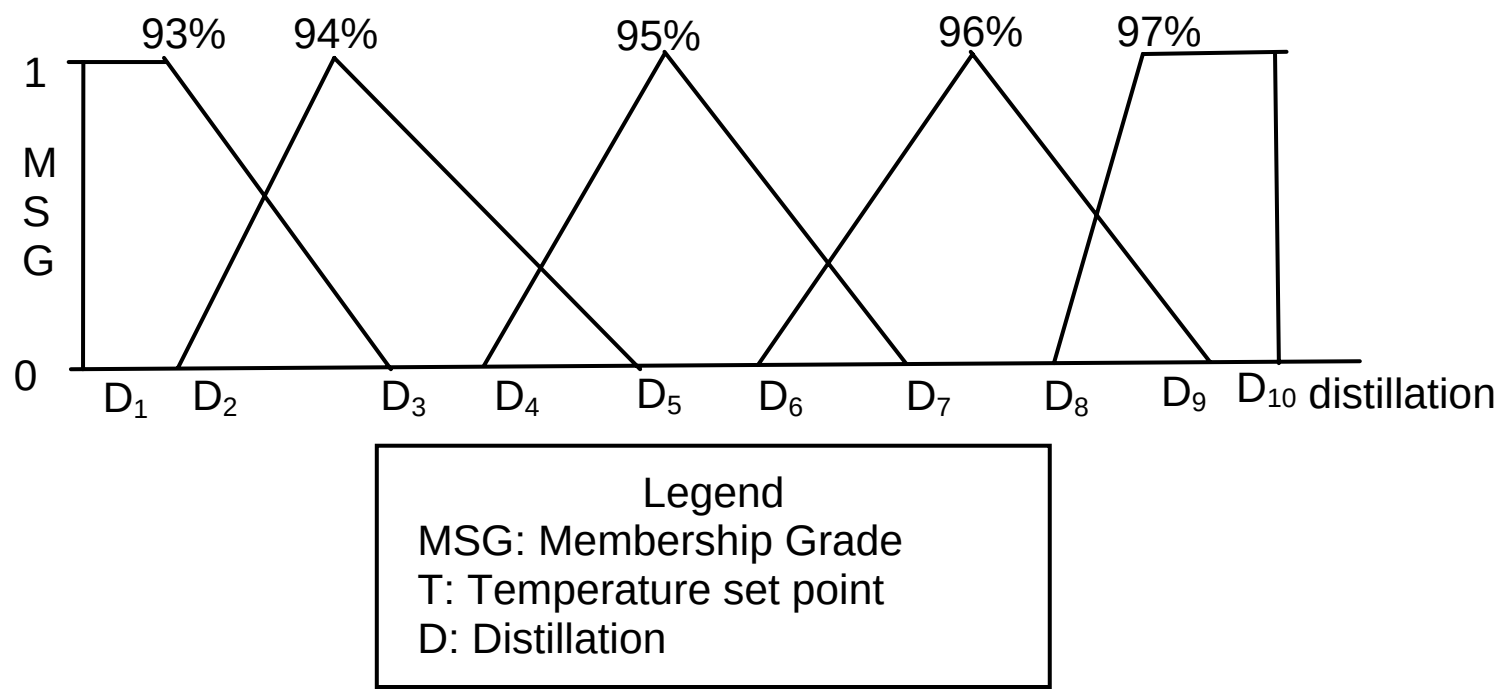

FIGURE 3.5.2: Membership grades of percentage of distillation

The membership grade varies from 0 to 1.
For getting the precise temperature set point for gasoil the throttle variables(the grade of membership) are quantified into five subsets. Here fuzzy rules are used to find the possible percentage of distillation for each temperature set point and the Centre Max-Min rule is used to find a throttle membership grade for existing fuzzy rules. To get the grade of membership to each existing fuzzy rule throttle variables are quantified into five subsets as follows.

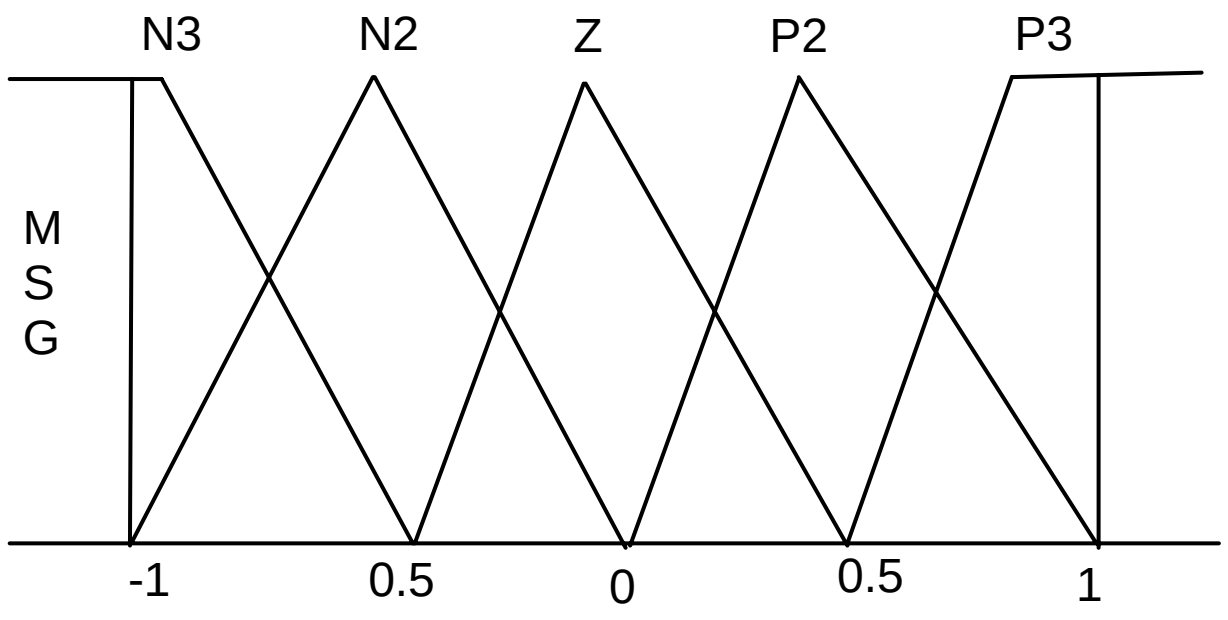

N3: Very Big Negative
N2: Big Negative
Z : Nominal
P2: Big Positive
P3: very Big Positive

FIGURE 3.5.3: Throttle values

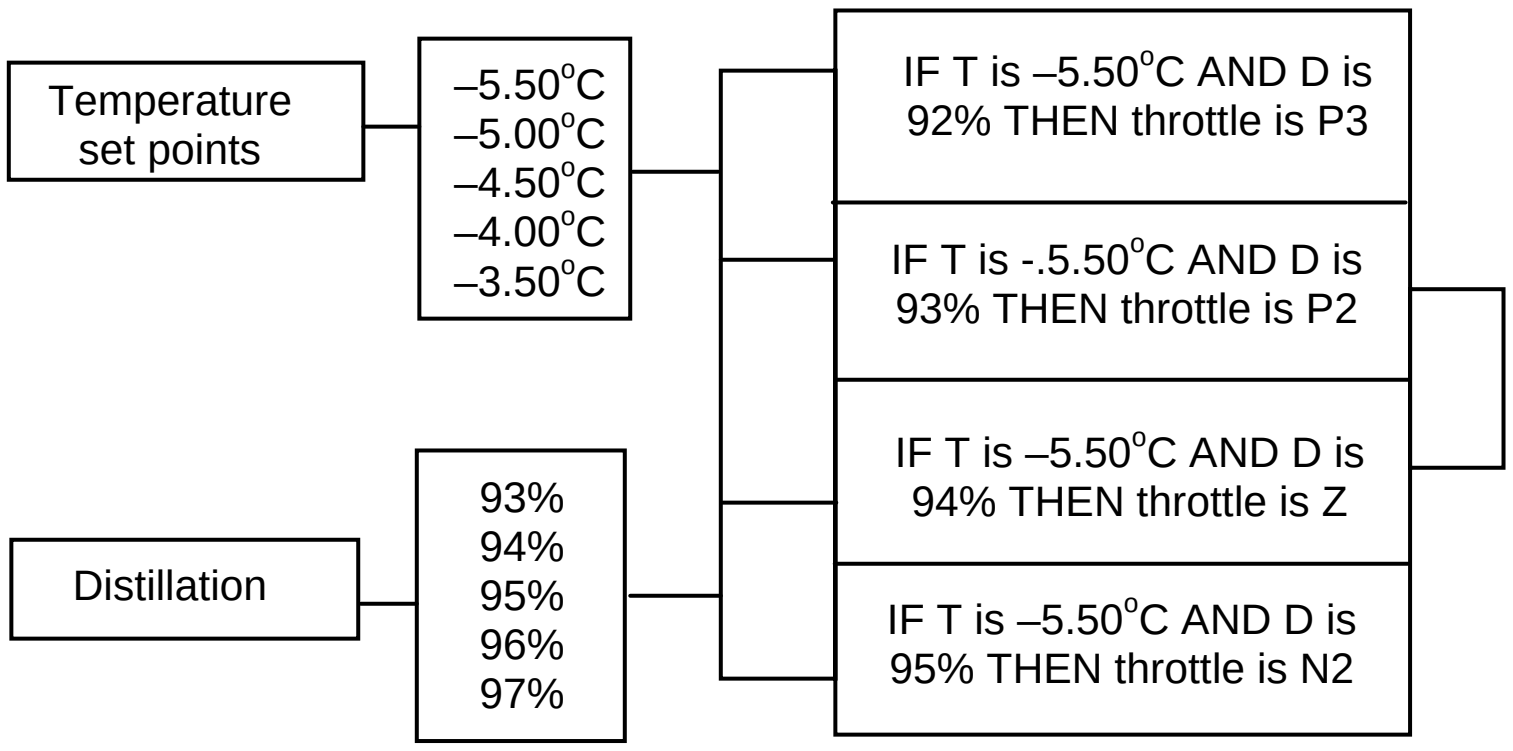

FIGURE 3.5.4: Fuzzy rules for the temperature set point –5.50°C

**The fuzzy rule for distillation of gasoil :**

Rule -1 : If T is -5.50 degree celsius temperature AND D is 92% THEN throttle is P3.

Rule -2 : If T is -5.50 degree celsius temperature AND D is 93% THEN throttle is P2.

Rule - 3 : If T is -5.50 degree celsius temperature AND D is 94% THEN throttle is Z.

Rule - 4 : If T is -5.50 degree celsius temperature AND D is 95% THEN throttle is N2.

We find the throttle value to the temperature set point -5.50$^0$C. In the above stated rules, only rule -2, and rule-3 are applicable that is only we get 93% and 94% of distillation respectively.

**Rule-2**

The throttle value to the temperature set point -5.50$^0$C for 93% distillation is calculated as follows.

Throttle = (0.3+0.4)/2=0.35,

The graphical representation of the membership grade of the temperature set point -5.50$^0$C for 93 percentage of distillation is as follows.

Graph 3.5.1: Graphical representation of Rule 2

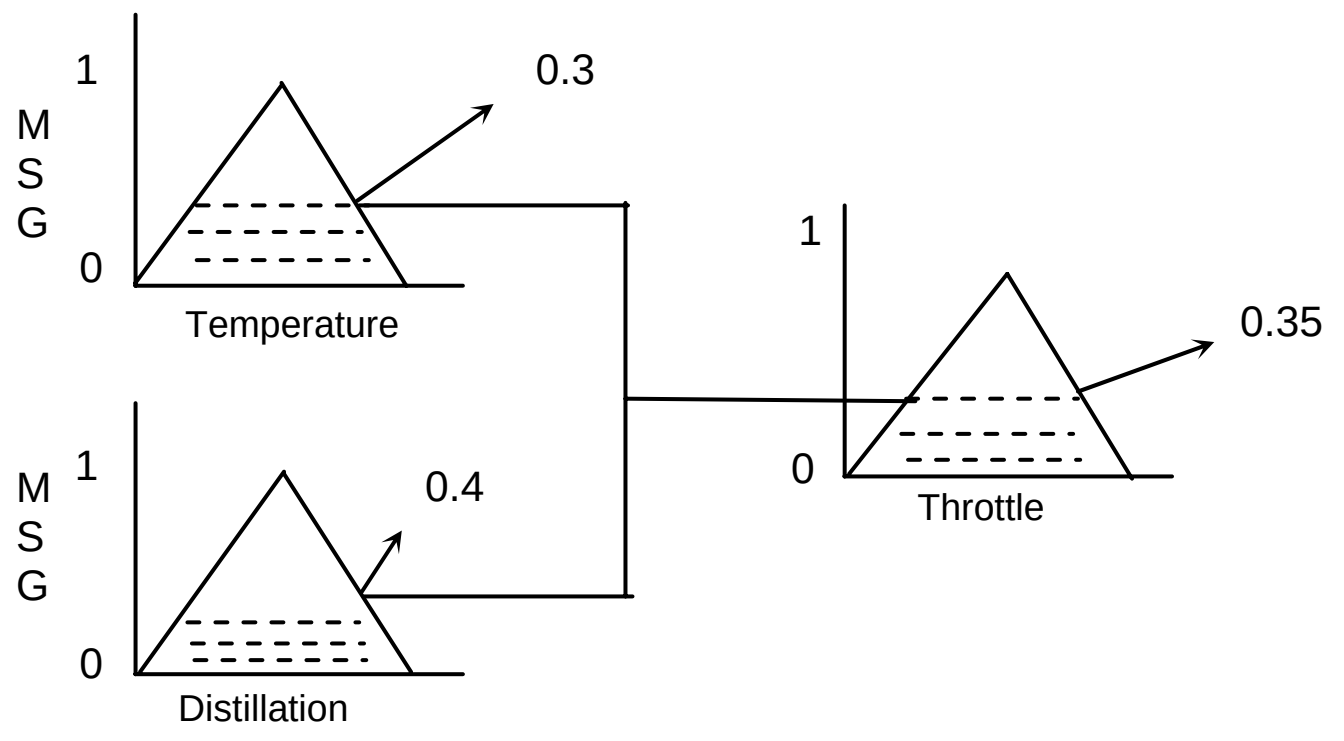

Graph 3.5.2: The two outputs are then defuzzified by center max-min rule

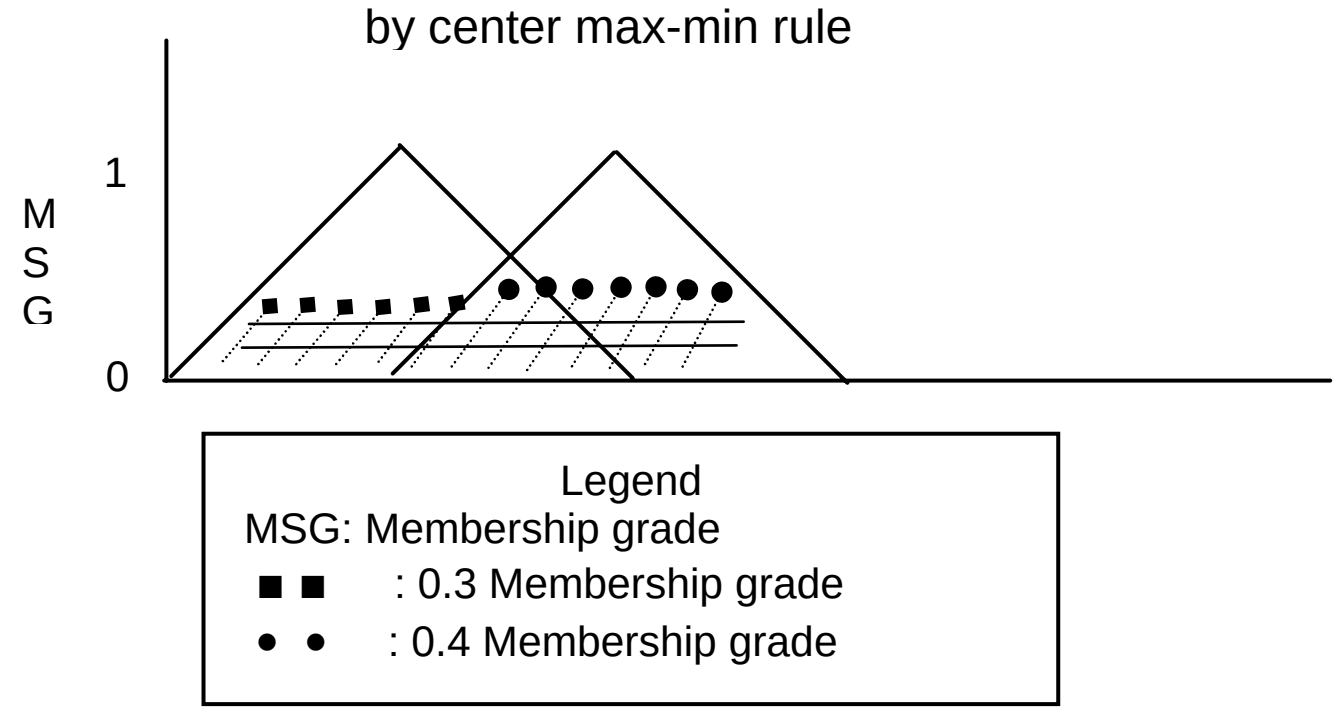

**Rule-3**

The throttle value to the temperature set point -5.50$^0$C for 94% distillation is as follows.

Throttle = (0.2 + 0.6)/2 = 0.40,

The graphical representation of the membership grade of the temperature set point –5.50$^0$C for 94 percentage of distillation is as follows.

Graph 3.5.3: Graphical representation of Rule 3

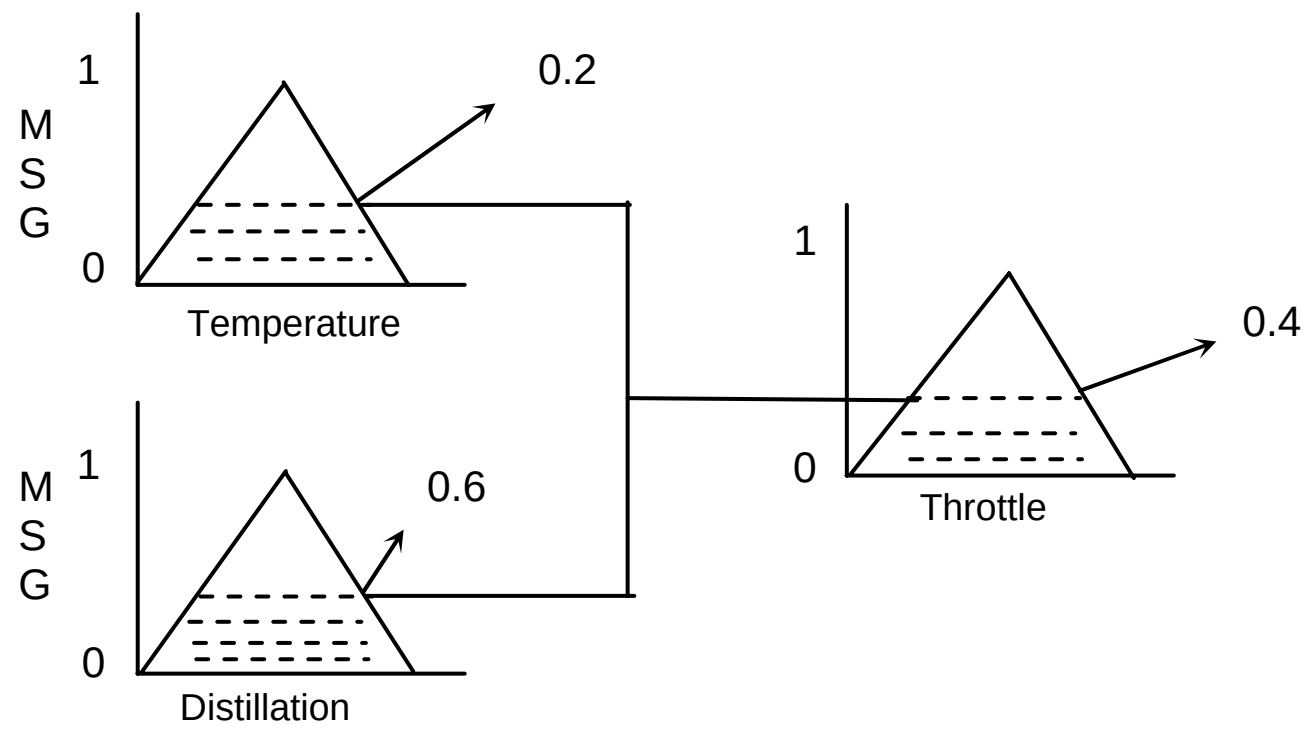

Graph 3.5.4: The two outputs are then defuzzified by center max-min rule

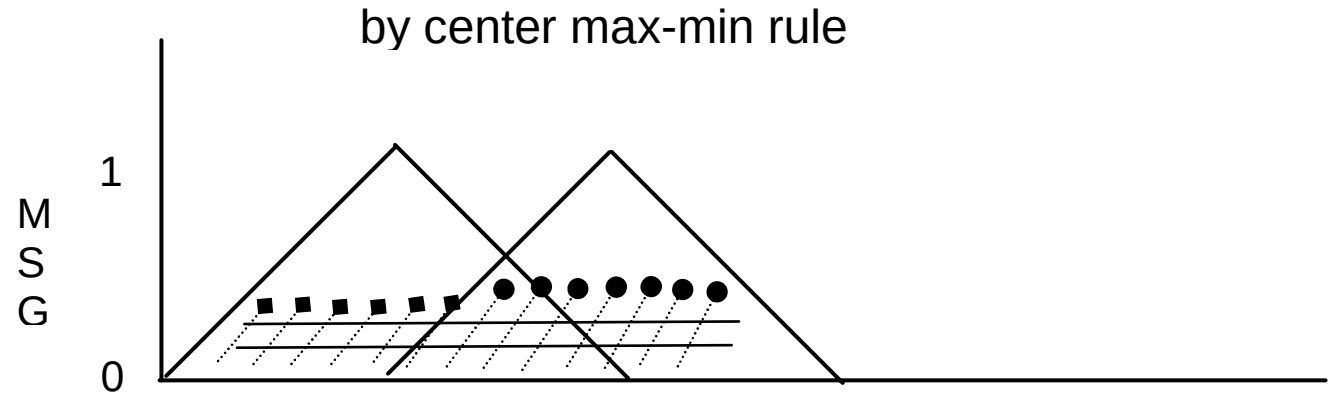

Legend
MSG: Membership grade
■ ■ : 0.34 Membership grade
• • : 0.57 Membership grade

Here, the Centre Max-Min rule is used to find a precise temperature set point by grade of membership(throttle)value.

**Using Center Max-Min rule to find precise temperature set point for gasoil**

Throttle(grade of membership) = m(P3) × Location(P2) + m(Z) + Location(N2)

= 0.35 × 0.5 + 0.4 × 0.5

= 0.375

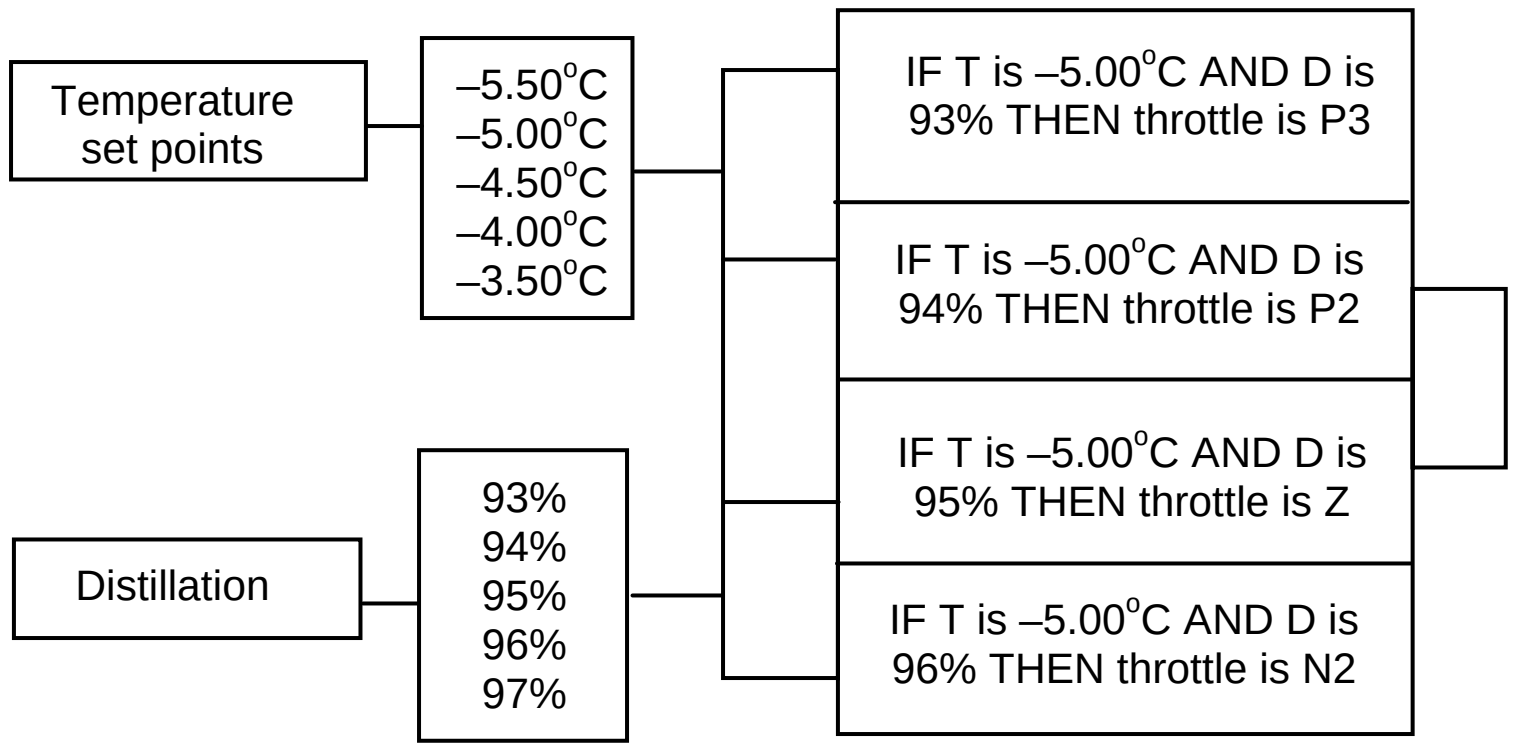

FIGURE 3.5.5: Fuzzy rules for the temperature set point -5.00°C

**The fuzzy rule for the distillation of gasoil :**

Rule -1: If T is -5.00 degree celsius temperature AND D is 93% THEN throttle is P3.

Rule -2: If T is -5.00 degree celsius temperature AND D is 94% THEN throttle is P2.

Rule -3: If T is -5.00 degree celsius temperature AND D is 95% THEN throttle is Z.

Rule -4: If T is -5.00 degree celsius temperature AND D is 96% THEN throttle is N2.

We calculate the throttle value to the temperature set point for $–5.00^0$ C. In the above stated rules, only rule -2, and rule-3 are applicable that is distillation is 94% and 95% of distillation respectively.

**Rule-2**

The throttle value to the temperature set point $-5.00^0$ C for 94% distillation is as follows.

Throttle = (0.6 + 0.5)/2 = 0.55

**Rule-3**

The throttle value to the temperature set point -5.00$^0$C for 95% distillation is as follows.

Throttle = (0.35 + 0.4)/2 = 0.375,

The graphical representation of the membership grade of the temperature set point -5.50$^0$C for 94 percentage of distillation is as follows.

**Using Center Max-Min rule to find precise temperature set point for gasoil**

Throttle(grade of membership) = m(P3) × Location(P2) + m(Z) + Location(N2)

Throttle = 0.55 × 1 + 0.375 × 0.5 / 0.55 + 0.375
= 0.789.

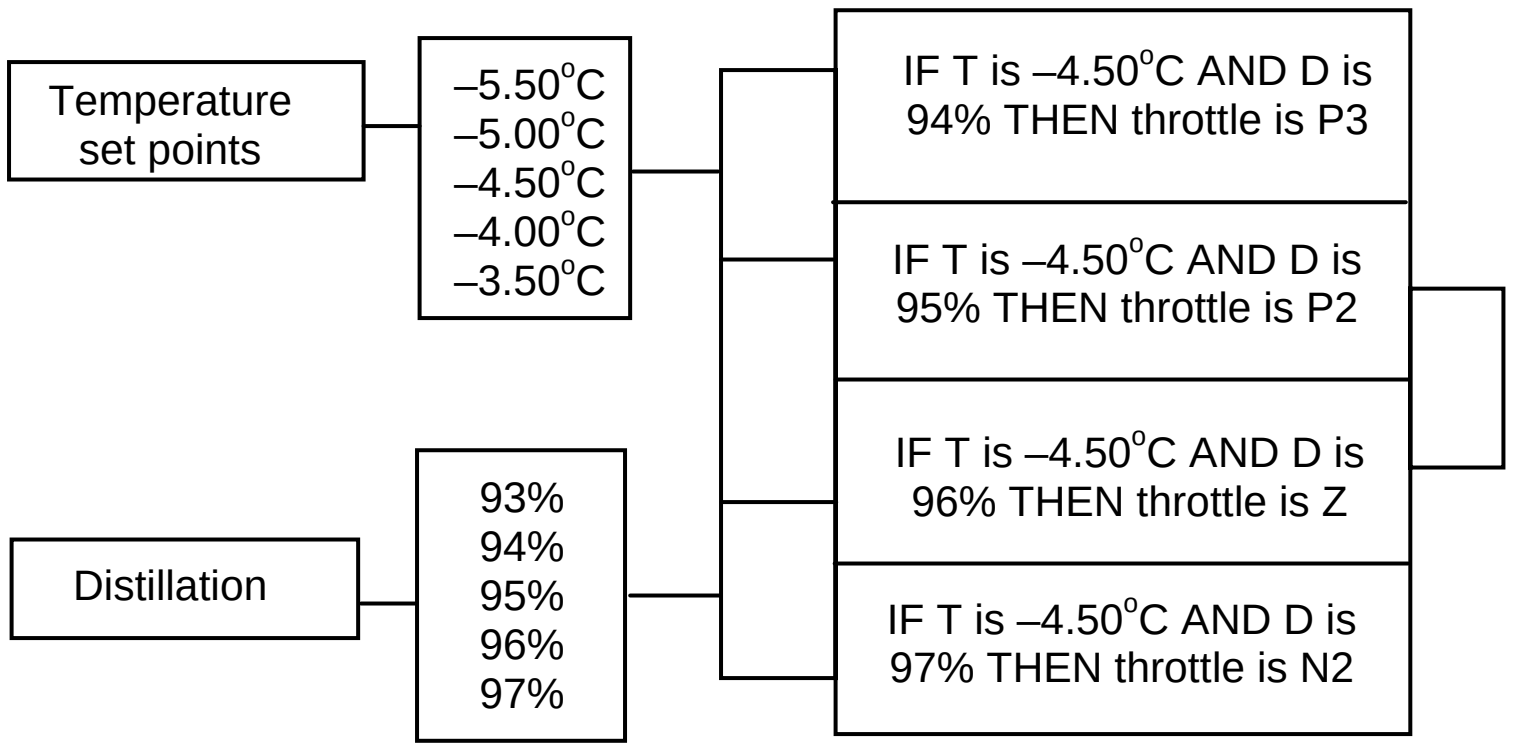

FIGURE 3.5.6: Fuzzy rules for the temperature set point -4.50$^o$C

**The fuzzy rule for the distillation of gasoil :**

Rule -1 : If T is -4.50 degree celsius temperature AND D is 94% THEN throttle is P3.

Rule -2 : If T is -4.50 degree celsius temperature AND D is 95% THEN throttle is P2.

Rule -3 : If T is -4.50 degree celsius temperature AND D is 96% THEN throttle is Z.

Rule -4 : If T is -4.50 degree celsius temperature AND D is 97% THEN throttle is N2.

We calculate the throttle value to the temperature set point for –4.50$^0$ C. In the above stated rules, only rule -2, and rule-3 are applicable that is distillation is 95% and 96% of distillation respectively.

**Rule-2**

The throttle value to the temperature set point -4.50$^0$ C for 95% distillation is calculated as follows.

Throttle = (0.35 + 0.4)/2 = 0.375

**Rule-3**

The throttle value to the temperature set point -4.50$^0$C for 96% distillation is calculated as follows.

Throttle = (0.3 + 0.25)/2 = 0.277,

**Using Center Max-Min rule to find precise temperature set point for gasoil**

Throttle(grade of membership) = m(P3) × Location(P2) + m(Z)+Location(N2)

Throttle = 0.375 × 0.5 + 0.277 × 0.5 / 0.375 + 0.277

= 0.5

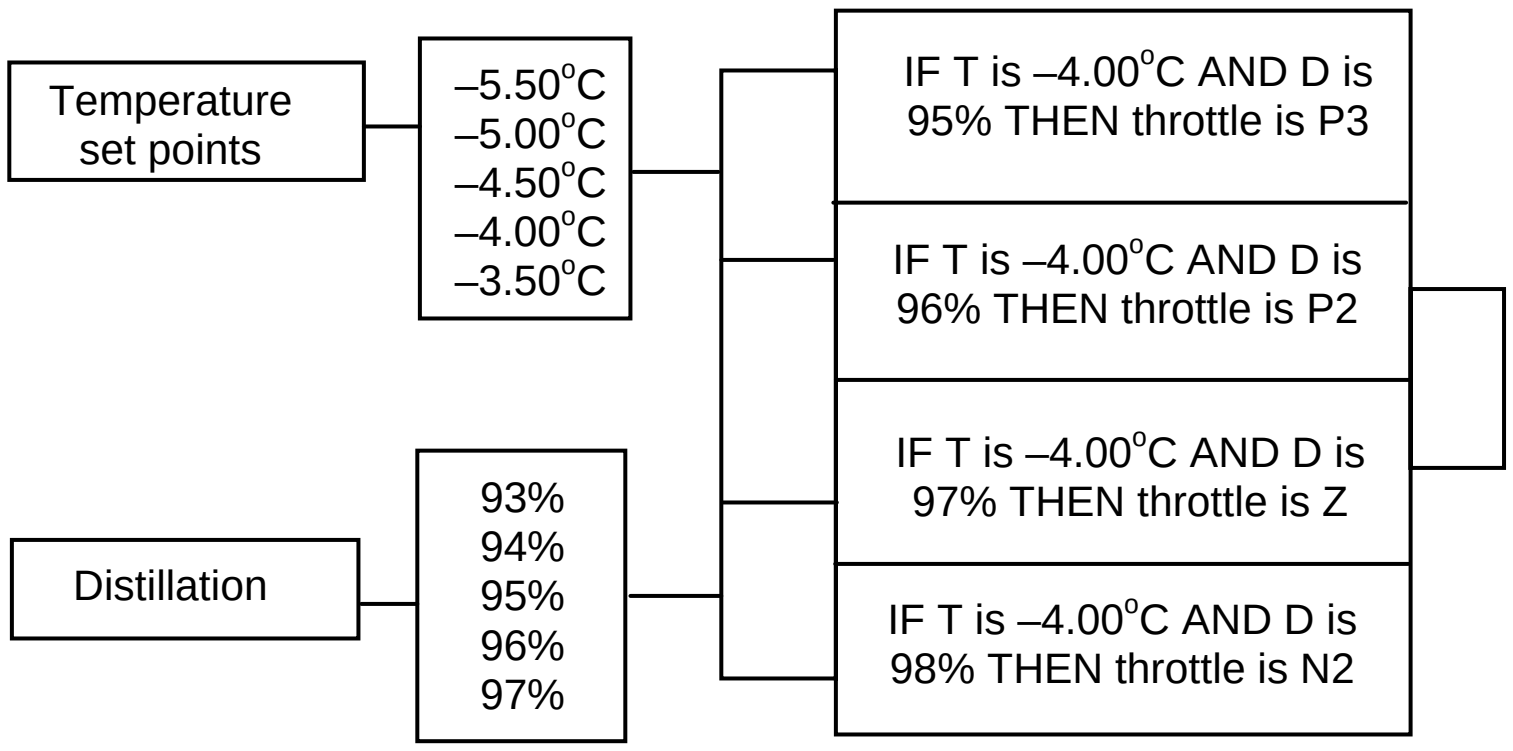

FIGURE 3.5.7: Fuzzy rules for the temperature set point -4.00°C

**The fuzzy rule for the distillation of gasoil :**

Rule -1 : If T is -4.00 degree celsius temperature AND D is 95% THEN throttle is P3.
Rule -2 : If T is -4.00 degree celsius temperature AND D is 96% THEN throttle is P2.
Rule -3 : If T is -4.00 degree celsius temperature AND D is 97% THEN throttle is Z.
Rule -4 : If T is -4.00 degree celsius temperature AND D is 98% THEN throttle is N2.

Consider the temperature is -4.00 degree Celsius and the distillation of gasoil being 96% and 97%. Here rule-1 and rule-4 are not applicable

**Rule-2**

The throttle value to the temperature set point $-4.00^0$ C for 96% distillation is calculated as follows.

$$\text{Throttle} = (0.25 + 0.3)/2 = 0.275$$

**Rule-3**

The throttle value to the temperature set point $-4.00^0$C for 97% distillation is calculated as follows.

$$\text{Throttle} = (0.4+0.4)/2 = 0.4,$$

**Using Center Max-Min rule to find precise temperature set point for gasoil**

Throttle(grade of membership) = m(P3) × Location(P2) + m(Z)+Location(N2)

$$\text{Throttle} = 0.275 \times 0.5 + 0.4 \times 0.5 / 0.275 + 0.4$$
$$= 0.501.$$

We have taken a range of temperature set points for distillation of gasoil from the crude oil refinery [Ebbesen (1992)] to find a precise temperature set point. This data is analysed with rules of fuzzy control theory. The fuzzy rules are expressed in terms of degree of membership grade to each temperature set point. Finally the ultimate membership grade is obtained using Center Max-Min rule for the distillation of gasoil.

The authors have observed that the highest membership grade for -5.00$^0$C using center max-min rule results in maximum distillation of gasoil and gives better quality.

### 3.6 Conclusions

Finding of precise set point temperatures for the distillation of kerosene, naphtha and gasoil have always remained to be uncertain in a crude oil refinery. Fuzzy control theory is able to predict the precise set point of temperature for kerosene, naphtha and gasoil, which guarantees the maximum percentage of distillation, and also the quality for long hours. By this method the random choice of temperature set point from the range of temperature set points, which affects the quality and quantity of crude is completely over come.

## Chapter Four

# STUDY OF FLOW RATES IN CHEMICAL PLANTS

This chapter has 3 sections. Use of FRE to estimate flow rates in chemical plants forms the section one of this chapter. In section two fuzzy neural networks are used to estimate velocity of flow distribution in a pipe network. The final section estimates the three-stage counter current extraction unit again using fuzzy neural networks.

## 4.1 Use of FRE in Chemical Engineering

The use of fuzzy relational equations (FRE) for the first time has been used in the study of flow rates in chemical plants. They have only used the concept of linear algebraic equations to study this problem and have shown that use of linear equations does not always guarantee them with solutions. Thus we are not only justified in using fuzzy relational equation but we are happy to state by adaptation of FRE we are guaranteed of solutions to the problem. We have adapted the fuzzy relational equations to the problem of estimation of flow rates in a chemical plant, flow rates in a pipe network and use of FRE in a 3 stage counter current exaction unit [44].

Experimental study of chemical plants is time consuming expensive and need intensive labor, researchers and engineers prefer only theoretical approach, which is inexpensive and effective. Only linear equations have been used to study: (1). A typical chemical plant having several inter-linked units (2). Flow distribution in a pipe network and (3). A three stage counter current extraction unit. Here, we tackle these problems in 2 stages. At the first stage we use FRE to obtain a solution. This is done by the method of partitioning the matrix as rows. If no solution exists by this method we as the second stage adopt Fuzzy Neural Networks by giving weightages. We by varying the weights arrive at a solution which is very close to the predicted value or the difference between the estimated value and the predicted value is zero. Thus by using fuzzy approach we see that we are guaranteed of a solution which is close to the predicted value, unlike the linear algebraic equation in which we may get a solution and even granted we get a solution it may or may not concur with the predicted value.

To attain both solution and accuracy we tackle the problems using Fuzzy relational equations at the first stage and if no solution is possible by this method we adopt neural networks at the second stage and arrive at a solution.

Consider the binary relation P(X, Y), Q(Y, Z) and R(X, Z) which are defined on the sets $X = \{x_i / i \in I\}$ $Y = \{y_i / j \in J\}$ and $Z\{z_k / k \in K\}$ where we assume that $I = N_n$, $J = N_r$ and $K = N_s$. Let the membership matrices of P, Q and R be denoted by $P = [p_{ij}]$, $Q = [q_{ik}]$ and $R = [r_{ik}]$ respectively, where $p_{ij} = P(x_i, y_j)$, $q_{ik} = Q(y_j, z_k)$ and $r_{ik} = R(x_i, z_k)$ for $i \in I$ (= $N_n$), $j \in J$ (= $N_m$) and $k \in K$ (= $N_s$). Entries in P, Q and R are taken from the interval [0, 1]. The three matrices constrain each other by the equation

$$P \circ Q = R \qquad (1)$$

(where $\circ$ denotes the max-min composition) known as the fuzzy relation equation (FRE) which represents the set of equation

$$\text{Max } p_{ij}q_{jk} = r_{ik} \qquad (2)$$

for all $i \in N_n$, $k \in N_s$. If after partitioning the matrix and solving the equation (1) yields maximum of $q_{jk} < r_{ik}$ for some $q_{jk}$, then this set of equation has no solution. So at this stage to solve the equation 2, we use feed-forward neural networks of one layer with n-neurons with m inputs shown in Figure 4.1.1.

Inputs of the neuron are associated with real numbers $W_{ij}$ referred as weights. The linear activation function f is defined by

$$f(a) = \begin{cases} 0 & \text{if } a < 0 \\ a & \text{if } a \in [0, 1] \\ 1 & \text{if } a > 1 \end{cases}$$

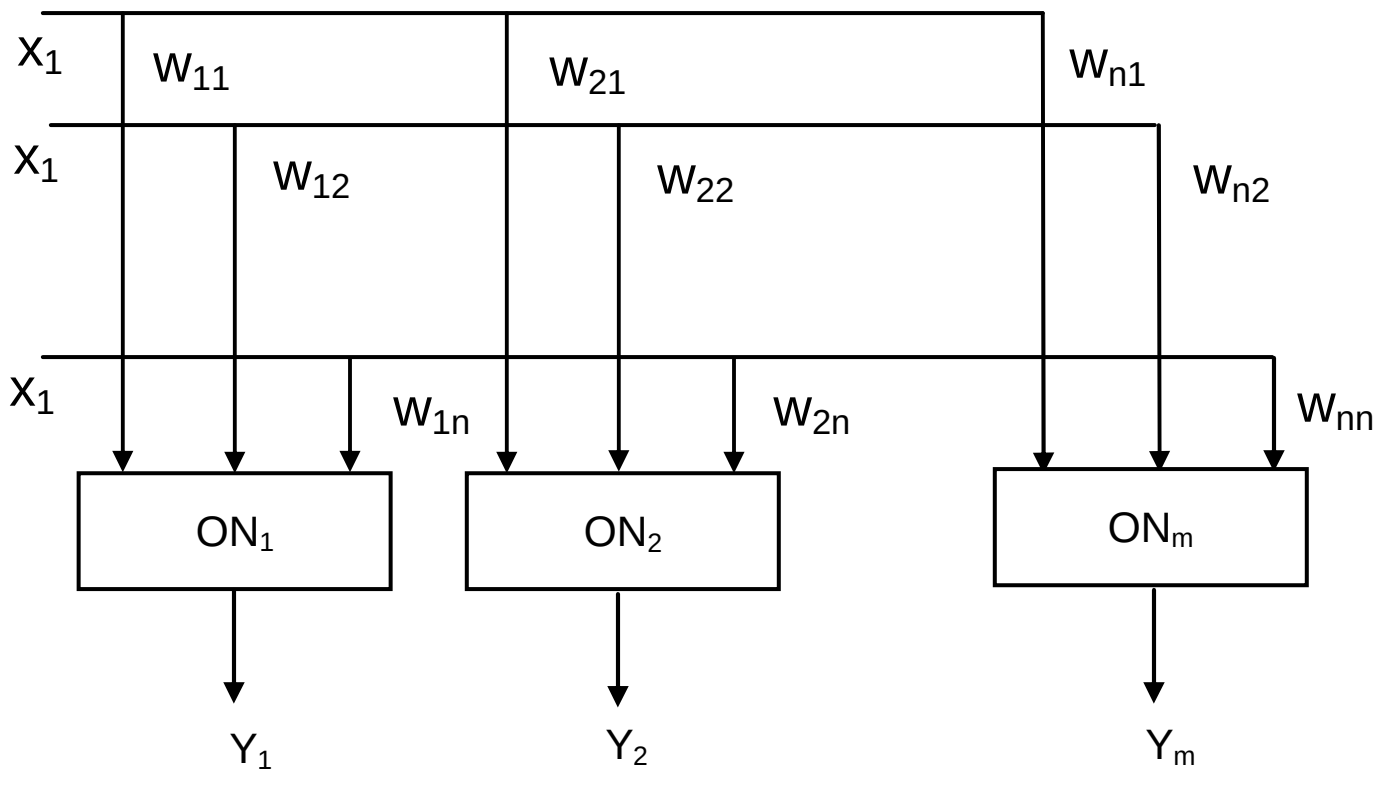

Figure: 4.1.1

The output $y_i = f(\max W_{ij}x_j)$, for all $i \in N_n$ and $j \in N_m$. Solution to (1) is obtained by varying the weights $W_{ij}$ so that the difference between the predicted value and the calculated value is zero.

***FRE to estimate flow rates in a chemical plants***

A typical chemical plant consists of several interlinked units. These units act as nodes. The flowsheet is given in Figure 4.1.2.

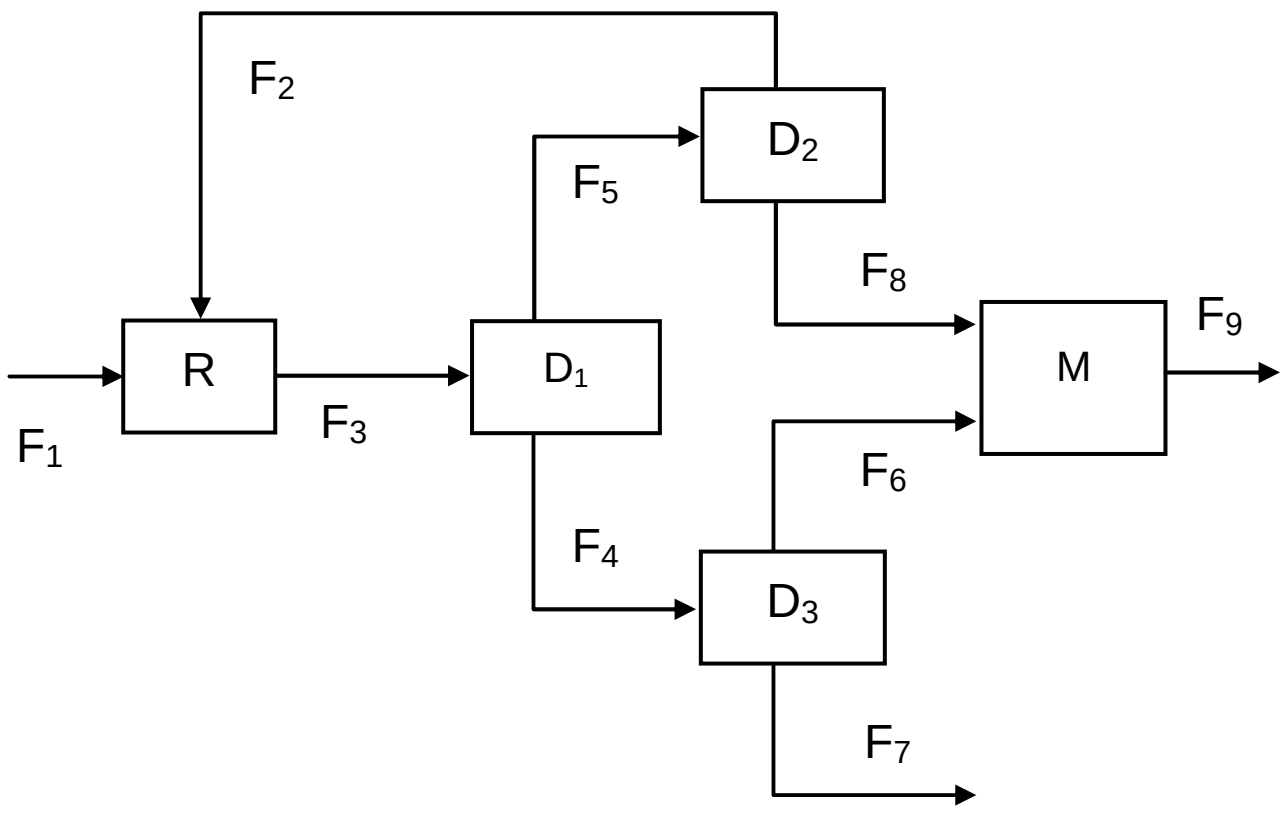

Figure: 4.1.2

An experimental approach would involve measuring the nine flow-rates to describe the state of the plant which would involve more money and labor.

While studying this problem in practice researchers have has neglected density variations across each stream. The mass balance equations across each node at steady state can be written as

$$\begin{aligned} F_3 - F_2 &= F_1, \\ F_2 - F_4 &= F_5, \\ F_4 - F_7 &= F_6, \\ F_2 + F_8 &= F_5, \\ F_8 &= F_9 - F_6. \end{aligned} \qquad (3)$$

Here $F_i$ represents the volumetric flow rate of the $i^{th}$ stream. In equation (3) at least four variables have to be specified or determined experimentally.

The remaining five can then be estimated from the equation (3), which is generated by applying the principle of conservation of mass to each unit. We assume $F_1$, $F_5$, $F_6$ and $F_9$ are experimentally measured, equation (3) reads with known values on the right-hand side as follows:

$$\begin{bmatrix} -1 & 1 & 0 & 0 & 0 \\ 0 & 1 & -1 & 0 & 0 \\ 0 & 0 & 1 & -1 & 0 \\ 1 & 0 & 0 & 0 & 1 \\ 0 & 0 & 0 & 0 & 1 \end{bmatrix} \begin{bmatrix} F_2 \\ F_3 \\ F_4 \\ F_7 \\ F_8 \end{bmatrix} = \begin{bmatrix} F_1 \\ F_5 \\ F_6 \\ F_5 \\ F_9 - F_6 \end{bmatrix} \quad (4)$$

$$P \circ Q = R \quad (5)$$

where, P, Q and R are explained. Using principle of conservation of mass balanced equation we estimate the flow rates of the five liquid stream. We in this problem aim to minimize the errors between the measured and the predicted value. We do this by giving suitable membership grades $p_{ij} \in [0, 1]$ and estimate the flow rates by using these $p_{ij}$'s in the equation 3. Now the equation 4 reads as follows:

$$\begin{bmatrix} p_{11} & p_{12} & 0 & 0 & 0 \\ 0 & p_{22} & p_{23} & 0 & 0 \\ 0 & 0 & p_{33} & p_{34} & 0 \\ p_{41} & 0 & 0 & 0 & p_{45} \\ 0 & 0 & 0 & 0 & p_{55} \end{bmatrix} \begin{bmatrix} F_2 \\ F_3 \\ F_4 \\ F_7 \\ F_8 \end{bmatrix} = \begin{bmatrix} F_1 \\ F_5 \\ F_6 \\ F_5 \\ F_9 - F_6 \end{bmatrix} \quad (6)$$

where $P = (p_{ij})$,
$Q = (q_{ik}) = [F_2\ F_3\ F_4\ F_7\ F_8]^t$ and
$R = (r_{ik}) = [F_1\ F_5\ F_6\ F_5\ F_9 - F_6]^t$.

We now apply the partitioning method of solution to equation (6). The partitioning of P correspondingly partitions R, which is give by a set of give subsets as follows:

$$[p_{11}\ p_{12}\ 0\ 0\ 0] \begin{bmatrix} F_2 \\ F_3 \\ F_4 \\ F_7 \\ F_8 \end{bmatrix} = \begin{bmatrix} F_1 \\ F_5 \\ F_9 \\ F_5 \\ F_9 - F_6 \end{bmatrix}, \ldots$$

$$[0\ 0\ \ 0\ \ 0\ p_{55}] \begin{bmatrix} F_2 \\ F_3 \\ F_4 \\ F_5 \\ F_8 \end{bmatrix} = \begin{bmatrix} F_1 \\ F_5 \\ F_6 \\ F_5 \\ F_9 - F_6 \end{bmatrix}.$$

Suppose the subsets satisfies the condition max $q_{ik} < r_{ik}$ then it has no solution. If it does not satisfy, this condition, then it has a final solution. If we have no solution we proceed to the second stage of solving the problem using Fuzzy Neural Networks.

When the FRE has no solution by the partition method, we solve these FRE using neural networks. This is done by giving weightages of zero elements as 0 and the modified FRE now reads as

$$P_1 \circ \begin{bmatrix} F_2 \\ F_3 \\ F_4 \\ F_7 \\ F_8 \end{bmatrix} = \begin{bmatrix} F_1 \\ F_5 \\ F_6 \\ F_5 \\ F_9 - F_6 \end{bmatrix}.$$

The linear activation function f defined earlier gives the output $y_i = f(\max W_{ij}\, x_j)$ $(i \in N_n)$ we calculate max $W_{ij}x_j$ as follows:

1. $W_{11}x_1 = 0.02F_2$, $W_{12}x_2 = 0F_2$, $W_{13}x_3 = 0F_2$ $W_{14}x_4 = 0.045F_2$, $W_{15}x_5 = 0F_2$

   $y_1 = f(\max_{j \in Nm} W_{ij}x_j) = f(0.02F_2, 0F_2, 0.045F_2, 0F_2)$

2. $W_{21}x_1 = 0.04F_3$, $W_{22}x_2 = 0.045F_3$, $W_{23}x_3 = 0F_3$, $W_{24}x_4 = 0.0F_3$, $W_{15}x_5 = 0F_3$

   $y_2 = f(\max_{j \in Nm} W_{ij}x_j) = f(0.04F_3, 0.045F_3, 0F_3, 0.0F_3, 0F_3)$

3. $W_{31}x_1 = 0.0F_4$, $W_{32}x_2 = 0.085F_4$, $W_{33}x_3 = 0.15F_4$, $W_{34}x_4 = 0.0F_4$ $W_{35}x_5 = 0F_4$

   $y_3 = f(\max_{j \in Nm} W_{ij}x_j) = f(0F_4, 0.085F_4, 0.15F_4, 0F_4, 0F_4)$

4. $W_{41}x_1 = 0.0F_7$, $W_{42}x_2 = 0F_7$, $W_{43}x_3 = 0.2F_7$, $W_{44}x_4 = 0.0F_7$, $W_{45}x_5 = 0F_7$

$$y_4 = f(\max_{j \in Nm} W_{ij}x_j) = f(0F_7, 0F_7, 0.2F_7, 0.0F_7, 0F_7)$$

5. $W_{51}x_1 = 0.0F_8$, $W_{52}x_2 = 0F_8$, $W_{53}x_3 = 0F_8$, $W_{54}x_4 = 0.45F_8$, $W_{55}x_5 = 0.5F_8$

$$y_5 = f(\max_{j \in Nm} W_{ij}x_j) = f(0F_8, 0F_8, 0F_8, 0.45F_8, 0.5F_8)$$

shown in Figure 4.1.2. Suppose the error does not reach 0 we change the weights till the error reaches 0. We continue the process again and again until the error reaches to zero.

Thus to reach the value zero we may have to go on giving different weightages (finite number of time) till say $s^{th}$ stage $P_s \circ Q_s$ whose linear activation function f, makes the predicted value to be equal to the calculated value. Thus by this method, we are guaranteed of a solution which coincides with the predicted value.

## 4.2 Fuzzy neural networks to estimate velocity of flow distribution in a pipe network

In flow distribution in a pipe network of a chemical plant, we consider liquid entering into a pipe of length T and diameter D at a fixed pressure $P_i$, The flow distributes itself into two pipes each of length $T_1(T_2)$ and diameter $D_1(D_2)$ given in Figure 4.2.1.

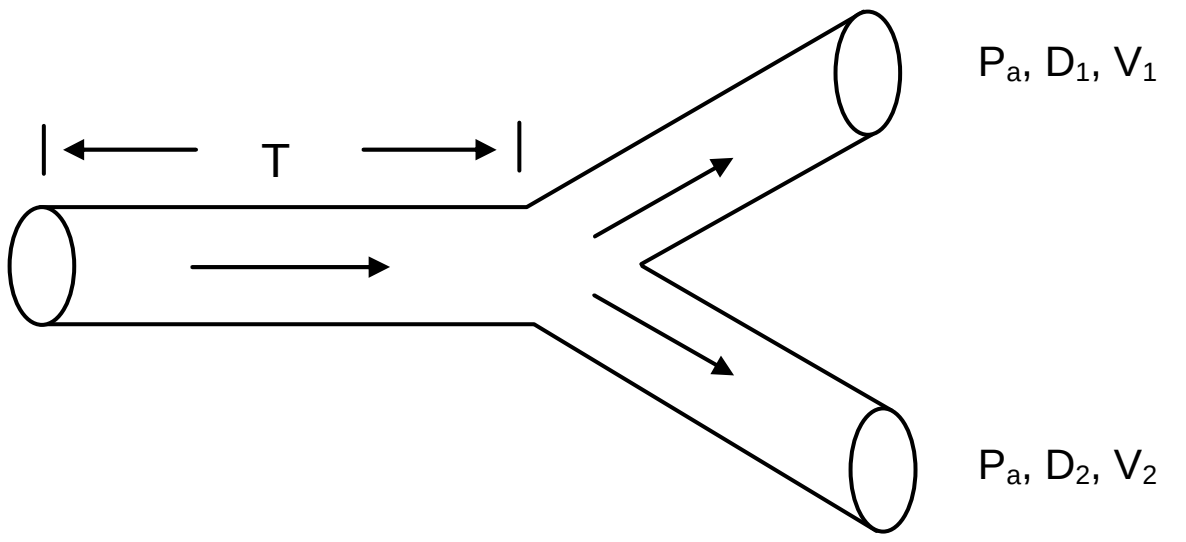

Figure: 4.2.1

The linear equation is based on Ohm's law, the drop in voltage V across a resistor R is given by the linear relation V =

iR (Ohm's law). The hydrodynamic analogue to the mean velocity v for laminar flow in a pipe is given by $\nabla_p$ = v $(32\mu T/D^2)$. This is classical-Poiseulle equation. In flow distribution in a pipe network, neglecting pressure losses at the junction and assuming the flow is laminar in each pipe, the macroscopic momentum balance and the mass balance at the junction yields,

$$P_1 - P_a = (32\mu T/D^2)v + (32\mu T_1 D_1^2)v_1,$$
$$P_i - P_a = (32\mu T/D^2)v + 32\mu T/D_2^2)v_2,$$
$$D^2 v = D_1^2 v_1 + v_2 D_2^2 \,. \qquad (1)$$

Hence $P_a$ is the pressure at which the fluid leaves the system at the two outlets. The set of three equation in (1) can be solved and we estimate v, $v_1$, $v_2$ for a fixed $(P_i - P_a)$. The system reads as

$$\begin{bmatrix} 32\mu T/D^2 & 32\mu T_1/D_1^2 & 0 \\ 32\mu T/D^2 & 0 & 32\mu T_2/D_2^2 \\ -D^2 & D_2^1 & D_2^2 \end{bmatrix} \begin{bmatrix} v \\ v_1 \\ v_2 \end{bmatrix} = \begin{bmatrix} p_i - p_a \\ p_i - p_a \\ 0 \end{bmatrix}.$$

We transform this equation into a fuzzy relation equation. We use a similar procedure described earlier and obtain the result by fuzzy relation equation. We get max (0.2v, 0.025v, 0.03v), max (0.035v, $0v_1$, $0.04v_1$), max ($0v_2$, $0.04v_2$, $0.045v_2$) by using neural networks for fuzzy relation equation described in [11]. Suppose the error does not reach to 0, we change the weights till the error reaches 0. We continue the process again and again till the error reaches zero.

### 4.3 Fuzzy neural networks to estimate three stage counter current extraction unit

Three-stage counter extraction unit is shown in Figure 4.3.1. The components A present in phase E (extract) along with a nondiffusing substance as being mixture.

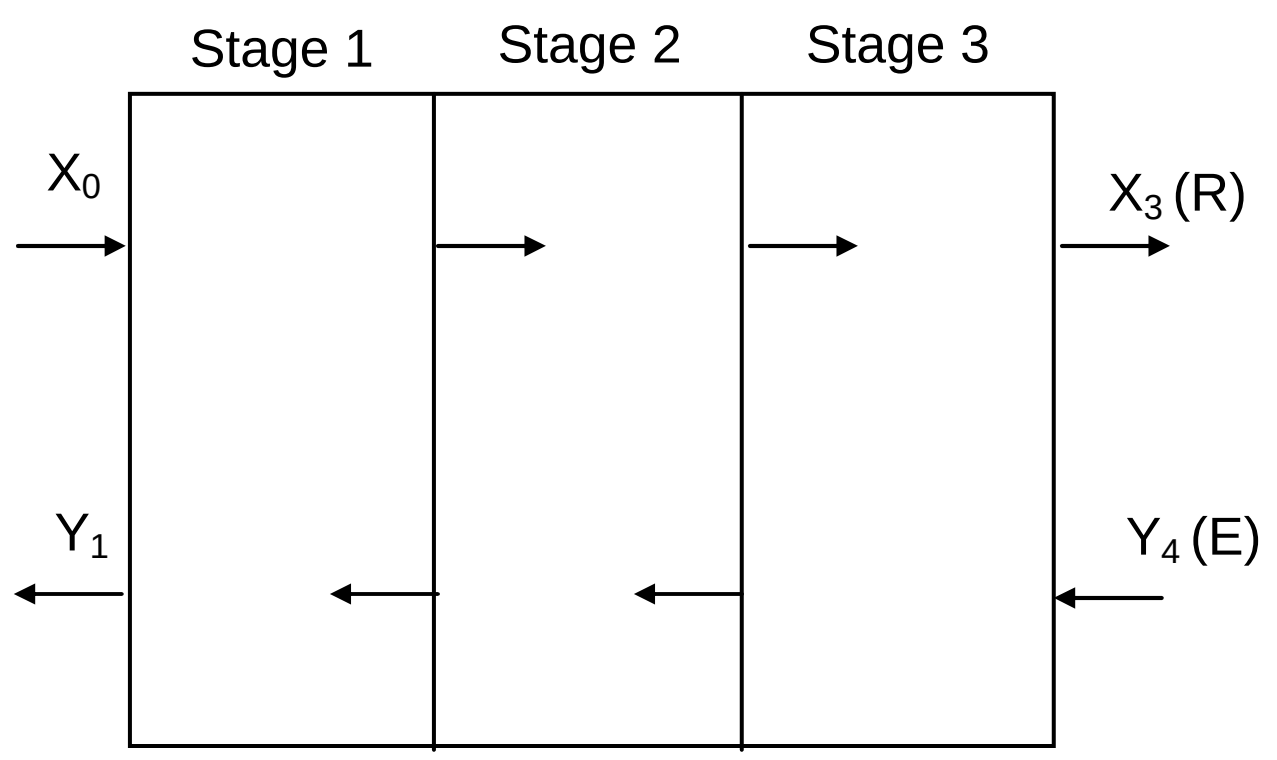

Figure: 4.3.1

*It is extracted into R by a nondiffusing solvent. The 3 extraction stage is given by the three equation.*

$$
\begin{aligned}
&E_sY_4 + R_sX_s = R_sX_3 + E_sY_3,\\
&EsY_3 + RsX_1 = E_s + R_sX_2,\\
&EsY_2 + R_sX_0 = E_sY_1 + R_sX_1 \qquad (1)
\end{aligned}
$$

$Y_i(X_i)$ = moles of A, The flow of each stage is denoted by $E_s(R_s)$ and this constant does not vary between the different stages. The assumption of a linear equilibrium relationship for the compositions leaving the $i^{th}$ stage equations

$$Y_i = KX_i \qquad (2)$$

for i = 1, 2, 3 reads as

$$
\begin{bmatrix}
R_s & E_s & 0 & -E_s & 0 & 0 \\
K & -1 & 0 & 0 & 0 & 0 \\
-R_s & 0 & R_s & E_s & 0 & -E_s \\
0 & 0 & K & -1 & 0 & 0 \\
0 & 0 & -R_2 & 0 & R_s & E_s \\
0 & 0 & 0 & 0 & K & -1
\end{bmatrix}
\begin{bmatrix}
X_1 \\ Y_1 \\ X_2 \\ Y_2 \\ X_3 \\ Y_3
\end{bmatrix}
=
\begin{bmatrix}
R_sX_0 \\ 0 \\ 0 \\ 0 \\ E_sY_4 \\ 0
\end{bmatrix}
$$

where $\{X_1, Y_1, X_2, Y_2, X_3, Y_3\}$ can be obtained for a given $E_s$, $R_s$ and K. Since use of linear algebraic equation does not result in the closeness of the measured and predicted value, we use neural networks for fuzzy relation equations to estimate the flow-rates of the stream, moles of the three-stage counter extraction unit and velocity of the flow distribution in a pipe network. As neural networks is a method to reduce the errors between the measured value and the predicted value. This allows varying degrees of set membership (weightages) based on a membership function defined over the range of value. The (weightages) membership function usually varies from 0 to 1. We use the similar modified procedure described earlier and get result by fuzzy relation equation. We get max $(0.2X_1, 0.25X_1, 0.3X_1, 0X_1, 0X_1, 0X_1)$, max $(0.35Y_1, 0.4Y_1, 0Y_1, 0Y_1, 0Y_1, 0Y_1)$ max $(0X_2, 0X_2, 0.45X_2, 0.5X_2, 0.55X_2, 0X_2)$, max $(0.6Y_2, 0Y_2, 0.65Y_5, 0.7Y_2, 0Y_2, 0Y_2)$ max $(0X_3, 0X_3, 0X_3, 0X_3, 0.75X_3, 0.8X_3)$, max $(0Y_3, 0Y_3, 0.85Y_3, 0Y_3, 0.9Y_3, 0.95Y_3)$ by neural networks for fuzzy relation equation. We continue this process until the error reaches zero or very near to zero.

Thus we see that when we replace algebraic linear equations by fuzzy methods to the problems described we are not only guaranteed of a solution, but our solution is very close to the predicted value.

**Chapter Five**

# MINIMIZATION OF WASTE GAS FLOW IN CHEMICAL INDUSTRIES

Chemical Industries and Automobiles are extensively contributing to the pollution of environment, Carbon monoxide, nitric oxide, ozone, etc., are understood as the some of the factors of pollution from chemical industries. The maintenance of clean and healthy atmosphere makes it necessary to keep the pollution under control which is caused by combustion waste gas. The authors have suggested theory to control waste gas pollution in environment by oil refinery using fuzzy linear programming. To the best of our knowledge the authors [43]are the first one to apply fuzzy linear programming to control or minimize waste gas in oil refinery.

An oil refinery consists of several inter linked units. These units act as production units, refinery units and compressors parts. These refinery units consume high-purity gas production units. But the gas production units produce high-purity gas along with a low purity gas. This low purity gas goes as a waste gas flow and this waste gas released in the atmosphere causes pollution in the environment. But in the oil refinery the quantity of this waste gas flow is an uncertainty varying with time and quality of chemicals used in the oil refinery. Since a complete eradication of waste gas in atmosphere cannot be made; here one aims to minimize the waste gas flow so that pollution in environment can be reduced to some extent. Generally waste gas flow is determined by linear programming method. In the study of minimizing the waste gas flow, some times the current state of the refinery may already be sufficiently close to the optimum. To over come this situation we adopt fuzzy linear programming method.

The fuzzy linear programming is defined by

$$\begin{aligned} &\text{Maximize} \quad z = cx \\ &\text{Such that} \quad Ax \leq b \\ &\qquad\qquad\quad\; x \leq 0 \end{aligned}$$

where the coefficients A, b and c are fuzzy numbers, the constraints may be considered as fuzzy inequalities with variables x and z. We use fuzzy linear programming to determine uncertainty of waste gas flow in oil refinery which pollutes the environment.

Oil that comes from the ground is called “Crude oil”. Before one can use it, oil has to be purified at a factory called a “refinery”, so as to convert into a fuel or a product for use. The refineries are high-tech factories, they turn crude oil into useful energy products. During the process of purification of crude oil in an oil refinery a large amount of waste gas is emitted to atmosphere which is dangerous to human life, wildlife and plant life. The pollutants can affect the health in various ways, by causing diseases such as bronchitis or asthma, contributing to cancer or birth defects or perhaps by damaging the body’s immune system which makes people more susceptible to a variety of other health risks. Mainly, this waste gas affects Ozone Layer. Ozone (or Ozone Layer) is 10-50 km above the surface of earth. Ozone provides a critical barrier to solar ultraviolet radiation, and protection from skin cancers, cataracts, and serious ecological disruption. Further sulfur dioxide and nitrogen oxide combine with water in the air to form sulfuric acid and nitric acid respectively, causing acid rain. It has been estimated that emission of 70 percentage of sulfur dioxide and nitrogen oxide are from chemical industries.

We cannot stop this process of oil refinery, since oil and natural gas are the main sources of energy. We cannot close down all oil refineries, but we only can try to control the amount of pollution to a possible degree. In this paper, the authors use fuzzy linear programming to reduce the waste gas from oil refinery. The authors describe the knowledge based system (KBS) that is designed and incorporate it in this paper to generate an on-line advice for operators regarding the proper distribution of gas resources in an oil refinery. In this system, there are many different sources of uncertainty including modeling errors, operating cost, and different opinions of experts on operating strategy. The KBS consists of sub-functions, like first sub-functions, second sub-functions, etc. Each and every sub-functions are discussed relative to certain specific problems.

For example: The first sub-function is mainly adopted to the compressor parts in the oil refineries. Till date they were using stochastic programming, flexibility analysis and process design problems for linear or non-linear problem to compressor parts in oil refinery. Here we adopt the sub function to study the proper distribution of gas resources in an oil refinery and also use fuzzy linear programming (FLP) to minimize the waste gas flow. By the term proper distribution of gas we include the study of both the production of high-purity gas as well as the amount of waste gas flow which causes pollution in environment.

In 1965, Lofti Zadeh [115, 116] wrote his famous paper formally defining multi-valued, or "fuzzy" set theory. He extended traditional set theory by changing the two-values indicator functions i.e., 0, 1 or the crisp function into a multi-valued membership function. The membership function assigns a "grade of membership" ranging from 0 to 1 to each object in the fuzzy set. Zadeh formally defined fuzzy sets, their properties, and various properties of algebraic fuzzy sets. He introduced the concept of linguistic variables which have values that are linguistic in nature (i.e. pollution by waste gas = {small pollution, high pollution, very high pollution}).

Fuzzy Linear Programming (FLP): FLP problems with fuzzy coefficients and fuzzy inequality relations as a multiple fuzzy reasoning scheme, where the past happening of the scheme correspond to the constraints of thee FLP problem. We assign facts (real data from industries) of the scheme, as the objective of the FLP problem. Then the solution process consists of two steps. In the fist step, for every decision variable, we compute the (fuzzy) value of the objective function via constraints and facts/objectives. At the second step an optimal solution to FLP problem is obtained at any point, which produces a maximal element to the set of objective functions (in the sense of the given inequality relation).

The Fuzzy Linear Programming (FLP) problem application is designed to offer advice to operating personnel regarding the distribution of Gas within an oil refinery (Described in Figure 5.1) in a way which would minimize the waste gas in environment there by reduce the atmospheric pollution .

GPUI, GPU2 and GPU3 are the gas production units and GGG consumes high purity gas and vents low purity gas. Gas from these production units are sent to some oil refinery units, like sulfur, methanol, etc. Any additional gas needs in the oil refinery must be met by the gas production unit GPU3.

The pressure swing adsorption unit (PSA) separates the GPU2 gas into a high purity product stream and a low purity tail stream (described in the Figure 5.1). $C_1$, $C_2$, $C_3$, $C_4$, $C_5$, are compressors. The flow lines that dead –end is an arrow represent vent to flare or fuel gas. This is the wasted gas that is to be minimized. Also we want to minimize the letdown flow from the high purity to the low purity header

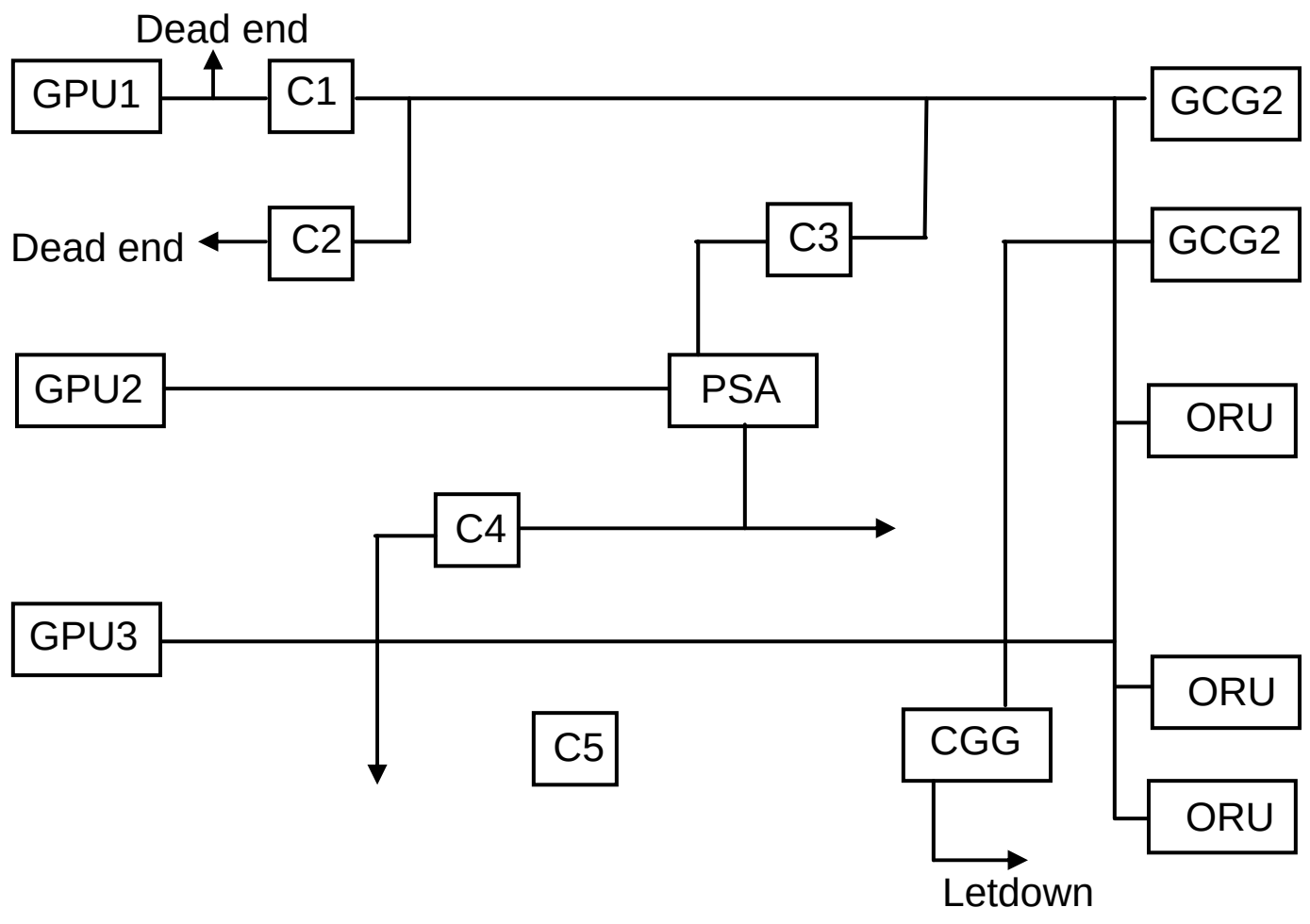

Figure: 5.1

FLP is a method of accounting for uncertainty is used by the authors for proper distribution of gas resources, so as to minimize the waste gas flow in atmosphere. FLP allows varying degrees of set membership based on a membership function defined over a range of values. The membership function usually varies from 0 to 1. FLP allow the representation of many different sources of uncertainty in the oil refinery. These sources may (or) may not be probabilistic in nature. The uncertainty is represented by membership functions describing the parameters in the optimization model. A solution is found that either maximizes a given feasibility measure and maximizes the wastage of gas flow. FLP is used here to characterize the neighborhood of solutions that defines the boundaries of acceptable operating states.

Fuzzy Linear Programming (FLP) can be stated as;

$$\left.\begin{array}{c} \max\text{imize } z = cx \\ \text{s.t } Ax \le b \\ x \ge 0 \end{array}\right] \qquad \ldots (*)$$

The coefficients A, b and c are fuzzy numbers, the constraints may be considered as fuzzy inequalities. The decision space is defined by the constraints with c, $x \in N$, $b \in R^m$ and $A \in R^m$, where N, $R^m$, and $R^{mxn}$ are reals.

The optimization model chosen by the knowledge based system (KBS) is determined online and is dependent on the refinery units. This optimization method is to reduce the amount of waste gas in pollution.

We aim to

1. The gas ($GCG_2$) vent should be minimized.
2. The let down flow should be minimized and
3. The make up gas produced by the as production unit (GPU3) should be minimized.

Generally the waste gas emitted by the above three ways pollute the environment. The objective function can be expressed as the sum of the individual gas waste flows. The constrains are given by some physical limitations as well as operator entries that describe minimum and maximum desired flows.

The obtained or calculated resultant values of the decision variables are interpreted as changes in the pressure swing adsorption feed, and the rate that gas is imported to CGG and gas production unit (GPU3). But in the optimization model there is uncertainty associated with amount of waste gas from oil refinery, and also some times the current state of the refinery may already be sufficiently close to the optimum.

For example to illustrate the problem, if the fuzzy constraints $x_1$, the objects are taken along the x-axis are shown in the figures 5.2 and 5.3, which represent the expression.

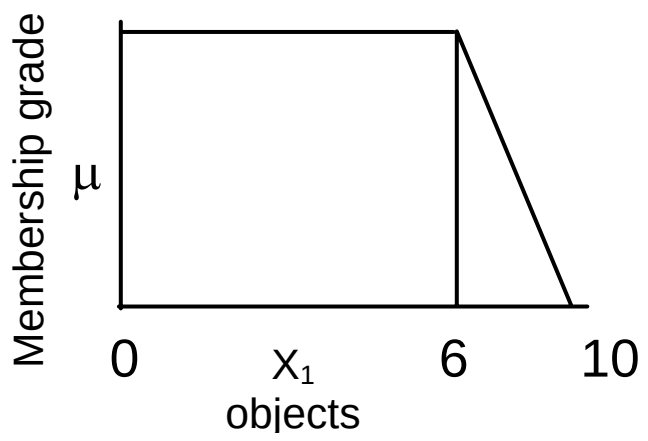

Figure: 5.2

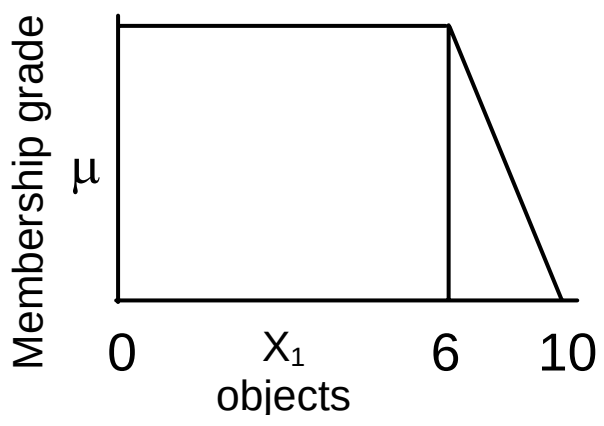

Figure: 5.3

$$x_1 \leq 8 \text{ (with tolerance } p = 2) \qquad (1)$$

The membership function μ are taken along the y-axis i.e. $\mu(x_1)$ lies in [0, 1] this can be interpreted as the confidence with which this constraint is satisfied (0 for low and 1 for high). The fuzzy inequality constraints can be redefined in terms of their α-cuts.

$$\{S_\alpha / \alpha \ \varepsilon \ [0, 1]\}, \text{ where } S_\alpha = \{\gamma / (\mu(\gamma) \geq \alpha)\}.$$

The parameter α is used to turn fuzzy inequalities into crisp inequalities. So we can rewrite equation (1)

$$x_1 \leq 6 + 2\,(2)\,(1 - \alpha)$$
$$x_1 \leq 6 + 4\,(1 - \alpha)$$

where $\alpha$ $\varepsilon$ [0, 1] expressed in terms of $\alpha$ in this way the fuzzy linear programming problem can be solved parametrically. The solution is a function on $\alpha$

$$x^* = f(\alpha) \qquad (2)$$

with the optimal value of the objective function determined by substitution in equation (1).

$$z^* = cx^* = g(\alpha). \qquad (3)$$

This is used to characterize the objective function. The result covers all possible solutions to the optimization problem for any point in the uncertain interval of the constraints.

The $\alpha$-cuts of the fuzzy set describes the region of feasible solutions in figures 5.2 and 5.3. The extremes ($\alpha = 0$ and $\alpha = 1$) are associated with the minimum and maximum values of $x^*$ respectively. The given equation (2) can also be found this, is used to characterize the objective function. The result covers all possible solutions to the optimization problem for any point in the uncertain interval of the constraints.

Fuzzy Membership Function to Describe Uncertainty: The feasibility of any decision ($\mu_D$) is given by the intersection of the fuzzy set describing the objective and the constraints.

$$\mu_D (x) = \mu_z(x) \wedge \mu_N (x)$$

where ^ represents the minimum operator, that is the usual operation for fuzzy set intersection. The value of $\mu_N$ can be easily found by intersecting the membership values for each of the constraints.

$$\mu_N (x) = \mu_1(x) \wedge \mu_2 (x) \wedge \ldots \wedge \mu_m (x).$$

The membership functions for the objective ($\mu_z$) however is not obvious z is defined in (2). Often, predetermined aspiration target values are used to define this function. Since reasonable values of this kind may not be available, the solution to the FLP equation (3) is used to characterize this function.

$$\mu_z(x) = \begin{cases} 1 & \text{if } z(x) \geq b(0) \\ \dfrac{z(x) - b(1)}{b(0) - b(1)} & \text{if } b(1) \leq z(x) \leq b(0) \\ 0 & \text{if } z(x) \leq b(0). \end{cases} \qquad (5)$$

The result is that the confidence value increases as the value of the objective value increases. This is reasonable because the goal is to maximize this function the limits on the function defined by reasonable value is obtained by extremes of the objective value.

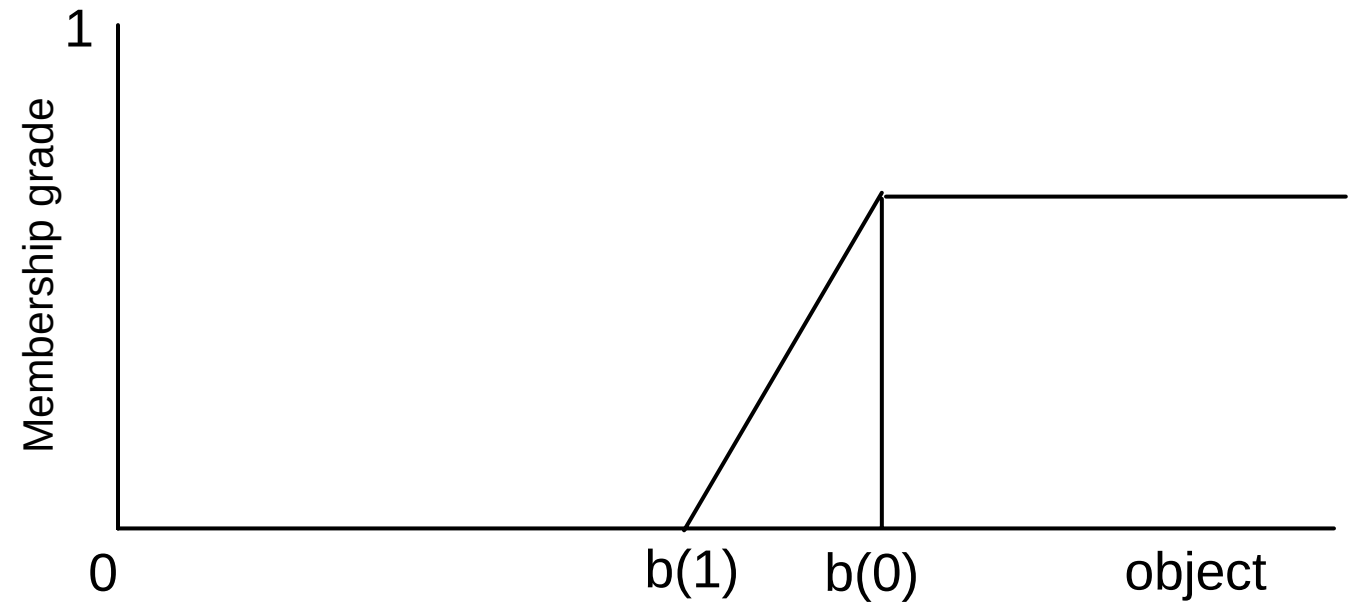

Figure: 5.4

These are the results generated by the fuzzy linear programming. Since both $\mu_N$ and $\mu_z$ have been characterized, now our goal is to describe the appropriateness of any operation state. Given any operating x, the feasibility can be specified based on the objective value, the constraints and the estimated uncertainty is got using equation (4). The value of $\mu_D$ are shown as the intersection of the two membership functions.

Defining the decision region based on the intersection we describe the variables and constraints of our problem. The variable $x_1$ represents the amount of gas fed to pressure swing adsorption from the gas production unit. The variable $x_2$ represents the amount of gas production that is sent to CGG. This problem can be represented according to equation (*). The constraints on the problem are subjected to some degree of uncertainty often some violation of the constraints within this range of uncertainty is tolerable. This problem can be represented according to equation (*). Using the given refinery data from the chemical plant.

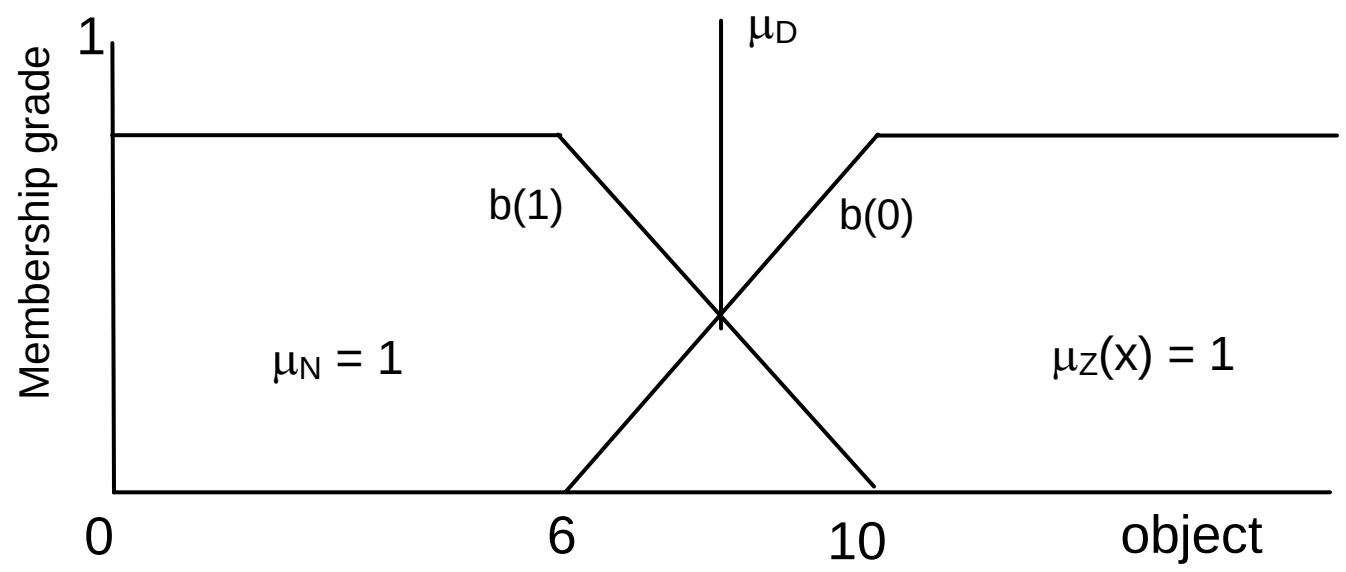

Figure: 5.5

c = [-0.544 3]

$$A = \begin{bmatrix} 1 & 0 \\ 0 & 1 \\ -0.544 & 1 \end{bmatrix},$$

$$b = \begin{bmatrix} 33.652 \\ 23.050 \\ 4.743 \end{bmatrix}$$

Using equation (*) we get
$Zc = -0.544\ x_1 + 3x_2$ it represents gas waste flow. The gas waste flow is represented by the following three equations:

i. $x_1 + 0x_2 \le 33.652$ is the total dead – end waste flow gas.
ii. $0x_1 + x_2 = 23.050$ is the total (GCG2) gas consuming gas – treaters waste flow gas.
iii. $-\ 0544\ x_1 + x_2 \le 4.743$ is the total let-down waste flow gas.

All flow rates are in million standard cubic feet per day. (i.e. 1 MMSCFD = 0.3277 $m^3/s$ at STP). The value used for may be considered to be desired from operator experts opinion. The third constraint represents the minimum let-down flow receiving to keep valve from sticking. The value to this limit cannot be given an exact value, therefore a certain degree of violation may be tolerable. The other constraints may be subject to some uncertainty as well as they represent the maximum allowable values for $x_1$ and $x_2$. In this problem we are going to express all constrains in terms of $\alpha$, $\alpha$, $\varepsilon$ [0, 1]. We have to chose a value of tolerance on the third constraint as $p_3 = 0.1$, then this constraint is represented parametrically as

$$a_3\, x \leq (b_3 - p_3) + 2p_3\,(1 - \alpha).$$

For example, if we use crisp optimization problem with the tolerance value p = 0.1 we obtain the following result:

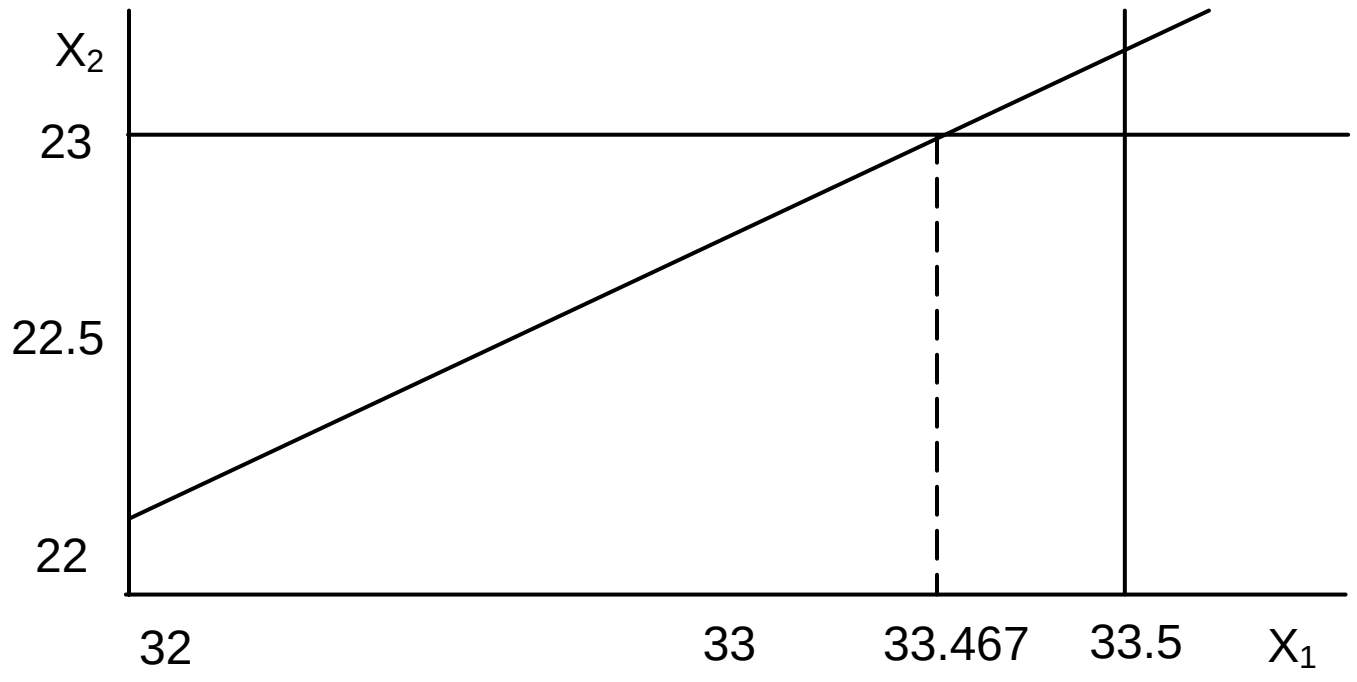

Figure: 5.6

where $x_1$ represents the amount of gas fed to PSA from gas production unit which is taken along the x axis, and $x_2$ amount of gas sent to CGG which is taken along the y axis,
we get $x_1 = 33.469$, when $x_2 = 23.050$

$$x^* = \begin{bmatrix} 33.469 \\ 23.050 \end{bmatrix}$$

z = 50.941. Finally we compare this result with our fuzzy linear programming method.

We replace two valued indicator function method by fuzzy linear programming.

Fuzzy Linear Programming is used now to maximize the objective function as well as minimize the uncertainty (waste flow gas). For that all of the constraints are expressed in terms of $\alpha$, $\alpha$, $\in$ [0, 1].

$$a_3\, x \le (b_3 - p_3) + 2p_3\,(1 - \alpha).\ \alpha \in [0, 1]$$

where $a_3$ is the third row in the matrix A. i.e. = $0.544x_1 + x_2 \le 4.843 - 0.2\,\alpha$, when the tolerance $p_3 = 0.3$, we fix the value of $\alpha$ $\varepsilon$ [0.9,1], when the tolerance $p_3 = 0.1$, we see $\alpha$ $\varepsilon$ [0.300, 0.600].

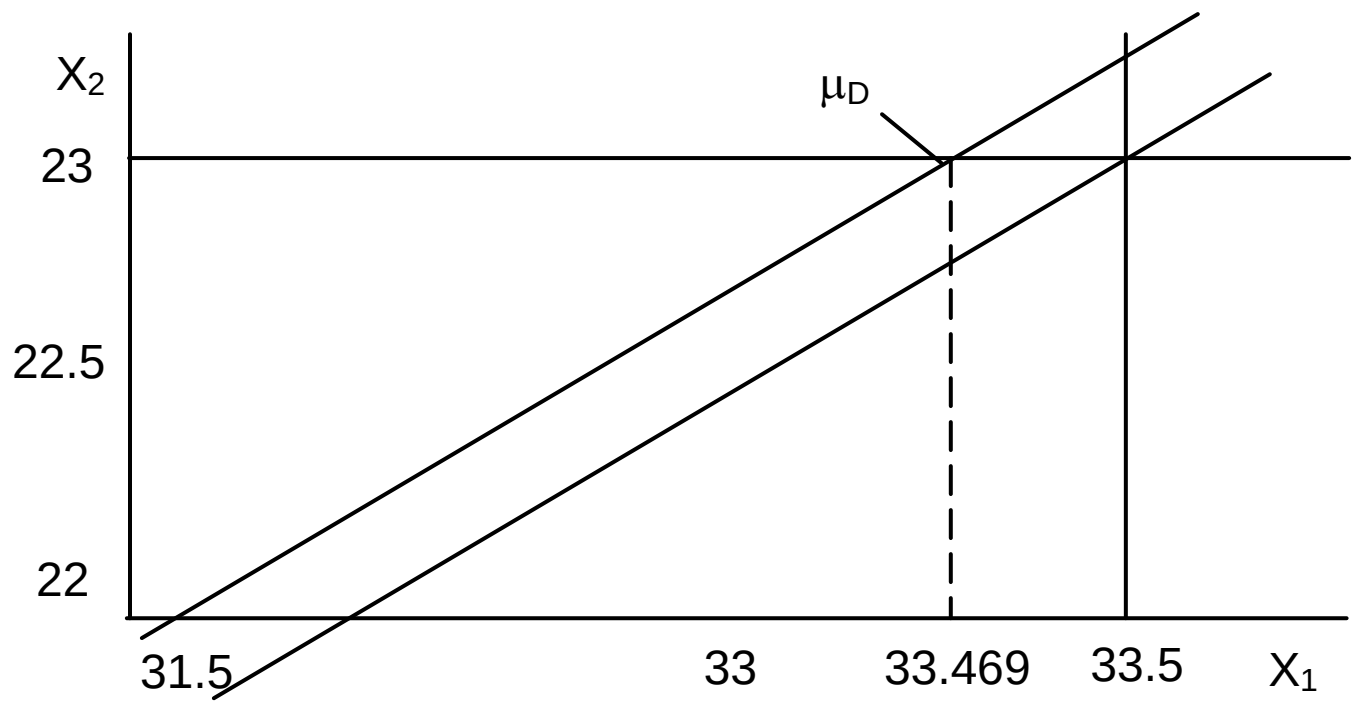

Figure: 5.7

where $x_1$ represents the amount of gas fed to PSA from gas production unit which is taken along the x axis, and $x_2$ amount of gas that is sent to CGG which is taken along the y axis,

When $x_2 = 23.050$ and $\alpha = 0.0$, we get $x_1 = 33.469$.

When $x_2 = 23.050$ and $\alpha = 0.4$, we get $x_1 = 33.616$

The set ($\mu_z$) is defined in equation 5. Fuzzy Linear Programming solution is

$$x^* = f(\alpha) = \begin{bmatrix} 33.469 \\ 23.050 \end{bmatrix}$$

this value is recommended as there is no changes in the operating policy.

So we have to chose the value for $\alpha$ as 0.6 for the tolerance $p_3 = 0.1$, we get the following graph where $x_1$ represents the amount of gas fed to PSA from gas production unit which is taken along the x axis, and $x_2$ amount of gas that sent to CGG which is taken along the y axis,

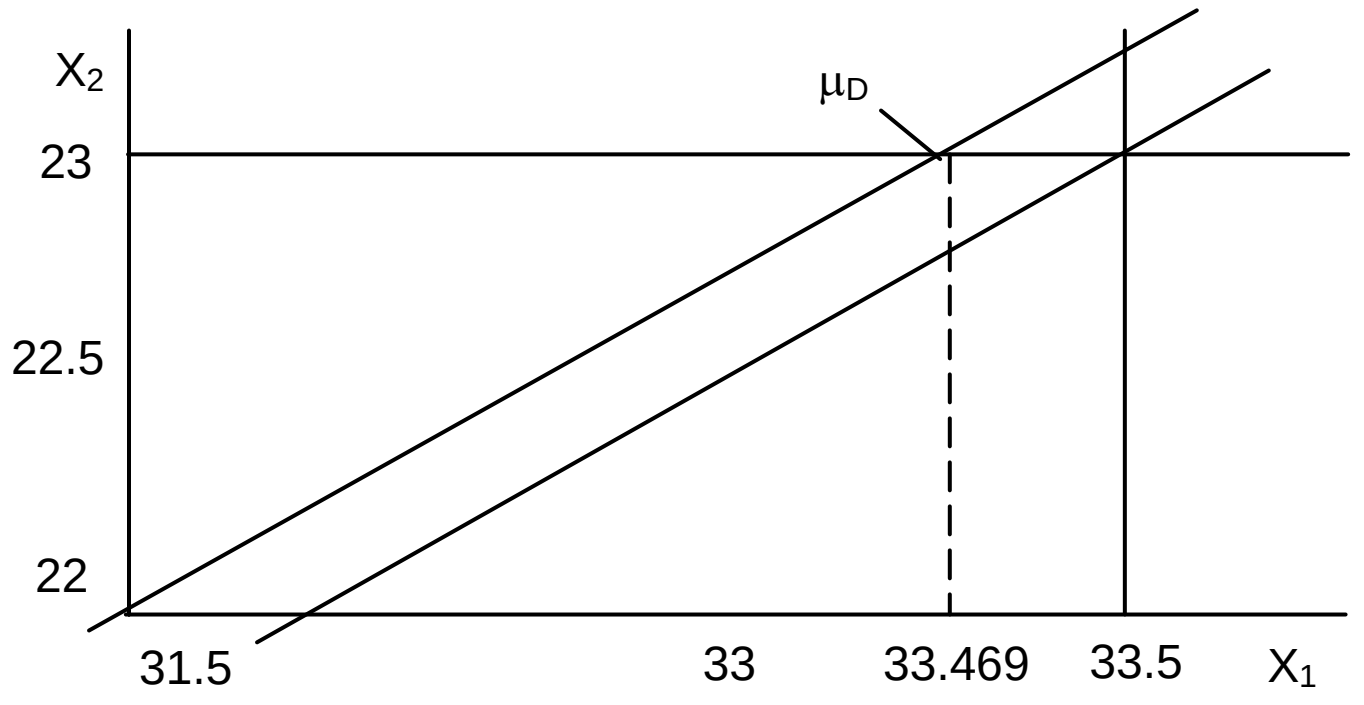

Figure: 5.8

when $x_2$ = 23.050 and $\alpha$ = 0.0 we get $x_1$ = 33.469

when $x_2$ = 23.050 and $\alpha$ = 0.6 we get $x_1$ = 33.689.

The operating region

$$x^* = f(0.6) = \begin{bmatrix} 33.689 \\ 23.050 \end{bmatrix}.$$

Now if the tolerance on the third constraint is increased to $p_3$ = 0.2. This results is the region shown in the following graph. As expected the region has increased to allow a larger range of operating states.

when $x_2$ = 23.050 and $\alpha$ = 0.0 we get $x_1$ = 33.285

when $x_2$ = 23.050 and $\alpha$ = 0.9 we get $x_1$ = 33.947.

The operating region is

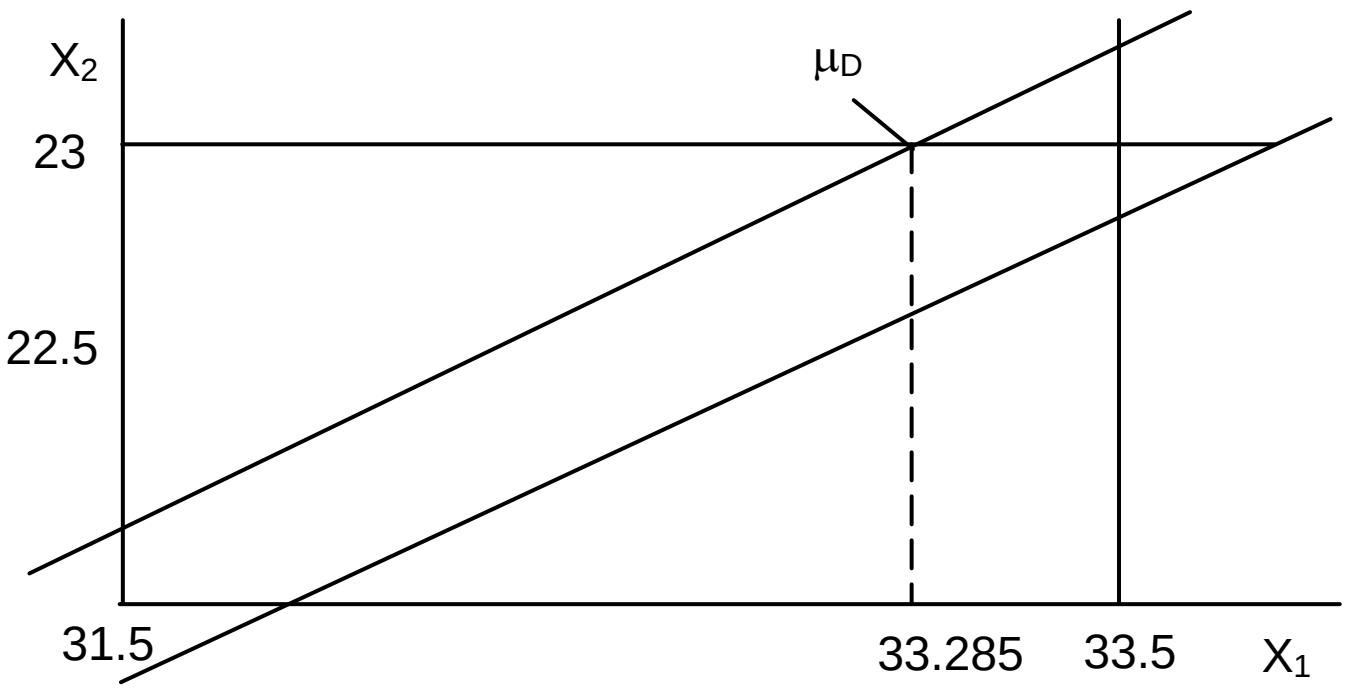

Figure: 5.9

where $x_1$ represents the amount of gas fed to PSA from gas production unit which is taken along the x axis, and $x_2$ amount of gas that is sent to CGG which is taken along the y axis.

$$x^* = f\,(0.9) = \begin{bmatrix} 33.947 \\ 23.050 \end{bmatrix}.$$

The fuzzy linear programming solution is

$$x^* = f\,(\alpha) = \begin{bmatrix} 33.285 \\ 23.050 \end{bmatrix}$$

$$z^* = 51.043.$$

Finally we have to take $\alpha$ ε [0.9, 1.00].
Choose $\alpha$ = 0 and when the tolerance $p_3$ = 0.3 we get the following graph when $x_2$ = 23.050 we get $x_1$ = 33.101.

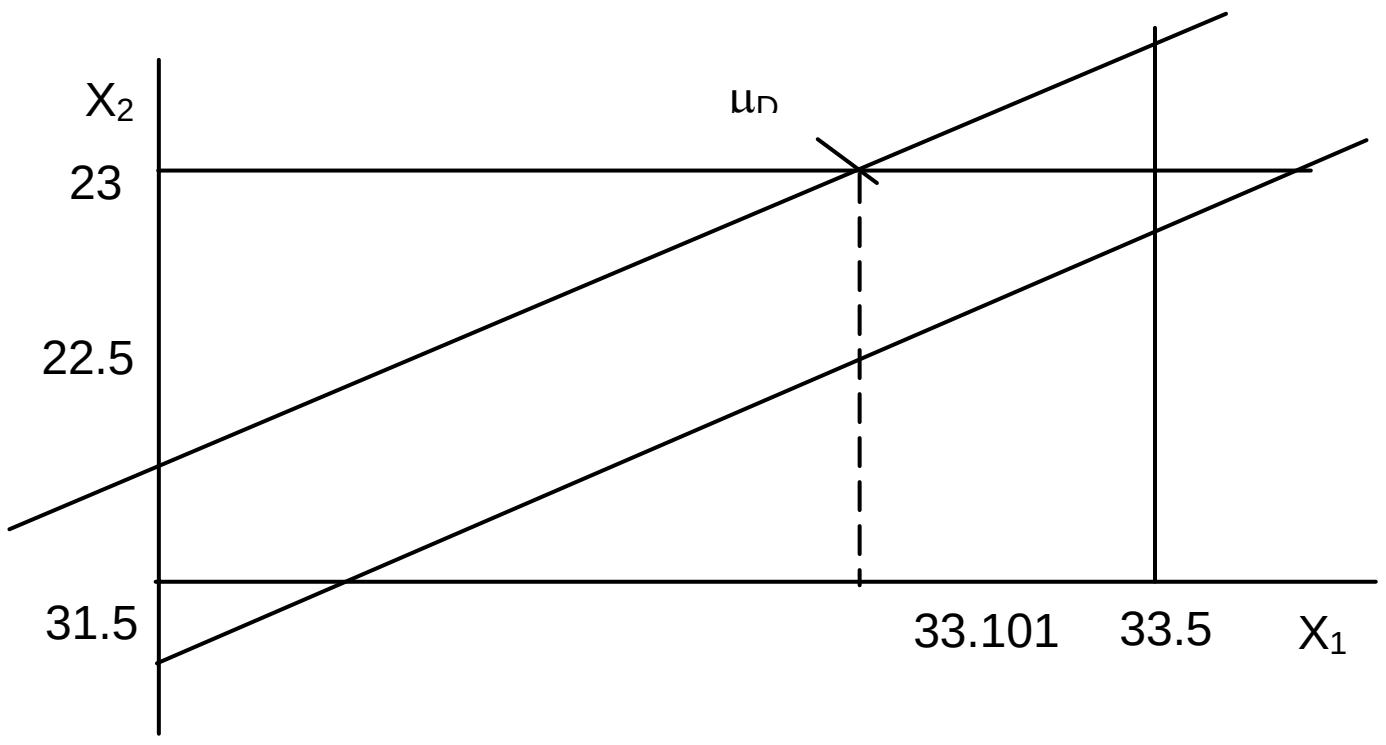

Figure: 5.10

where $x_1$ represents the amount of gas fed to PSA from gas production unit; and $x_2$ amount of gas that is sent to CGG.

When $\alpha$ = 1 and $x_2$ = 23.050 we get $x_1$ = 34.204. The operating region is

$$x^* = f\,(1.0) = \begin{bmatrix} 33.204 \\ 23.050 \end{bmatrix}.$$

The fuzzy linear programming solutions are

$$x^* = f(\alpha) = \begin{bmatrix} 33.101 \\ 23.050 \end{bmatrix}.$$

The fuzzy linear programming solutions are

$$z^* = g(\alpha) = 51.143.$$

We chose maximum value from the Fuzzy Linear Programming method i.e. $z^* = 51.143$.

Thus when we work by giving varying membership functions and use fuzzy linear programming we see that we get the minimized waste gas flow value as 33.101 in contrast to 33.464 measured in million standard cubic feet per day and the maximum gas waste flow of system is determined to be 51.143 in contrast to their result of 50.941 measured in million standard cubic feet per day. Since the difference we have obtained is certainly significant, this study when applied to any oil refinery will minimize the waste gas flow to atmosphere considerably and reduce the pollution.

**Chapter Six**

# USE OF NEUTROSOPHIC RELATIONAL EQUATIONS IN CHEMICAL ENGINEERING

This chapter has 2 sections. Section one gives introduction to Neutrosophic Relational Equations (NRE) and section two gives use of NRE in chemical Engineering.

## 6.1 Introduction to Neutrosophic Relation and their properties

In this section we introduce the notion of neutrosophic relational equations and fuzzy neutrosophic relational equations and analyze and apply them to real-world problems, which are abundant with the concept of indeterminacy. We also mention that most of the unsupervised data also involve at least to certain degrees the notion of indeterminacy.

Throughout this section by a neutrosophic matrix we mean a matrix whose entries are from the set N = $[0, 1] \cup I$ and by a fuzzy neutrosophic matrix we mean a matrix whose entries are from N' = $[0, 1] \cup \{nI / n \in (0,1]\}$.

Now we proceed on to define binary neutrosophic relations and binary neutrosophic fuzzy relation.

A binary neutrosophic relation $R_N(x, y)$ may assign to each element of X two or more elements of Y or the indeterminate *I*. Some basic operations on functions such as the inverse and composition are applicable to binary relations as well. Given a neutrosophic relation $R_N(X, Y)$ its domain is a neutrosophic set on $X \cup I$ domain R whose membership function is defined by

$$\text{domR}(x) = \max_{y \in X \cup I} R_N(x, y)$$

for each $x \in X \cup I$.

That is each element of set $X \cup I$ belongs to the domain of R to the degree equal to the strength of its strongest relation to any member of set $Y \cup I$. The degree may be an indeterminate *I* also. Thus this is one of the marked difference between the binary fuzzy relation and the binary neutrosophic relation. The range of $R_N(X,Y)$ is a neutrosophic relation on Y, ran R whose membership is defined by

$$\text{ran } R(y) = \max_{x \in X} R_N(x, y)$$

for each $y \in Y$, that is the strength of the strongest relation that each element of Y has to an element of X is equal to the degree of that element's membership in the range of R or it can be an indeterminate *I*.

The height of a neutrosophic relation $R_N(x, y)$ is a number h(R) or an indeterminate *I* defined by

$$h_N(R) = \max_{y \in Y \cup I} \max_{x \in X \cup I} R_N(x, y).$$

That is $h_N(R)$ is the largest membership grade attained by any pair (x, y) in R or the indeterminate *I*.

A convenient representation of the neutrosophic binary relation $R_N(X, Y)$ are membership matrices $R = [\gamma_{xy}]$ where $\gamma_{xy} \in R_N(x, y)$.

Another useful representation of a binary neutrosophic relation is a neutrosophic sagittal diagram. Each of the sets X, Y represented by a set of nodes in the diagram, nodes corresponding to one set are clearly distinguished from nodes representing the other set. Elements of $X' \times Y'$ with non-zero membership grades in $R_N(X, Y)$ are represented in the diagram by lines connecting the respective nodes. These lines are labeled with the values of the membership grades.

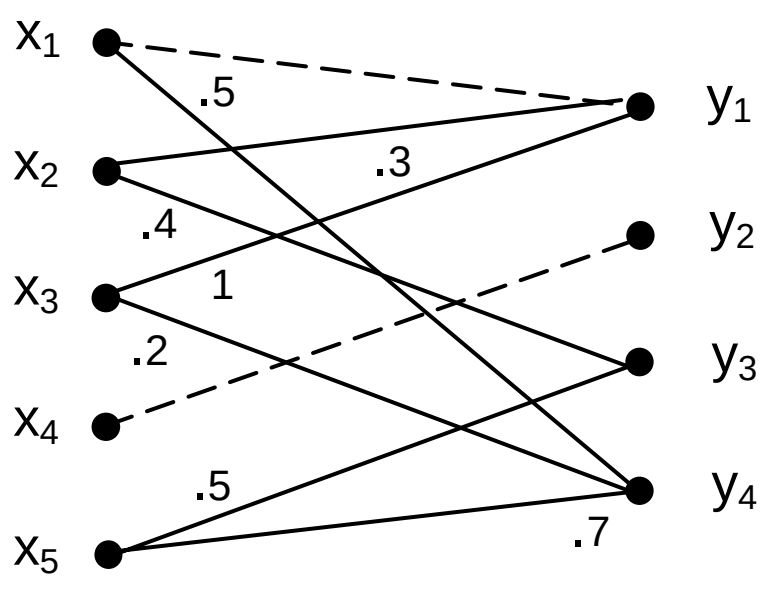

FIGURE: 6.1.1

An example of the neutrosophic sagittal diagram is a binary neutrosophic relation $R_N(X, Y)$ together with the membership neutrosophic matrix is given below.

$$\begin{array}{c} \\ x_1 \\ x_2 \\ x_3 \\ x_4 \\ x_5 \end{array} \begin{array}{c} \begin{array}{cccc} y_1 & y_2 & y_3 & y_4 \end{array} \\ \begin{bmatrix} I & 0 & 0 & 0.5 \\ 0.3 & 0 & 0.4 & 0 \\ 1 & 0 & 0 & 0.2 \\ 0 & I & 0 & 0 \\ 0 & 0 & 0.5 & 0.7 \end{bmatrix} \end{array}$$

The inverse of a neutrosophic relation $R_N(X, Y) = R(x, y)$ for all $x \in X$ and all $y \in Y$. A neutrosophic membership matrix $R^{-1} = [r^{-1}_{yx}]$ representing $R^{-1}_N(Y, X)$ is the transpose of the matrix R for $R_N(X, Y)$ which means that the rows of $R^{-1}$ equal the columns of R and the columns of $R^{-1}$ equal rows of R. Clearly $(R^{-1})^{-1} = R$ for any binary neutrosophic relation.

Consider any two binary neutrosophic relations $P_N(X, Y)$ and $Q_N(Y, Z)$ with a common set Y. The standard composition of these relations which is denoted by $P_N(X, Y) \bullet Q_N(Y, Z)$ produces a binary neutrosophic relation $R_N(X, Z)$ on $X \times Z$ defined by

$$R_N(x, z) = [P \bullet Q]_N(x, z) = \max_{y \in Y} \min[P_N(x, y), Q_N(x, y)]$$

for all $x \in X$ and all $z \in Z$.

This composition which is based on the standard $t_N$-norm and $t_N$-co-norm, is often referred to as the max-min composition. It can be easily verified that even in the case of binary neutrosophic relations

$$[P_N(X, Y) \bullet Q_N(Y, Z)]^{-1}$$
$$= Q_N^{-1}(Z, Y) \bullet P_N^{-1}(Y, X). [P_N(X, Y) \bullet Q_N(Y, Z)] \bullet R_N(Z, W)$$
$$= P_N(X, Y) \bullet [Q_N(Y, Z) \bullet R_N(Z, W)],$$

that is, the standard (or max-min) composition is associative and its inverse is equal to the reverse composition of the inverse relation. However, the standard composition is not commutative, because $Q_N(Y, Z) \bullet P_N(X, Y)$ is not well defined when $X \neq Z$. Even if $X = Z$ and $Q_N (Y, Z) \circ P_N (X, Y)$ are well defined still we can have $P_N (X, Y) \circ Q (Y, Z) \neq Q (Y, Z) \circ P (X, Y)$.

Compositions of binary neutrosophic relation can the performed conveniently in terms of membership matrices of the relations. Let $P = [p_{ik}]$, $Q = [q_{kj} ]$ and $R = [r_{ij}]$ be membership matrices of binary relations such that $R = P \circ Q$. We write this using matrix notation

$$[r_{ij}] = [p_{ik}] \text{ o } [q_{kj}]$$

where $r_{ij} = \max_k \min (p_{ik}, q_{kj})$.

A similar operation on two binary relations, which differs from the composition in that it yields triples instead of pairs, is known as the relational join. For neutrosophic relation $P_N (X, Y)$ and $Q_N (Y, Z)$ the relational join $P * Q$ corresponding to the neutrosophic standard max-min composition is a ternary relation $R_N (X, Y, Z)$ defined by $R_N (x, y, z) = [P * Q]_N (x, y, z) = \min [P_N (x, y), Q_N (y, z)]$ for each $x \in X$, $y \in Y$ and $z \in Z$.

This is illustrated by the following Figure 6.1.2.

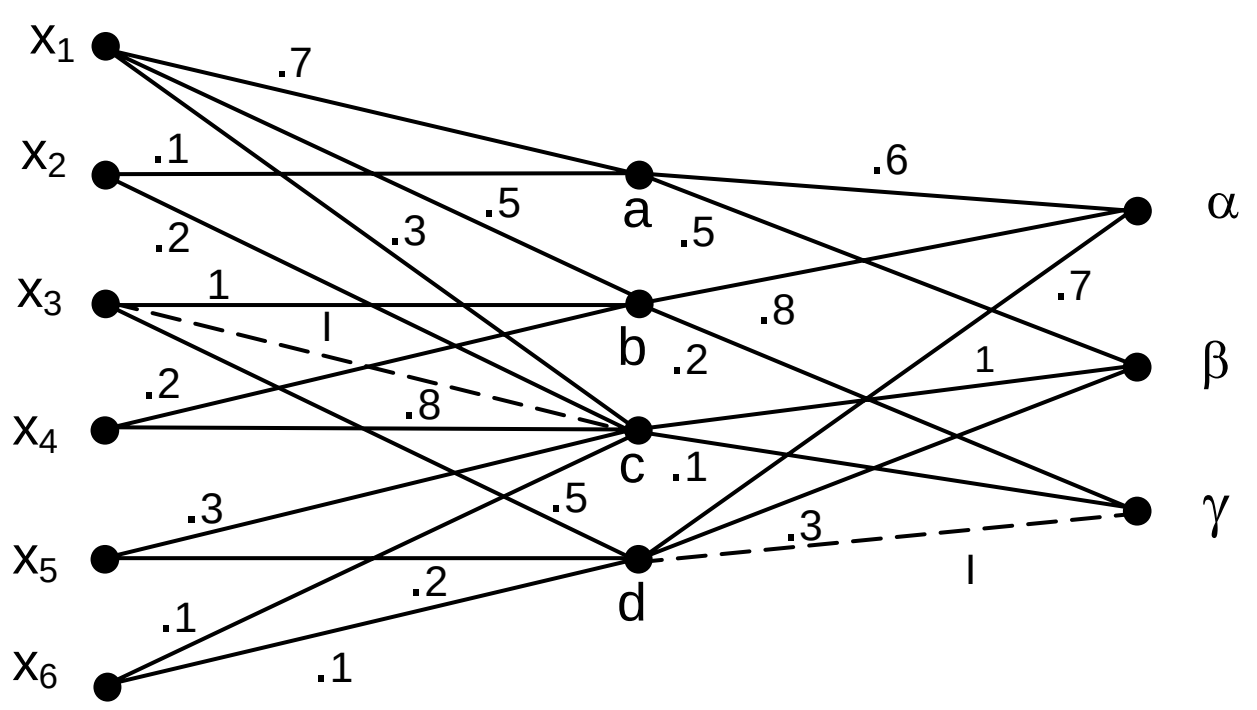

FIGURE: 6.1.2

In addition to defining a neutrosophic binary relation there exists between two different sets, it is also possible to define neutrosophic binary relation among the elements of a single set X. A neutrosophic binary relation of this type is denoted by $R_N(X, X)$ or $R_N(X^2)$ and is a subset of $X \times X = X^2$.

These relations are often referred to as neutrosophic directed graphs or neutrosophic digraphs. [42]

Neutrosophic binary relations $R_N$ (X, X) can be expressed by the same forms as general neutrosophic binary relations. However they can be conveniently expressed in terms of simple diagrams with the following properties.

I. Each element of the set X is represented by a single node in the diagram.
II. Directed connections between nodes indicate pairs of elements of X for which the grade of membership in R is non zero or indeterminate.
III. Each connection in the diagram is labeled by the actual membership grade of the corresponding pair in R or in indeterminacy of the relationship between those pairs.

The neutrosophic membership matrix and the neutrosophic sagittal diagram is as follows for any set X = {a, b, c, d, e}.

$$\begin{bmatrix} 0 & I & .3 & .2 & 0 \\ 1 & 0 & I & 0 & .3 \\ I & .2 & 0 & 0 & 0 \\ 0 & .6 & 0 & .3 & I \\ 0 & 0 & 0 & I & .2 \end{bmatrix}$$

Neutrosophic membership matrix for x is given above and the neutrosophic sagittal diagram is given below.

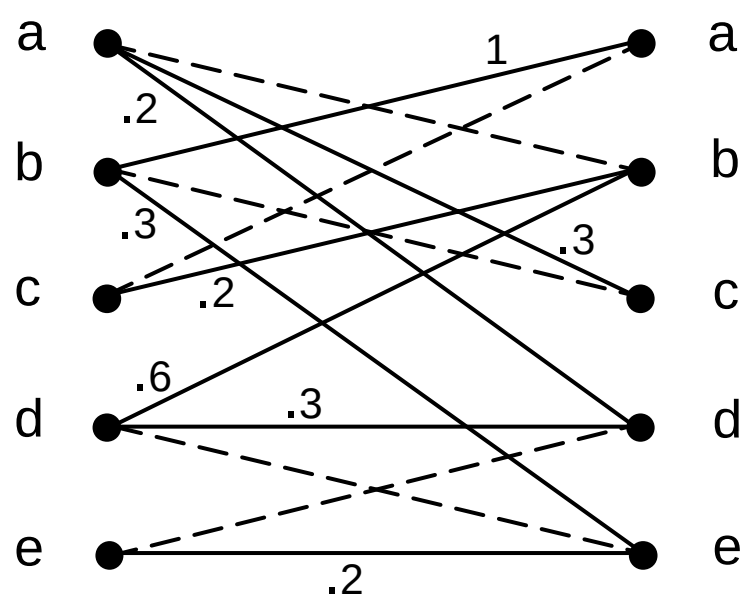

Figure: 4.1.3

Neutrosophic diagram or graph is left for the reader as an exercise.

The notion of reflexivity, symmetry and transitivity can be extended for neutrosophic relations $R_N$ (X, Y) by defining them in terms of the membership functions or indeterminacy relation.

*Thus $R_N$ (X, X) is reflexive if and only if $R_N$ (x, x) = 1 for all x $\in$ X.* If *this is not the case for some x $\in$ X the relation is irreflexive.*

*A weaker form of reflexivity, if for no x in X; $R_N$(x, x) = 1 then we call the relation to be anti-reflexive referred to as $\in$-reflexivity, is sometimes defined by requiring that*

$$R_N (x, x) \geq \in \text{ where } 0 < \in < 1.$$

*A fuzzy relation is symmetric if and only if*

$$R_N (x, y) = R_N (y, x) \text{ for all } x, y, \in X.$$

*Whenever this relation is not true for some x, y* $\in X$ *the relation is called asymmetric. Furthermore when* $R_N (x, y) > 0$ *and* $R_N (y, x) > 0$ *implies that x = y for all x, y* $\in X$ *the relation R is called anti-symmetric.*

*A fuzzy relation* $R_N (X, X)$ *is transitive (or more specifically max-min transitive) if*

$$R_N (x, z) \geq \max_{y \in Y} min [R_N (x, y), R_N (y, z)]$$

*is satisfied for each pair* $(x, z) \in X^2$*. A relation failing to satisfy the above inequality for some members of X is called non-transitive and if* $R_N (x, x) < \max_{y \in Y} min [R_N(x, y), R_N(y, z)]$ *for all* $(x, x) \in X^2$*, then the relation is called anti-transitive*

*Given a relation* $R_N(X, X)$ *its transitive closure* $\overline{R}_{NT} (x, X)$ *can be analyzed in the following way.*

The transitive closure on a crisp relation $R_N$ (X, X) is defined as the relation that is transitive, contains

$$R_N (X, X) < \max_{y \in Y} \min [R_N (x, y) R_N (y, z)]$$

for all $(x, x) \in X^2$, then the relation is called anti-transitive. Given a relation $R_N$ (x, x) its transitive closure $\overline{R}_{NT}$ (X, X) can be analyzed in the following way.

The transitive closure on a crisp relation $R_N$ (X, X) is defined as the relation that is transitive, contains $R_N$ and has the fewest possible members. For neutrosophic relations the last requirement is generalized such that the elements of transitive closure have the smallest possible membership grades, that still allow the first two requirements to be met.

Given a relation $R_N$ (X, X) its transitive closure $\overline{R}_{NT}$ (X, X) can be determined by a simple algorithm.

Now we proceed on to define the notion of neutrosophic equivalence relation.

**DEFINITION 4.1.1:** *A crisp neutrosophic relation* $R_N(X, X)$ *that is reflexive, symmetric and transitive is called a neutrosophic*

*equivalence relation. For each element x in X, we can define a crisp neutrosophic set $A_x$ which contains all the elements of X that are related to x by the neutrosophic equivalence relation.*

*Formally $A_x$ = [ y | (x, y) ∈ $R_N$ (X, X)}. $A_x$ is clearly a subset of X. The element x is itself contained in $A_x$, due to the reflexivity of R because R is transitive and symmetric each member of $A_x$ is related to all other members of $A_x$. Further no member of $A_x$ is related to any element of X not included in $A_x$. This set $A_x$ is clearly referred to as a neutrosophic equivalence class of $R_N$ (X, x) with respect to x. The members of each neutrosophic equivalence class can be considered neutrosophic equivalent to each other and only to each other under the relation R.*

But here it is pertinent to mention that in some X; (a, b) may not be related at all to be more precise there may be an element a ∈ X which is such that its relation with several or some elements in X \ {a} is an indeterminate. The elements which cannot determine its relation with other elements will be put in as separate set.

A neutrosophic binary relation that is reflexive, symmetric and transitive is known as a neutrosophic equivalence relation.

Now we proceed on to define Neutrosophic intersections neutrosophic t – norms ($t_N$ – norms)

Let A and B be any two neutrosophic sets, the intersection of A and B is specified in general by a neutrosophic binary operation on the set N = [0, 1] ∪ *I*, that is a function of the form

$$i_N: N \times N \to N.$$

For each element x of the universal set, this function takes as its argument the pair consisting of the elements membership grades in set A and in set B, and yield the membership grade of the element in the set constituting the intersection of A and B. Thus,

$$(A \cap B)(x) = i_N [A(x), B(x)] \text{ for all } x \in X.$$

In order for the function $i_N$ of this form to qualify as a fuzzy intersection, it must possess appropriate properties, which ensure that neutrosophic sets produced by $i_N$ are intuitively acceptable as meaningful fuzzy intersections of any given pair of neutrosophic sets. It turns out that functions known as $t_N$-norms, will be introduced and analyzed in this section. In fact the class of $t_N$- norms is now accepted as equivalent to the class of neutrosophic fuzzy intersection. We will use the terms $t_N$ – norms and neutrosophic intersections inter changeably.

Given a $t_N$ – norm, $i_N$ and neutrosophic sets A and B we have to apply:

$$(A \cap B)(x) = i_N [A(x), B(x)]$$

for each $x \in X$, to determine the intersection of A and B based upon $i_N$.

However the function $i_N$ is totally independent of x, it depends only on the values A (x) and B(x). Thus we may ignore x and assume that the arguments of $i_N$ are arbitrary numbers $a, b \in [0, 1] \cup I = N$ in the following examination of formal properties of $t_N$-norm.

A neutrosophic intersection/ $t_N$-norm $i_N$ is a binary operation on the unit interval that satisfies at least the following axioms for all $a, b, c, d \in N = [0, 1] \cup I$.

$1_N$ $\quad i_N(a, 1) = a$

$2_N$ $\quad i_N(a, I) = I$

$3_N$ $\quad b \leq d$ implies $i_N(a, b) \leq i_N(a, d)$

$4_N$ $\quad i_N(a, b) = i_N(b, a)$

$5_N$ $\quad i_N(a, i_N(b, d)) = i_N(a, b), d)$.

We call the conditions $1_N$ to $5_N$ as the axiomatic skeleton for neutrosophic intersections / $t_N$ – norms. Clearly $i_N$ is a continuous function on $N \setminus I$ and $i_N(a, a) < a \ \forall a \in N \setminus I$

$$i_N(I\ I) = I.$$

If $a_1 < a_2$ and $b_1 < b_2$ implies $i_N(a_1, b_1) < i_N(a_2, b_2)$. Several properties in this direction can be derived as in case of t-norms.

The following are some examples of $t_N$ –norms

1. $i_N$ (a, b) = min (a, b)
   $i_N$ (a, *I*) = min (a, *I*) = *I* will be called as standard neutrosophic intersection.
2. $i_N$ (a, b) = ab for a, b ∈ N \ *I* and $i_N$ (a, b) = *I* for a, b ∈ N where a = *I* or b = *I* will be called as the neutrosophic algebraic product.
3. Bounded neutrosophic difference.
   $i_N$ (a, b) = max (0, a + b – 1) for a, b ∈ *I*
   $i_N$ (a, *I*) = *I* is yet another example of $t_N$ – norm.
   1. Drastic neutrosophic intersection
   2.

$$i_N(a, b) = \begin{cases} a \text{ when } b = 1 \\ b \text{ when } a = 1 \\ I \text{ when } a = I \\ \quad \text{or } b = I \\ \quad \text{or } a = b = I \\ 0 \quad \text{otherwise} \end{cases}$$

As *I* is an indeterminate adjoined in $t_N$ – norms. It is not easy to give then the graphs of neutrosophic intersections. Here also we leave the analysis and study of these $t_N$ – norms (i.e. neutrosophic intersections) to the reader.

The notion of neutrosophic unions closely parallels that of neutrosophic intersections. Like neutrosophic intersection the general neutrosophic union of two neutrosophic sets A and B is specified by a function

$$\mu_N: N \times N \to N \text{ where } N = [0, 1] \cup I.$$

The argument of this function is the pair consisting of the membership grade of some element x in the neutrosophic set A and the membership grade of that some element in the neutrosophic set B, (here by membership grade we mean not only the membership grade in the unit interval [0, 1] but also the indeterminacy of the membership). The function returns the membership grade of the element in the set A ∪ B.

Thus (A ∪ B) (x) = $\mu_N$ [A (x), B(x)] for all x ∈ X. Properties that a function $\mu_N$ must satisfy to be initiatively acceptable as neutrosophic union are exactly the same as properties of functions that are known. Thus neutrosophic union will be called as neutrosophic t-co-norm; denoted by $t_N$ – co-norm.

A neutrosophic union / $t_N$ – co-norm $\mu_N$ is a binary operation on the unit interval that satisfies at least the following conditions for all a, b, c, d $\in$ N = [0, 1] $\cup$ *I*

$C_1$ $\mu_N (a, 0) = a$

$C_2$ $\mu_N (a, I) = I$

$C_3$ $b \leq d$ implies

$\mu_N (a, b) \leq \mu_N (a, d)$

$C_4$ $\mu_N (a, b) = \mu_N (b, a)$

$C_5$ $\mu_N (a, \mu_N (b, d))$

$= \mu_N (\mu_N (a, b), d).$

Since the above set of conditions are essentially neutrosophic unions we call it the axiomatic skeleton for neutrosophic unions / $t_N$-co-norms. The addition requirements for neutrosophic unions are

i. $\mu_N$ is a continuous functions on N \ {*I*}

ii. $\mu_N (a, a) > a$.

iii. $a_1 < a_2$ and $b_1 < b_2$ implies $\mu_N (a_1. b_1) < \mu_N (a_2, b_2)$; $a_1, a_2, b_1, b_2 \in N \setminus \{I\}$

We give some basic neutrosophic unions.

Let $\mu_N : [0, 1] \times [0, 1] \to [0, 1]$

$\mu_N (a, b) = \max (a, b)$

$\mu_N (a, I) = I$ is called as the standard neutrosophic union.

$\mu_N (a, b) = a + b - ab$ and

$\mu_N (a, I) = I$ .

This function will be called as the neutrosophic algebraic sum.

$$\mu_N (a, b) = \min (1, a + b) \text{ and } \mu_N (a, I) = I$$

will be called as the neutrosophic bounded sum. We define the notion of neutrosophic drastic unions

$$\mu_N(a, b) = \begin{cases} a \text{ when } b = 0 \\ b \text{ when } a = 0 \\ I \text{ when } a = I \\ \quad \text{or } b = I \\ 1 \text{ otherwise.} \end{cases}$$

Now we proceed on to define the notion of neutrosophic Aggregation operators. Neutrosophic aggregation operators on neutrosophic sets are operations by which several neutrosophic sets are combined in a desirable way to produce a single neutrosophic set.

Any neutrosophic aggregation operation on n neutrosophic sets ($n \geq 2$) is defined by a function $h_N: N^n \to N$ where $N = [0, 1] \cup I$ and $N^n = \underbrace{N \times \ldots \times N}_{n-\text{times}}$ when applied to neutrosophic sets $A_1, A_2, \ldots, A_n$ defined on X the function $h_N$ produces an aggregate neutrosophic set A by operating on the membership grades of these sets for each $x \in X$ (Here also by the term membership grades we mean not only the membership grades from the unit interval [0, 1] but also the indeterminacy *I* for some $x \in X$ are included). Thus

$$A_N(x) = h_N(A_1(x), A_2(x), \ldots, A_n(x))$$

for each $x \in X$.

### 6.2 Use of NRE in chemical engineering

The use of FRE for the first line has been used in the study of flow rates in chemical plants. In this study we are only guaranteed of a solution but when we use NRE in study of flow rates in the chemical plants we are also made to understand that certain flow rates are indeterminates depending on the leakage, chemical reactions and the new effect due to chemical reactions which may change due to change in the density/ viscosity of the fluid under study their by changing the flow rates while analyzing as a mathematical model. So in the study of flow rates in chemical plants some indeterminacy are also related with it. FRE has its own limitation for it cannot involve in its analysis the indeterminacy factor.

We have given analysis in chapter 2 using FRE. Now we suggest the use of NRE and bring out its importance in the determination of flow rates in chemical plants.

Consider the binary neutrosophic relations $P_N$ (X, Y) $Q_N$ (Y, Z) and R (X, Z) which are defined on the sets X, Y and Z. Let the membership matrices of P, Q and R be denoted by P = $[p_{ij}]$, Q = $[q_{jk}]$ and R = $[r_{ij}]$ respectively where $p_{ij} = P(x_i, y_j)$, $q_{jk} = Q(y_j, r_k)$ and $r_{ik} = R(x_i, z_k)$ for $i \in I = N_n$, $j \in J = N_m$ and $k \in K = N_k$ entries of P, Q and R are taken for the interval $[0\ 1] \times FN$. The three neutrosophic matrices constrain each other by the equation

$$P \circ Q = R \qquad (1)$$

where 'o' denotes the max-min composition (1) known as the Neutrosophic Relational Equation (NRE) which represents the set of equation

$$\max p_{ij}\, q_{jk} = r_{ik}. \qquad (2)$$

For all $i \in N_n$ and $k \in N_s$. After partitioning the matrix and solving the equation (1) yields maximum of $q_{jk} < r_{ik}$ for some $q_{jk}$, then this set of equation has no solution so to solve equation (2) we invent and redefine a feed – forward neural networks of one layer with n-neurons with m inputs. The inputs are associated with $w_{ij}$ called weights, which may be real, or indeterminates from R*I*. The neutrosophic activation function $f_N$ is defined by

$$f_N(a) = \begin{cases} 0 \text{ if } a < 0 \\ a \text{ if } a \in [01] \\ 1 \text{ if } a > 1 \\ aI \text{ if } a \in FN \\ I \text{ if } aI > I \\ 0 \text{ if in } aI,\ a < 0. \end{cases}$$

The out put $y_i = f_N$ (max $w_{ij}$ $x_j$). Now the NRE is used to estimate the flow rates in a chemical plant. In places where the indeterminacy is involved the expert can be very careful and use methods to overcome indeterminacy by adopting more and more constraints which have not been given proper representation and their by finding means to eliminate the indeterminacy involved in the weights.

In case of impossibility to eliminate these indeterminacy one can use the maximum caution in dealing with these values which are indeterminates so that all types of economic and time loss can be met with great care. In the flow rate problem the use of NRE mainly predicts the presence of the indeterminacy which can be minimized using $f_N$; where by all other in-descripancies are given due representation.

We suggest the use of NRE for when flow rates are concerned in any chemical plant the due weightage must be given the quality of chemicals or raw materials which in many cases are not up to expectations, leakage of pipe, the viscosity or density after chemical reaction time factor, which is related with time temperature and pressure for which usually due representations, is not given only ideal conditions are assumed. Thus use of NRE may prevent accident, economic loss and other conditions and so on.

# REFERENCE

1. Alexander H., *Biodegradation of Chemicals of Environmental concern*, Science, 211, (1981), 132-8

2. Atkinson, B. and Mavituna, F., *Biochemical Engineering and Biotechnology Handbook*, Macmiliam, 1983.

3. Blanco, A., Delgado, M., and Requena, I., Solving Fuzzy Relational Equations by max-min Neural Network, *Proc. 3rd IEEE Internet Conf. On Fuzzy Systems*, Orlando (1994) 1737-1742.

4. Buckley, J.J., and Hayashi, Y., Fuzzy Neural Networks: A Survey, *Fuzzy Sets and Systems*, 66 (1994) 1-13.

5. Chung, F., and Lee, T., A New Look at Solving a System of Fuzzy Relational Equations, *Fuzzy Sets and Systems*, 99 (1997) 343-353.

6. Di Nola, A., Pedrycz, W., Sessa, S., and Sanchez, E., Fuzzy Relation Equations Theory as a Basis of Fuzzy Modeling: An Overview, *Fuzzy Sets and Systems*, 40 (1991) 415-429.

7. Di Nola, A., Pedrycz, W., Sessa, S., and Wang, P.Z., Fuzzy Relation Equation under a Class of Triangular Norms: A Survey and New Results, *Stochastica*, 8 (1984) 99-145.

8. Di Nola, A., Sessa, S., Pedrycz, W., and Sanchez, E., *Fuzzy Relational Equations and their Application in Knowledge Engineering*, Kluwer Academic Publishers, Dordrecht, 1989.

9. Driankov, D, Hellendoorn, H. and Reinfrank, M., *An introduction to fuzzy control*, Narosa Publising houses New Delhi, 1997.

10. Ezio, B. and Y.Kung. Some properties of singular value decomposition and their applications to digital signal processing, *Elsevier Science Publishers* 18 (1989) 277 – 289

11. George, J. Klir and Bo Yuan, *Fuzzy Sets and Fuzzy logic*, Prentice Hall of India, New Delhi (1997).

12. Grant Kessler, R Cement Kiln Dust (CKD) methods for reduction and control, *IEEE Transactions on industry applications*, 3 1 (2) (1995) 407 – 412.

13. Hans Peter Schmauder, *Methods in Biotechnology*, Taylor and Francis, 2002

14. Holland, J., *Adaptation in Natural and Artificial Systems*, The University of Michigan Press, Ann Arbor, 1975.

15. John E. Smith, *Biotechnology*, 3rd ed., Cambridge Univ. Press, 1996.

16. Kagei, S., Fuzzy Relational Equation with Defuzzification Algorithm for the Largest Solution, *Fuzzy Sets and Systems*, 123 (2001) 119-127.

17. King, P.P., Biotechnology an industrial view, *J. of chemical Technology and Biotechnology*, 32, (1982), 2-8.

18. Klema, V.C. and A.J. Laub, The singular value decomposition : Its computation and some applications, *IEEE Transactions on Automatic Control*, 25(2) (1980) 164 – 176.

19. Kosko, B., *Neural Networks and Fuzzy Systems: Dynamical Approach to Machine Intelligence*, Prentice-Hall, Englewood Cliffs NJ, 1992.

20. Kreft. W., The interruption of material cycles taking account of integrated further utilization in the cement plant, *Zement Kalk Gips*, (1987) 257 – 259.

21. Lars E. Ebbesen, Main fractionate crude switch control, *Comput. Chem. Eng.*, 16, (1992), S165-S171.

22. Li, X., and Ruan, D., Novel Neural Algorithm Based on Fuzzy S-rules for Solving Fuzzy Relation Equations Part I, *Fuzzy Sets and Systems*, 90 (1997) 11-23.

23. Li, X., and Ruan, D., Novel Neural Algorithms Based on Fuzzy S-rules for Solving Fuzzy Relation Equations Part II, *Fuzzy Sets and Systems*, 103 (1999) 473-486.

24. Liu, F., and Smarandache, F., *Intentionally and Unintentionally. On Both, A and Non-A, in Neutrosophy*. http://lanl.arxiv.org/ftp/math/papers/0201/0201009.pdf

25. Lu, J., An Expert System Based on Fuzzy Relation Equations for PCS-1900 Cellular System Models, *Proc. South-eastern INFORMS Conference*, Myrtle Beach SC, Oct 1998.

26. Mizumoto, M., and Zimmermann, H.J., Comparisons of Fuzzy Reasoning Methods, *Fuzzy Sets and Systems*, 8 (1982) 253-283.

27. Neundorf, D., and Bohm, R., Solvability Criteria for Systems of Fuzzy Relation Equations, *Fuzzy Sets and Systems*, 80 (1996) 345-352.

28. Pavlica, V., and Petrovacki, D., About Simple Fuzzy Control and Fuzzy Control Based on Fuzzy Relational Equations, *Fuzzy Sets and Systems*, 101 (1999) 41-47.

29. Pedrycz, W, Algorithms for solving fuzzy relational equations in a probabilistic setting, *Fuzzy Sets and Systems*, 28 (1998) 183 – 202.

30. Pedrycz, W., *Fuzzy Control and Fuzzy Systems*, Wiley, New York, 1989.

31. Pedrycz, W., Inverse Problem in Fuzzy Relational Equations, *Fuzzy Sets and Systems*, 36 (1990) 277-291.

32. Pedrycz, W., Processing in Relational Structures: Fuzzy Relational Equations, *Fuzzy Sets and Systems*, 25 (1991) 77-106.

33. Pedrycz, W., s-t Fuzzy Relational Equations, *Fuzzy Sets and Systems*, 59 (1993) 189-195.

34. Ramathilagam, S., *Mathematical Approach to the Cement Industry problems using Fuzzy Theory*, Ph.D. Dissertation, Guide: Dr. W. B. Vasantha Kandasamy, Department of Mathematics, Indian Institute of Technology, Madras, November 2002.

35. Sanchez, E., Resolution of Composite Fuzzy Relation Equation, *Inform. and Control*, 30 (1976) 38-48.

36. Smarandache, F., *A Unifying Field in Logics: Neutrosophic Logic. Neutrosophy, Neutrosophic Set, Neutrosophic Probability and Statistics*, third edition, Xiquan, Phoenix, 2003.

37. Smarandache, F., *Collected Papers III*, Editura Abaddaba, Oradea, 2000. http://www.gallup.unm.edu/~smarandache/CP3.pdf

38. Smarandache, F., Definitions Derived from Neutrosophics, In Proceedings of the *First International Conference on Neutrosophy, Neutrosophic Logic, Neutrosophic Set, Neutrosophic Probability and Statistics*, University of New Mexico, Gallup, 1-3 December 2001.

39. Stamou, G.B., and Tzafestas, S.G., Neural Fuzzy Relational Systems with New Learning Algorithm, *Mathematics and Computers in Simulation*, 51 (2000) 301-314.

40. Stamou, G.B., and Tzafestas, S.G., Resolution of Composite Fuzzy Relation Equations based on Archimedean Triangular Norms, *Fuzzy Sets and Systems*, 120 (2001) 395-407.

41. Sugeno, M., *Industrial Applications of Fuzzy Control*, Elsevier, New York, 1985.

42. Vasantha Kandasamy, W.B., and Smarandache, F., *Fuzzy Cognitive Maps and Neutrosophic Cognitive Maps*, Xiquan, Phoenix, 2003.

43. Vasantha Kandasamy, W.B., Neelakantan, N.R., and Kannan, S.R., Operability Study on Decision Tables in a Chemical Plant using Hierarchical Genetic Fuzzy Control Algorithms, *Vikram Mathematical Journal*, 19 (1999) 48-59.

44. Vasantha Kandasamy, W.B., Neelakantan, N.R., and Kannan, S.R., Replacement of Algebraic Linear Equations by Fuzzy Relation Equations in Chemical Engineering, In *Recent Trends in Mathematical Sciences*, Proc. of Int. Conf. on Recent Advances in Mathematical Sciences held at IIT Kharagpur on Dec. 20-22, 2001, published by Narosa Publishing House, (2001) 161-168.

45. Vasantha Kandasamy, W.B., Neelakantan, N.R., and Ramathilagam, S., Use of Fuzzy Neural Networks to Study the Proper Proportions of Raw Material Mix in Cement Plants, *Varahmihir J. Math. Sci.*, 2 (2002) 231-246.

46. Vasantha Kandasamy, W.B., and Smarandache, F., *Fuzzy Relational Equations and Neutrosophic Relational Equations*, Hexis, Church Rock, 2004.

47. Yager, R.R., On Ordered Weighted Averaging Operators in Multi Criteria Decision Making, *IEEE Trans. Systems, Man and Cybernetics*, 18 (1988) 183-190.

48. Yen, J., Langari, R., and Zadeh, L.A., *Industrial Applications of Fuzzy Logic and Intelligent Systems*, IEEE Press, New York 1995.

49. Zadeh, L.A., A Theory of Approximate Reasoning, *Machine Intelligence*, 9 (1979) 149- 194.

50. Zadeh, L.A., Similarity Relations and Fuzzy Orderings, *Inform. Sci.*, 3 (1971) 177-200.

51. Zimmermann, H.J., *Fuzzy Set Theory and its Applications*, Kluwer, Boston, 1988.

# INDEX

**I**

**M**

**N**

**S**

**T**

# ABOUT THE AUTHORS

**Dr.W.B.Vasantha Kandasamy** is an Associate Professor in the Department of Mathematics, Indian Institute of Technology Madras, Chennai. In the past decade she has guided 12 Ph.D. scholars in the different fields of non-associative algebras, algebraic coding theory, transportation theory, fuzzy groups, and applications of fuzzy theory of the problems faced in chemical industries and cement industries.

She has to her credit 640 research papers. She has guided over 64 M.Sc. and M.Tech. projects. She has worked in collaboration projects with the Indian Space Research Organization and with the Tamil Nadu State AIDS Control Society. This is her 36$^{th}$ book.

On India's 60th Independence Day, Dr.Vasantha was conferred the Kalpana Chawla Award for Courage and Daring Enterprise by the State Government of Tamil Nadu in recognition of her sustained fight for social justice in the Indian Institute of Technology (IIT) Madras and for her contribution to mathematics. (The award, instituted in the memory of Indian-American astronaut Kalpana Chawla who died aboard Space Shuttle Columbia). The award carried a cash prize of five lakh rupees (the highest prize-money for any Indian award) and a gold medal.

She can be contacted at vasanthakandasamy@gmail.com
You can visit her on the web at: http://mat.iitm.ac.in/~wbv

---

**Dr. Florentin Smarandache** is a Professor of Mathematics and Chair of Math & Sciences Department at the University of New Mexico in USA. He published over 75 books and 150 articles and notes in mathematics, physics, philosophy, psychology, rebus, literature.

In mathematics his research is in number theory, non-Euclidean geometry, synthetic geometry, algebraic structures, statistics, neutrosophic logic and set (generalizations of fuzzy logic and set respectively), neutrosophic probability (generalization of classical and imprecise probability). Also, small contributions to nuclear and particle physics, information fusion, neutrosophy (a generalization of dialectics), law of sensations and stimuli, etc. He can be contacted at smarand@unm.edu